\input amstex
\document
\let\co\colon
 \magnification1200

\expandafter\edef\csname :RestoreCatcodes\endcsname{%
   \catcode`\noexpand\noexpand\noexpand \^=\the\catcode`\^%
}
\catcode`\^=7
\expandafter\edef\csname :RestoreCatcodes\endcsname{%
   \csname :RestoreCatcodes\endcsname
   \catcode`\noexpand \_=\the\catcode`\_%
   \catcode`\noexpand :=\the\catcode`:%
   \catcode`\noexpand @=\the\catcode`@%
   \catcode`\noexpand /=\the\catcode`/%
   \catcode`\noexpand &=\the\catcode`&%
   \catcode`\noexpand \^^M=\the\catcode`\^^M%
   \catcode`\noexpand \^^I=\the\catcode`\^^I%
   \let\expandafter\noexpand
       \csname:RestoreCatcodes\endcsname=\noexpand\undefined}
\catcode`\:11  \catcode`\@11
   \let\:wlog\wlog \def\wlog#1{}
   \def\:wrn#1#2{\immediate\write\sixt@@n{--DraTeX warning--
      \ifcase #1
    DraTex.sty already loaded              
\or \string\Draw\space within \string\Draw 
\or Changing definition of \string#2      
\or No intersection points: #2            
\or Improper rotation of axes: #2         
\or (#2) in \string\DSeg\space is a point 
\fi}}
\def\:err#1#2{\errmessage{--DraTeX error-- \ifcase #1
     \string#2\space meaningless in three dimensions 
\or  \string#2\space meaningless in two dimensions   
\or  No \string\MarkLoc(#2)            
\or  \string#2 in three dimensions     
\or  Too many parameters in definition 
\or  \string\MoveFToOval(#2)? 
\fi}}
   \ifx\:Xunits\:undefined \else \:wrn0{} \fi
   \catcode`\ 9  \catcode`\^^M9 \catcode`\^^I9
      \def\:UBorder{0}
\newdimen\:LBorder \newdimen\:RBorder\chardef\:eight=8
\mathchardef\:cccvx=360
\newdimen\:mp    \:mp   0.1\p@
\newdimen\:mmp   \:mmp  0.01\p@
\newdimen\:mmmp  \:mmmp 0.001\p@
\newdimen\:XC    \:XC   90\p@
\newdimen\:CVXXX \:CVXXX180\p@
\newdimen\:CCCVX \:CCCVX\:cccvx\p@ \newdimen\:TeXLoc
\newbox\:box\newif\if:IIID  \newdimen\:Z   \newdimen\:Zunits
\newdimen\:Ex   \newdimen\:Ey  \newdimen\:Ez

\def\:AbsVal#1{ \ifdim#1<\z@-\fi #1 }
\def\:abs#1{\ifdim #1<\z@ #1-#1 \fi}
\def\:AbsDif#1#2#3{  #1#2   \advance#1  -#3
   \ifdim #1<\z@ #1-#1 \fi}
\def\:diff#1#2#3{ #1#2  \advance#1 -#3 }
\def\:average#1#2#3{
   #1#2  \advance#1  #3   \divide#1 \tw@}\def\:Opt#1#2#3#4{
   \def\:temp{
      \ifx      \:next\ifnum \def\:next{#3#1#4#2}
      \else\ifx \:next#1     \def\:next{#3}
      \else                  \def\:next{#3#1#4#2}\fi\fi \:next}
   \futurelet\:next\:temp}\def\Define#1{\:multid#1
   \:Opt(){\:Define#1}0}

\def\:DraCatCodes{\catcode`\ 9   \catcode`\^^M9
   \catcode`\^^I9  \catcode`\&13  \catcode`\~13 }

\def\:Define#1(#2){\begingroup  \:DraCatCodes  \::Define#1(#2)}

\def\::Define#1(#2)#3{\endgroup
   \let\:NextDefine\NextDefine
   \let\NextDefine\relax
   \ifcase#2\relax
      \def#1{#3}\or
      \:TxtPar\def#1(##1){#3}\or
         \:TxtPar\def#1(##1,##2){#3}\or
   \:TxtPar\def#1(##1,##2,##3){#3}\or
   \:TxtPar\def#1(##1,##2,##3,##4){#3}\or
   \:TxtPar\def#1(##1,##2,##3,##4,##5){#3}\or
   \:TxtPar\def#1(##1,##2,##3,##4,##5,##6){#3}\or
   \:TxtPar\def#1(##1,##2,##3,##4,##5,##6,##7){#3}\or
   \:TxtPar\def#1(##1,##2,##3,##4,##5,##6,##7,##8){#3}\or 
      \:TxtPar\def#1(##1,##2,##3,##4,##5,##6,##7,##8,##9){#3}\or
      \:err4{}\fi      \let\:TxtPar\relax  \:NextDefine}

\let\NextDefine\relax\let\:TxtPar\relax
\def\WarningOn{\def\:multid##1{
   \ifx ##1\:undefined \else \:wrn2##1\fi}}
\def\:gobble#1{}
\def\WarningOff{\let\:multid\:gobble}     \WarningOff
\Define\Indirect{\futurelet\:next\:Indirect}

\Define\:Indirect{\:theDoReg
   \ifx \:next<
      \def\:temp{\let\DoReg\:DoReg}
      \def\:next<##1>{\expandafter\:temp\csname :<##1>\endcsname}
   \else
      \def\:next##1<##2>{
         \expandafter\ifx \csname :<##2> \endcsname \relax
               \def\:next{##1}     \fi
\:indrwrn\Define     \:indrwrn\Object
\:indrwrn\Table      \:indrwrn\IntVar     \:indrwrn\DecVar 
         \def\:temp{\let\DoReg\:DoReg##1}
         \expandafter\:temp \csname :<##2> \endcsname}
   \fi      \:next}  \def\:indrwrn#1{    \def\:temp{#1}
   \ifx \:next\:temp \def\:wrn##1##2{\let\:wrn\::wrn} \fi}
\let\::wrn\:wrn 
\Define\:Hline{
   \setbox\:box\hbox{\vrule height0.5\:thickness
                    depth 0.5\:thickness width\:x}
   {\:d\:X \advance\:d \wd\:box
 \advance\:X -\:TeXLoc   \global\:TeXLoc\:d
 \vrule width\:X depth\z@ height\z@
 \raise \:Y \box\:box} }\Define\:Vline{
   \setbox\:box\hbox{\vrule width\:thickness
                 \ifdim \:y>\z@ height\:y  depth\z@
                 \else          height\z@  depth-\:y \fi}
   \advance\:X  -0.5\:thickness
   {\:d\:X \advance\:d \wd\:box
 \advance\:X -\:TeXLoc   \global\:TeXLoc\:d
 \vrule width\:X depth\z@ height\z@
 \raise \:Y \box\:box} }\Define\:MvTo(2){\:X#1\:Xunits \:Y#2\:Yunits}
\Define\:Mv(2){\advance\:X  #1\:Xunits
               \advance\:Y  #2\:Yunits}
\def\:DLn(#1,#2,{\:MvTo(#1,#2) \:LnTo(}\Define\:LnTo(2){
   \:x#1\:Xunits \advance\:x  -\:X
   \:y#2\:Yunits \advance\:y  -\:Y
   \:Ln(\:x\du,\:y\du) }\Define\:Ln(2){
   \:x#1\:Xunits \:y#2\:Yunits
   { \ifdim \:x<\z@
        \advance\:X \:x   \:x-\:x
        \advance\:Y \:y   \:y-\:y
     \fi
     \:yy\:AbsVal\:y
     \:dd\:yy \advance\:dd \:x
     \ifdim \:dd>\:mmmp
        \ifdim \:x>\:yy
           { \ifdim\:X<\:LBorder \global\:LBorder\:X\fi
\advance \:X  \:x
\ifdim \:X>\:RBorder \global\:RBorder\:X\fi\advance \:Y  0.5\:thickness
\ifdim \:Y>\:UBorder \xdef\:UBorder{\the\:Y} \fi
\advance \:Y  \:y
\ifdim \:Y>\:UBorder \xdef\:UBorder{\the\:Y} \fi }
           \let\:Yunitsy\:x  \let\:Xunitsx\:y  \let\:temp\:Hline
           \:dd0.6\:yy    \:divide\:dd\:x  
        \else
           { \advance \:X  -0.5\:thickness
\ifdim \:X<\:LBorder \global\:LBorder\:X \fi
\advance \:X  \:thickness
\advance \:X  \:x
\ifdim \:X>\:RBorder \global\:RBorder\:X\fi
\ifdim \:Y>\:UBorder \xdef\:UBorder{\the\:Y} \fi
\advance \:Y  \:y
\ifdim \:Y>\:UBorder \xdef\:UBorder{\the\:Y} \fi }
           \let\:Yunitsy\:y  \let\:Xunitsx\:x  \let\:temp\:Vline
           \:dd0.6\:x
\ifdim \:x>\:mmp    \:divide\:dd\:yy   \fi
        \fi
        \advance\:dd  0.4\p@
\:ragged\:Cons\:dd\:ragged
        \:HVLn
      \fi }
   \advance\:X  \:x   \advance\:Y  \:y}\Define\:HVLn{
   \:xx\:AbsVal\:Xunitsx  \:divide\:xx\:ragged
   \advance\:xx  0.99\p@   \:K\:InCons\:xx \relax
   \ifnum \:K>\z@
      \divide\:Xunitsx \:K  \advance\:K \@ne
      \divide\:Yunitsy \:K
   \else \:K\@ne  \fi  \:NextLn}

\Define\:NextLn{
   \ifnum\:K=\z@  \let\:NextLn\relax
   \else  { \:temp }  \advance\:K  \m@ne
      \advance\:X  \:x   \advance\:Y  \:y
   \fi   \:NextLn}\newdimen\:ragged

\Define\Ragged(1){    \:ragged#1\p@  \:ragged0.1\:ragged }
\Ragged(7.5)
\Define\PaintUnderCurve(4){{
   \:Z\:Y   \def\:next{\Curve(#1,#2,#3,#4)}
   \MoveToLoc(#1)  \:d\:X   \MoveToLoc(#4)
   \advance\:d -\:X
   \ifdim \:AbsVal\:d<\:mmp  \def\:next{}
   \else
      \def\:CrvLnTo(##1,##2){
         \:x \:X   \:y  \:Y    \:X\:DJ  \:Y\:yyyy
         \:xx\:X   \:dddd \:Y  \:X\:x  \:Y\:y
         { \advance\:Y  \:dddd   \divide\:Y \tw@
           \advance\:Z -\:Y
           \advance\:Y   0.5\:Z
           \:dddd  \:AbsVal \:Z  \:d\z@
           \def\:CrvLnTo{\:LnTo}
           \:yy\:Y  \:dd\:dddd  \:ddd\:dddd
           \::paint  }}
    \fi \:next       }}
\Define\DoCurve(1){ \let\:StopCurve\:SlowCurve
  \def\:CMv(##1){  \:x\:X \:y\:Y   \:MvTo(##1)
      \advance\:x -\:X   \advance\:y -\:Y
      \:xxx \:x    \:yyy\:y}
   \:DoCurve{\Curve(#1)}
   \let\:StopCurve\:FastCurve}

\def\:DoCurve#1(#2)#3{{\XSaveUnits
   \def\:next{#1}    \:MvTo(#2,#2)
   \:x\:AbsVal\:X  \:y\:Y  \:ddd\z@  \:length
   \:Z\:d   \:divide\:Z{1.41421\p@}
   \edef\:tempa{\the\:DoDist}   \global\:DoDist\z@
   \def\:CrvLnTo(##1){ \MarkLoc(1^)    \:CMv(##1)
      { \MarkLoc(2^)   \:ddd\z@    \:length
\:dd  \:DoDist  \global\advance\:DoDist  \:d
\:ddd \:DoDist  \:divide\:ddd\:Z
\DoReg\:InCons\:ddd  \:Z\DoReg\:Z

      \ifdim \:Z>\:dd
         \advance\:Z -\:DoDist
\advance\:dd -\:DoDist
\:divide\:Z\:dd
\advance\:X \:Cons\:Z\:xxx
\advance\:Y \:Cons\:Z\:yyy   \:DoRot 
         \def\:CrvLnTo{\:LnTo}
         \def\:OvalLn{\:Ln}  \XRecallUnits    #3 \fi}}
   \:next    \xdef\:DoDim{\:Cons\:DoDist}
   \global\:DoDist\:tempa     }
   \let\DoDim\:DoDim}

\newdimen\:DoDist\def\:DoRot{ \DSeg\RotateTo(1^,2^) }
\def\DoLine(#1,#2)(#3)#4{
   \MarkLoc($1)  \Move(#1,#2)
   \def\:next{   { \MarkLoc($2)
     \DSeg\RotateTo($1,$2)   \let\:DoRot\relax
     \edef\:RecallRagged{\the\:ragged} \MoveTo(#3,#3)
     \:x\:AbsVal\:X  \:y\:Y  \:ddd\z@  \:length
     \:ragged\:d   \divide\:ragged \tw@
     \DoCurve($1,$1,$2,$2)(#3)
        {\:ragged\:RecallRagged #4}  }
   \let\DoDim\:DoDim}  \:next }
\def\Table#1{\begingroup  \:DraCatCodes   \:multid#1
   \:DefineData#1}

\def\:DefineData#1#2{\endgroup
   \let\:temp~  \def~{\noexpand~}
   \edef#1{\noexpand\:DoPoly\expandafter\noexpand\csname :\string#1\endcsname}
   \expandafter\edef  \csname :\string#1\endcsname
       ##1{\noexpand\ifcase##1(#2)\noexpand\fi}
   \let~\:temp   \:DoNextPoly   \:DoNextPoly}

\def\:OR{\let\:or\or}  \:OR \catcode`\&13  \def&{)\noexpand\:or(}

\def\TableData{\begingroup  \:DraCatCodes  \:TableData}

\def\:TableData#1#2#3{\endgroup   \Table\:temp{#3}
   \:K\z@   \:J\z@    \def\:tempa(##1){\advance\:J \@ne }
   \:temp(0,999){\:tempa}    \let\:tempa&      \def\:temp{\def#1}
   \def&##1&{
      \ifnum  \:K<\:J
         \advance\:K \@ne
         \ifnum  \:K=\@ne   \def#1{#2(##1)}
         \else
            \:IIIexpandafter\:temp\expandafter{
               #1 & #2(##1) }
         \fi
      \else  \let\:next\relax \fi
      \:next}
   \let\:next&     &#3&&   \let&\:tempa }

\catcode`\&4  \def\:DoPoly#1(#2,#3)#4{
   \expandafter\let \csname :Back\the\:level\endcsname\:or
\expandafter\edef\csname :DoVars\the\:level\endcsname{
   \:DoB\the\:DoB}
\advance\:level  \@ne
   \:DoB#3  \advance\:DoB -#2
   \def\:PolyOr(##1){
      \ifnum  \:DoB=\z@  \:OR
      \else   #4(##1)   \advance\:DoB \m@ne    \fi}
   \:OR
   \def\:temp{\let\:or\:PolyOr #4}
   \:IIIexpandafter\:temp#1{#2}
   \advance\:level  \m@ne
\csname :DoVars\the\:level\endcsname
\def\:temp{\let\:or}
\expandafter\:temp\csname :Back\the\:level\endcsname }
\Define\PaintRect(2){\def\:next{{   \MarkLoc(^)
   \MoveToLoc(^) \Move( #1,0) \MarkLoc(^1) \Move(0,#2)
   \MarkLoc(^2)  \Move(-#1,0) \MarkLoc(^3)
   \PaintQuad(^,^1,^2,^3)        }}\:next}

\Define\PaintRectAt(4){\def\:next{{   \MoveTo(#1,#2)
   \:K#3   \advance\:K -#1
   \:DoB#4 \advance\:DoB -#2  \PaintRect(\:K,\:DoB)}}
   \:next}\def\::paint{
   \ifdim \:d<\:ragged      \advance\:xx -\:X
      \:yyy\:Y \:xxx\:dddd
      \advance\:Y \:yy  \divide\:Y \tw@
      \:average\:dddd\:dd\:ddd
      \def\:next{\:brush(\:xx,\z@)\:Y\:yyy\:dddd\:xxx}
   \else  \divide\:d \tw@
      \:average\:x\:X\:xx
      \:average\:y\:Y\:yy
      \:average\:dddd\:dd\:ddd
   \fi   \:next}\Define\:paint{{   \:AbsDif\:d\:xx\:x
   \ifdim \:d<\:mmp  \let\::paint\relax
   \else
      \ifdim  \:y >\:yyy   \:dd\:y   \:y \:yyy   \:yyy \:dd  \fi
      \ifdim  \:yy>\:yyyy  \:dd\:yy  \:yy\:yyyy  \:yyyy\:dd  \fi
      \:AbsDif\:dd\:yyyy\:yyy      \:AbsDif\:ddd\:yy\:y
      \ifdim \:dd<\:ddd  \:dd\:ddd  \fi
      \ifdim \:d >\:dd   \:d\:dd  \fi
      \advance\:d \:d      \:diff\:dd\:y\:yyy
      \:diff\:ddd\:yy\:yyyy
      \advance\:y    -0.5\:dd    \:abs\:dd
      \advance\:yy   -0.5\:ddd   \:abs\:ddd
      \:X\:x  \:Y\:y
      \:average\:dddd\:dd\:ddd
      \def\:next{ \:lpaint \:rpaint }
   \fi
   \::paint }}

\Define\:lpaint{ { \:xx\:x  \:yy \:y  \:ddd\:dddd  \::paint} }
\Define\:rpaint{ { \:X \:x  \:Y  \:y  \:dd \:dddd  \::paint} }
\Define\PaintQuad(4){\def\:next{{\Units(1pt,1pt)
   \MoveToLoc(#1)  \:x  \:X  \:y  \:Y
   \MoveToLoc(#2)  \:xx \:X  \:yy \:Y
   \MoveToLoc(#3)  \:xxx\:X  \:yyy\:Y
   \MoveToLoc(#4)  \:xxxx  \:X  \:yyyy \:Y
   
   \:paintQuad  }}\:next}

\def\:paintQuad{{
   \:SetVal\:a\:x\:y\:xx\:yy\:xxxx\:yyyy
\:SetVal\:b\:xx\:yy\:xxx\:yyy\:x\:y
\:SetVal\:c\:xxx\:yyy\:xx\:yy\:xxxx\:yyyy
\:SetVal\:cc\:xxxx\:yyyy\:xxx\:yyy\:x\:y
\def\:A{\:a} \def\:B{\:b} \def\:C{\:c} \def\:D{\:cc}
\:sort\:B\:A
\:sort\:C\:B  \:sort\:B\:A
\:sort\:D\:C  \:sort\:C\:B  \:sort\:B\:A
\let\:temp\relax
\:IsTriang\:A\:B
\:IsTriang\:B\:C
\:IsTriang\:C\:D
\:temp  
   \:Quad\:A\:temp>   \:xxxx\:xx \:yyyy\:yy \:Z\:d
\:Quad\:D\:next<   
   \:PrePaint(\:xx,\:yy,\:d,\:xxxx,\:yyyy,\:Z)
   \:temp  \:next }}
\Define\:PrePaint(6){
   \:x#1 \:y#2 \:yyy#3 \:xx#4 \:yy#5 \:yyyy#6
   \:paint }\def\:Quad#1#2#3{
         \:GetVal#1\:x\:y0
   \:GetVal#1\:xx\:yy1
   \:GetVal#1\:xxx\:yyy2
   \ifdim \:xx#3\:xxx
      \:ddd\:xx   \:xx\:xxx   \:xxx\:ddd
      \:ddd\:yy   \:yy\:yyy   \:yyy\:ddd
   \fi
          \def#2{}
   \:diff\:dd\:xxx\:xx
   \ifdim \ifdim\:AbsVal\:dd<\:mp  \:yy=\:yyy   \else  \z@>\z@  \fi
      
      \:d\:yy
   \else
      \ifdim \:AbsVal\:dd>\:mp
         \:diff\:dd\:xxx\:x   \:diff\:ddd\:yyy\:y
         \:divide\:ddd\:dd    \:diff\:dd\:xx\:xxx
         \:ddd\:Cons\:ddd\:dd   \advance\:ddd \:yyy
         \:d\:ddd
      \else \:d\:yyy \fi
      \edef#2{  \noexpand\:PrePaint
         (\the\:x ,\the\:y ,\the\:y,
         \the\:xx,\the\:yy,\the\:d)    }
   \fi}\def\:SetVal#1#2#3#4#5#6#7{
   \edef#1{(\the#2,\the#3,\the#4,\the#5,\the#6,\the#7)}}

\def\:sort#1#2{
   \ifdim \:IIIexpandafter\:field#1 <
          \:IIIexpandafter\:field#2
      \let\:temp#1  \let#1#2  \let#2\:temp
   \fi  }

\Define\:field(6){#1}

\def\:GetVal#1#2#3{
   \:IIIexpandafter\::GetVal #1#2#3}

\def\::GetVal(#1,#2,#3,#4,#5,#6)#7#8#9{
      \ifcase #9 #7#1   #8#2\or #7#3   #8#4\or #7#5   #8#6 \fi}
\def\:IsTriang#1#2{
   \ifdim \:IIIexpandafter\:field#1 =
          \:IIIexpandafter\:field#2
      \ifdim \:IIIexpandafter\:fieldB#1 =
             \:IIIexpandafter\:fieldB#2
         \def\:temp{ \:FixTria }
   \fi \fi  }

\def\:FixTria{
   \edef\:temp{\:IIIexpandafter\:FrsII\:B}
   \ifdim \:IIIexpandafter\:field\:A =
          \:IIIexpandafter\:field\:B
      \ifdim \:IIIexpandafter\:fieldB\:A =
             \:IIIexpandafter\:fieldB\:B
          \edef\:temp{\:IIIexpandafter\:FrsII\:C}
   \fi\fi
   \edef\:A{\:IIIexpandafter\:FrsII\:A}
   \edef\:D{\:IIIexpandafter\:FrsII\:D}
   \edef\:temp{
      \def\noexpand\:a{(\:A,\:temp,\:D)}
      \def\noexpand\:b{(\:temp,\:A,\:D)}
      \def\noexpand\:c{\noexpand\:b}
      \def\noexpand\:cc{(\:D,\:A,\:temp)}}
   \:temp
   \def\:A{\:a}  \def\:B{\:b}  \def\:C{\:c}  \def\:D{\:cc}  }

\Define\:fieldB(6){#2}
\Define\:FrsII(6){#1,#2}
\def\:IIIexpandafter{\expandafter\expandafter\expandafter}

\Define\DrawRect(2){ \Line( #1,0) \Line(0, #2)
                     \Line(-#1,0) \Line(0,-#2)}

\Define\DrawRectAt(4){{
   \MoveTo(#1,#2) \LineTo(#3,#2) \LineTo(#3,#4)
   \LineTo(#1,#4) \LineTo(#1,#2)}}
\def\du#1{ \ifx#1\:Xunits   \else\ifx#1\:Yunits
      \else\ifx#1\:Zunits   \else #1
      \fi \fi \fi}
\let\:svdu=\du
\Define\XSaveUnits{
   \expandafter\edef\csname XRecallUnits\the\:level\endcsname{
      \:StoreUnits}
    \advance\:level  \@ne}

\Define\XRecallUnits{
   \advance\:level \m@ne
   \csname XRecallUnits\the\:level \endcsname}

\Define\SaveUnits{
     }

\Define\:StoreUnits{       \:Xunits \the\:Xunits
   \:Yunits \the\:Yunits \:Zunits \the\:Zunits
   \:Xunitsx\the\:Xunitsx \:Xunitsy\the\:Xunitsy
\:Yunitsx\the\:Yunitsx \:Yunitsy\the\:Yunitsy }
\Define\:SearchDir{
   \ifdim \:x<\z@  \:x-\:x  \:y-\:y
   \edef\:tempA{\advance\:ddd  \ifdim \:y<\z@ - \fi\:CVXXX}
\else         \def\:tempA{} \fi
\ifdim \:y<\z@
   \edef\:tempA{\:ddd-\:ddd \advance\:ddd  \:CCCVX \:tempA}
   \:y-\:y
\fi
\ifdim \:y>\:x
   \:ddd\:y \:y\:x \:x\:ddd
   \edef\:tempA{\advance\:ddd  -\:XC \:ddd-\:ddd \:tempA}
\fi
   \:divide\:y\:x  \:d57.29578\:y
   \:ddd\:d       \:sqr\:y  \:K  \@ne
   \Do(1,30){\advance\:K \tw@
      \:d-\:Cons\:y\:d  \:dd\:d  \divide\:dd \:K
      \advance\:ddd \:dd  }
   \advance\:ddd -0.49182\:dd
   \:tempA }\Define\Curve(4){\def\:next{{  \XSaveUnits \Units(1pt,1pt)
   \MoveToLoc(#1) \:DI \:X  \:DK \:Y
   \MoveToLoc(#2) \:ddd \:X  \:Ez \:Y
   \MoveToLoc(#3) \:dd\:X  \:Ey\:Y
   \MoveToLoc(#4) \:DJ\:X  \:yyyy\:Y   \:Curve    }}\:next}
\Define\:FastCurve{   \:AbsDif\:d\:DI\:ddd  \:AbsDif\:dddd\:DK\:Ez
   \advance\:d \:dddd
   \ifnum \:d<\:ragged
      \:AbsDif\:d\:DJ\:dd \:AbsDif\:dddd\:yyyy\:Ey
      \advance\:d \:dddd                     \fi}
\let\:StopCurve\:FastCurve
\Define\:SlowCurve{
   \:AbsDif\:d\:DI\:DJ  \:AbsDif\:dddd\:DK\:yyyy
   \advance\:d \:dddd                     }
\Define\:Curve{   \:StopCurve
   \ifnum \:d<\:ragged   \:X\:DI \:Y\:DK
                    \:CrvLnTo(\:Cons\:DJ,\:Cons\:yyyy)
   \def\:SubCurves{}  \fi
   \:SubCurves}

\def\:CrvLnTo{\:LnTo}\Define\:SubCurves{
   \:average\:yy\:DI\:ddd     \:average\:Ex\:DK\:Ez
   \:average\:ddd\:ddd\:dd    \:average\:Ez\:Ez\:Ey
   \:average\:dd\:dd\:DJ  \:average\:Ey\:Ey\:yyyy
   \:average\:Zunits\:yy\:ddd    \:average\:Vdirection\:Ex\:Ez
   \:average\:ddd\:ddd\:dd    \:average\:Ez\:Ez\:Ey
   \:average\:DL\:Zunits\:ddd  \:average\:xxxx\:Vdirection\:Ez
   { \:ddd \:yy  \:Ez \:Ex   \:dd\:Zunits
     \:Ey\:Vdirection \:DJ\:DL \:yyyy\:xxxx \:Curve }
   \:DI \:DL \:DK \:xxxx               \:Curve   }
\def\MoveToCurve[#1]{
   \Define\:BiSect(3){\MoveToLoc(##1)
      \CSeg[#1]\Move(##1,##2)
   \MarkLoc(##3) }\:MvToCrv}

\Define\:MvToCrv(4){  \:BiSect(#1,#2,:a) \:BiSect(#2,#3,:b)
   \:BiSect(#3,#4,:c) \:BiSect(:a,:b,:A) \:BiSect(:b,:c,:B)
   \:BiSect(:A,:B,:Q)}
\Define\:OvalDir(3){   \:CosSin{#3\p@}
   \:Zunits\p@   \:d#2\:Zunits   \:x\:Cons\:d\:x
   \:Zunits\p@   \:d#1\:Zunits   \:y\:Cons\:d\:y
   \:SearchDir   }\Define\DrawOval   (2){ \DrawOvalArc(#1,#2)(0,\:cccvx) }
\Define\PaintOval  (2){ \PaintOvalArc(#1,#2)(0,\:cccvx) }
\Define\DrawCircle (1){ \DrawOval(#1,#1) }
\Define\PaintCircle(1){ \PaintOval(#1,#1) }
\def\DrawOvalArc(#1,#2)(#3,#4){{
       \:xxxx#4\p@  \advance\:xxxx -#3\p@
       \ifdim \:xxxx=\z@ \else
   \let\:SinOne\:SinB  \:OvalDir(#1,#2,#3)  \:DJ\:ddd
\:OvalDir(#1,#2,#4)  \:diff\:DI\:ddd\:DJ
\ifdim\:DI<\z@ \advance\:DI  \:CCCVX \fi
\ifdim \:xxxx<\:CCCVX \else \:DI\:CCCVX \fi
\:InitOval(#1,#2)  \:CosSin\:DJ 
   \:xxxx\:x  \:yyyy\:y   \:xx\:X  \:yy\:Y
   \advance\:X \:Xx\:x  \advance\:X \:Yx\:y
   \advance\:Y \:Xy\:x  \advance\:Y \:Yy\:y
   \let\:Xunits\empty  \let\:Yunits\empty
   \Do(1,\:InCons\:DI){
      \:dd\:X  \:ddd\:Y  \:X\:xx  \:Y\:yy
      \:AdvOv\:xxx\:yyy\:xxxx\:yyyy
      \:X\:dd  \:Y\:ddd
      \advance\:xxx -\:X   \advance\:yyy -\:Y
      \:d\:AbsVal\:xxx
      \advance\:d  \:AbsVal\:yyy
      \ifdim \:d>\:ragged
         \:OvalLn(\:xxx,\:yyy)  \fi  }
   \:OvalDir(#1,#2,#4)    \:CosSin\:ddd
   \advance\:xx \:Xx\:x  \advance\:xx \:Yx\:y
   \advance\:yy \:Xy\:x  \advance\:yy \:Yy\:y
   \advance\:xx -\:X     \advance\:yy -\:Y
   \:OvalLn(\:xx,\:yy) \fi }}

\def\:OvalLn{\:Ln}\def\DoOvalArc(#1)(#2){   \:xx\:X  \:yy\:Y
   \def\:CMv(##1){  \:Mv(\:xxx,\:yyy)
      \:x\:xxx  \:y\:yyy}
   \:DoCurve{             \:X\:xx  \:Y\:yy
       \def\:DoRot{  \let\:Xunits\:XunitsReg
                     \let\:Yunits\:YunitsReg
                     \DSeg\RotateTo(1^,2^)     }
       \let\::OvalLn\:CrvLnTo
\Define\:OvalLn(2){ \:dd\:AbsVal####1
   \advance\:dd \:AbsVal####2 \:divide\:dd\:ragged
   \:J\:InCons\:dd  \advance\:J  \@ne
   \divide####1  \:J  \divide####2  \:J
   \Do(1,\:J){\::OvalLn(####1,####2)}}
       \DrawOvalArc(#1)(#2)}}
\def\NextTable{\begingroup  \:DraCatCodes \:NextTable}

\def\:NextTable#1{\endgroup
  \def\:DoNextPoly{#1\NextTable{}}}
\NextTable{}\def\:AdvOv#1#2#3#4{
   \:d\:CosOne#3 \advance\:d -\:SinOne#4
   #4\:CosOne#4    \advance#4     \:SinOne#3   #3\:d
   \divide#3 \:eight  \divide#4 \:eight
   #1\:X   #2\:Y
   \:d\:Xx#3  \advance\:d \:Yx#4  \advance#1  \:d
   \:d\:Xy#3  \advance\:d \:Yy#4  \advance#2  \:d  }

\def\:CosOne{7.99878}   \def\:SinB{0.13962}
\def\PaintOvalArc(#1,#2)(#3,#4){{ \ifdim #3\p@=#4\p@
                      \let\:next\relax  \else
   \:d\:AbsVal{#1\:Xunits} \advance\:d \:AbsVal{#2\:Yunits}
   \ifdim \:d<3\:ragged   \divide \:d \tw@   \PenSize(\:d)
\:Mv(-0.5\:d\du,0)  \:Ln(\:d\du,0)
   \else    \:InitOval(#1,#2)
      \MarkLoc(o$)   \RotateTo(#3) \MoveFToOval(#1,#2)
      \:Ex\:X  \:Ey\:Y   \edef\:FirstOvDir{\:Cons\:ddd\p@}
      \MoveToLoc(o$) \RotateTo(#4) \MoveFToOval(#1,#2)
      \:Ez \:X  \:Vdirection \:Y   \edef\:LastOvDir{\:Cons\:ddd\p@}
      \MoveToLoc(o$)
      \if:rotated
      \:Zunits\p@     \:Zunits#1\:Zunits
      \:xx\:Cons\:Zunits\:Xunitsx
      \:Zunits\p@     \:Zunits#2\:Zunits
      \:yy\:Cons\:Zunits\:Yunitsx
      \:x\:xx  \:y\:yy  \:ddd\z@  \:length
      \:ddd\:d  \:divide\:xx\:ddd   \:divide\:yy\:ddd
\else \:xx\p@  \:yy\z@   \fi
      \:AbsDif\:d{#3\p@}{#4\p@}
      \ifdim \:d>359\p@ \:Ez-\:Xx\:xx  \advance\:Ez -\:Yx\:yy
\:Vdirection-\:Xy\:xx  \advance\:Vdirection -\:Yy\:yy
\advance\:Ez \:X  \advance\:Vdirection \:Y
\:setpaint\:PaintOvOv<> 
      \else
         \:x\:xx \:y\:yy \:SearchDir
\:xxx\:FirstOvDir  \advance\:xxx -\:ddd
\ifdim \:xxx<\z@    \advance\:xxx \:CCCVX  \fi
\:yyy\:LastOvDir   \advance\:yyy -\:ddd
\ifdim \:yyy<\z@    \advance\:yyy \:CCCVX  \fi
\:J\z@
\ifdim \:xxx<\:yyy  \ifdim        \:yyy<\:CVXXX
            \:Pntovln\:FirstOvDir\:FirstOvDir\:LastOvDir
\else \ifdim  \:xxx>\:CVXXX
            \:Pntovln\:FirstOvDir\:FirstOvDir\:LastOvDir
\else \:yyy-\:yyy  \advance\:yyy \:CCCVX
      \ifdim  \:xxx<\:yyy
            \:setpaint\:PntLeftOvOv><  \:PntMovln\:FirstOvDir\:xx
      \else \:FxLx  \:setpaint\:PntLeftOvOv><  \:Pntmovln\:LastOvDir
\fi  \fi  \fi
\else               \ifdim        \:xxx<\:CVXXX
            { \:setpaint\:PaintOvOv<> }  \:FxLx  \:setpaint\:PntLeftOvOv><
\:Usrch  \:PaintMidOvLn\:LastOvDir\:FirstOvDir
\else \ifdim  \:yyy>\:CVXXX
            {  \:FxLx \:setpaint\:PaintOvOv<>  }  \:setpaint\:PntLeftOvOv><
 \:Dsrch  \:PaintMidOvLn\:LastOvDir\:FirstOvDir
\else \:xxx-\:xxx  \advance\:xxx \:CCCVX
      \ifdim  \:yyy<\:xxx
            \:setpaint\:PaintOvOv<>    \:PntMovln\:FirstOvDir\:xx
      \else
            \:FxLx  \:setpaint\:PaintOvOv<>   \:Pntmovln\:LastOvDir
\fi  \fi  \fi  \fi
   \fi \fi \fi}}\def\:setpaint#1#2#3{{\aftergroup#1
\:d\:Xx\:xx   \advance\:d \:Yx\:yy
\ifdim \:d<\z@  \aftergroup#3
\else           \aftergroup#2  \fi}}
\def\:FxLx{\:d\:Ex  \:Ex\:Ez  \:Ez\:d}

\def\:Pntovln#1{
   \let\:SinOne\:SinB     \:CosSin#1
   \:xxxx\:Ex  \:yyyy\:Ey  \:Z\z@
   \:PaintOvLn}

\def\:PntMovln{
   \:Dsrch  \:FxLx   \:xx\:ddd  \:PaintOvLn}

\def\:Pntmovln{
   \:Usrch  \:FxLx  \:xx\:ddd   \:PaintOvLn\:xx}
\def\:PaintMidOvLn#1#2{
   \:FxLx   \:xx\:ddd  \:PaintOvLn#1#2
   \:xxx\:Ez  \:yyy\:Vdirection
   \:ddd\:Xy\:x  \advance\:ddd \:Yy\:y  \advance\:ddd \:Y
   \:PaintSlice}

\def\:PntLeftOvOv{
   \:xx-\:xx  \:yy-\:yy  \:PaintOvOv}

\def\:Dsrch{      \def\:SinOne{-\:SinB}
   \:xx\:x  \:yy\:y   \:SearchDir
   \:x\:xx  \:y\:yy
   \:d\:Ey  \:Ey\:Vdirection  \:Vdirection\:d }

\def\:Usrch{
   \let\:SinOne\:SinB
   \:x\:xx  \:y\:yy   \:SearchDir
   \:x\:xx  \:y\:yy}\def\:PaintOvLn#1#2{
   \:diff\:dd\:Ey\:Vdirection  \:diff\:ddd\:Ex\:Ez
   \ifdim \:AbsVal\:ddd>\:mp
      \:divide\:dd\:ddd      \:d#2
      \advance\:d -#1
      \ifnum \:d<\z@  \advance\:d \:CCCVX  \fi
      \:DoB\:InCons\:d  \let\:next\:PntDo  \:next
   \fi}

\def\:PntDo{
   \ifnum\:DoB=\z@ \let\:next\relax
   \else
      \::AdvOv\:x\:y
      \ifdim \:d>\:ragged
         \:ddd\:xxx  \advance\:ddd -\:Ez
         \:ddd\:Cons\:dd\:ddd
         \advance\:ddd \:Vdirection  \:PaintSlice
      \fi
      \advance\:DoB \m@ne
   \fi  \:next}\Define\:PaintSlice{   \:AbsDif\:dddd\:yyy\:ddd
   \advance\:yyy \:ddd  \divide\:yyy \tw@
   { \advance\:dddd \:Z  \divide\:dddd \tw@
      \:X\:xxxx  \:Y\:yyy
      \:xx\:xxx  \advance\:xxx -\:xxxx
      \:d\z@
      \:yy\:Y  \:dd\:dddd  \:ddd\:dddd
      \::paint     }
   \:xxxx\:xxx  \:yyyy\:yyy  \:Z\:dddd    \:J\z@}
\def\::AdvOv#1#2{  \:AdvOv\:xxx\:yyy#1#2
   \:AbsDif\:d\:xxxx\:xxx       \advance\:J \@ne
   \ifnum \:J=\sixt@@n   \multiply\:d \@cclvi
   \fi }\def\:PaintOvOv#1{   \def\:hdir{#1}
   \:xxx\:Xx\:xx  \advance\:xxx \:Yx\:yy
   \:yyy\:Xy\:xx  \advance\:yyy \:Yy\:yy
   \advance\:xxx \:X   \advance\:yyy \:Y
   \:Z\z@  \:xxxx\:xxx  \:yyyy\:yyy
   \:x\:xx   \:y\:yy   \:DoB\z@   \:J\z@
   \let\:next\:scanOvOv  \:next }
\Define\:scanOvOv{  \advance\:DoB \@ne
   \advance\:J \@ne
   \:AbsDif\:d\:xxx\:Ez
   \ifdim \ifdim      \:xxx\:hdir\:Ez    -
          \else\ifdim \thr@@\:d<\:ragged -
          \else\ifnum \:DoB>358            -
          \fi \fi \fi                     \p@<\z@
      \:ddd\:yyyy  \advance\:ddd -0.5\:dddd
      \:yyy\:yyyy   \advance\:yyy   0.5\:dddd
      \:xxx\:Ez \:PaintSlice
      \let\:next\relax
   \else
     \:d\:xx \advance\:d -\:x
     \:xxx\:yy \advance\:xxx -\:y
     \ifdim -\:Xx\:d\:hdir\:Yx\:xxx
        \let\:SinOne\:SinB
        \::AdvOv\:xx\:yy
        \ifdim  \ifdim       \:d >\:ragged -
                \else \ifnum \:J>\sixt@@n -
                \fi\fi                      \p@<\z@
           \:J\z@  \def\:SinOne{-\:SinB}
           \:AdvOv\:dd\:ddd\:x\:y  \:PaintSlice  \fi
     \else
        \def\:SinOne{-\:SinB}
        \::AdvOv\:x\:y
        \ifdim  \ifdim       \:d >\:ragged -
                \else \ifnum \:J>\sixt@@n -
                \fi\fi                      \p@<\z@
           \:J\z@   \let\:SinOne\:SinB
           \:AdvOv\:dd\:ddd\:xx\:yy  \:PaintSlice  \fi
   \fi \fi
   \:next}\Define\SetBrush{\:Opt[]\:SetBrush{}}

\def\:SetBrush[#1](#2,#3)#4{    \def\:temp{#4}
   \ifx \:temp\empty
      \def\:brush{   \let\:Xunits\empty \let\:Yunits\empty
                     \:thickness\:dddd  \:Ln  }
   \else       \def\:BruShape{#4}
      \:dd#2\:Xunits      \:ddd#3\:Yunits
      \edef\:Grd{ \:dd\the\:dd  \:ddd\the\:ddd }
      \MarkLoc($$)  \def\:temp{#1}
\ifx \:temp\empty  \:X\z@  \:Y\z@
\else  \MoveTo(#1) \fi
\edef\:BrOrg{ \:x\the\:X  \:y\the\:Y }
\MoveToLoc($$)
      \def\:brush(##1,##2){ \::brush }  \fi  }

\SetBrush(,){}\def\::brush{{  \SetBrush(,){}
\let\:Xunits\:XunitsReg   \let\:Yunits\:YunitsReg
\advance\:Y -0.5\:dddd   \:yy\:Y
\advance\:yy \:dddd
\:BrOrg  \:Grd  \advance\:xx  \:X
\ifdim \:xx<\:X  \:d\:X \:X\:xx \:xx\:d  \fi
   \:GridPt\:X\:x\:dd
   \:GridPt\:Y\:y\:ddd            \:x\:X
   \:DoBrush       }}\Define\:DoBrush{
   \ifdim       \:Y>\:yy  \let\:DoBrush\relax
   \else \ifdim \:X>\:xx \advance\:Y \:ddd     \:X\:x
   \else { \:BruShape }  \advance\:X \:dd
   \fi  \fi  \:DoBrush  }\def\:GridPt#1#2#3{   \:xxxx#1
   \advance#1 -#2  \:divide#1#3
   #1\:InCons#1#3  \advance#1 #2
   \ifdim #1=\:xxxx
   \else  \ifdim \:xxxx>#2 \advance#1 #3 \fi \fi  }
\newcount\:IntId    \edef\:IntCount{0\space}
\newcount\:DecId    \edef\:DecCount{0\space}

\Define\:NewCount{\alloc@ 0\count \countdef \insc@unt }
\Define\:NewDimen{\alloc@ 1\dimen \dimendef \insc@unt }

\def\:NewVar#1#2#3#4#5{ \:multid#1
   \def\:temp{    \csname \string#4\the#4\endcsname\z@
      \edef#1{\noexpand#3  \csname \string#4\the#4\endcsname}}
   \def\:next{  \xdef#5{\the#4\space}   #2#1 \def\:next{\global\let}
     \expandafter\:next \csname \string#4\the#4\endcsname#1  \:temp }
   \advance#4  \@ne
   \ifnum #4 > #5 \expandafter\:next \else \expandafter\:temp \fi }

\def\IntVar#1{\:NewVar#1\:NewCount\:IntOp\:IntId\:IntCount}
\def\DecVar#1{\:NewVar#1\:NewDimen\:DecOp\:DecId\:DecCount}

\DecVar\Q  \DecVar\R   \DecVar\T
\IntVar\I  \IntVar\J   \IntVar\K
\def\WriteVal#1{\immediate\write\sixt@@n{...\string#1=#1;}}

\newdimen\:X   \newdimen\:Y
\newdimen\:x   \newdimen\:y   \newdimen\:d
\newdimen\:xx  \newdimen\:yy  \newdimen\:dd
\newdimen\:xxx \newdimen\:yyy \newdimen\:ddd
\newdimen\:xxxx\newdimen\:yyyy\newdimen\:dddd
\newcount\:J  \newcount\:K
\newdimen\:DI   \newdimen\:DJ
\newdimen\:DK   \newdimen\:DL   \newtoks\:t
\def\:IntFromPt#1#2{
   \:d#2\relax
   \advance\:d  \ifdim\:d<-0.5\p@-\fi  0.5\p@
   #1\:d    \divide#1  65536\relax}
\def\:temp{\catcode`\p12  \catcode`\t12}
\def\:Cons{\catcode`\p11  \catcode`\t11}
\:temp  \def\:Frac#1pt{#1}
        \def\:rnd#1.#2pt{#1}  \:Cons

\def\:Cons#1{\expandafter\:Frac\the#1}
\def\:sqr#1{#1\expandafter\:Frac\the#1#1}
\def\:InCons#1{\expandafter\:rnd\the#1}\def\:Val#1{#1;}
\let\Val\:Val\def\:IntOp#1#2{\csname :Op#2\endcsname#1}
\let\:SvIntOp\:IntOp
\def\:PreIntOp{\let\:IntOp\empty
   \let\Val\empty}
\def\:PostIntOp{\let\:IntOp\:SvIntOp
   \let\Val\:Val}

\expandafter\def\csname :Op;\endcsname#1{ \the#1}
\expandafter\def\csname :Op=\endcsname#1#2;{
   \:PreIntOp#1#2\:PostIntOp}
\expandafter\def\csname :Op+\endcsname#1#2;{
   \:PreIntOp\advance #1  #2\:PostIntOp}
\expandafter\def\csname :Op-\endcsname#1#2;{
   \:PreIntOp\advance #1  -#2\:PostIntOp}
\expandafter\def\csname :Op/\endcsname#1#2;{
   \:PreIntOp\divide#1   #2\:PostIntOp}
\expandafter\def\csname :Op*\endcsname#1#2;{
   \:PreIntOp\multiply#1   #2\:PostIntOp}
\def\:DecOp#1#2{ \csname :xOp#2\endcsname#1}
\let\:SvDecOp\:DecOp
\def\:PreDecOp{\let\:IntOp\the \def\:DecOp{\:Cons}
   \let\Val\empty   \let\:du\empty}
\def\:PostDecOp{\let\:IntOp\:SvIntOp \let\Val\:Val
   \let\:DecOp\:SvDecOp  \let\:du\::du  }

\def\::du#1{\p@
   \ifx#1\p@ \let\:temp\relax
   \else     \def\:temp{\du{#1}}
   \fi\:temp}                    \:PostDecOp

\expandafter\def\csname :xOp;\endcsname#1{ \:Cons#1}
\expandafter\def\csname :Op[\endcsname#1#2];{
   \:PreDecOp \:dd#2\p@  \:IntFromPt#1\:dd
                                 \:PostDecOp  }
\expandafter\def\csname :xOp=\endcsname#1#2;{
   \:PreDecOp#1#2\p@\:PostDecOp               }
\expandafter\def\csname :xOp(\endcsname#1#2){
   \:PreDecOp#1#2\p@\:PostDecOp               }
\expandafter\def\csname :xOp+\endcsname#1#2;{
   \:PreDecOp\advance #1  #2\p@\:PostDecOp  }
\expandafter\def\csname :xOp-\endcsname#1#2;{
   \:PreDecOp\advance #1  -#2\p@\:PostDecOp }
\expandafter\def\csname :xOp*\endcsname#1#2;{
   \:PreDecOp#1 #2#1\:PostDecOp              }
\expandafter\def\csname :xOp/\endcsname#1#2;{
   \:PreDecOp  \:divide#1{#2\p@}  \:PostDecOp }
\let\IF\ifnum    \let\THEN\relax
\let\ELSE\else   \let\FI\fi

\def\EqText(#1,#2){
   \z@=\z@ \fi  \def\:temp{#1}
                \def\:next{#2}    \ifx \:temp\:next }

\def\:IfInt#1(#2,#3){ \z@=\z@ \fi
   \:IntOp\:K=#2;  \:IntOp\:J=#3; \ifnum  \:K#1\:J }

\def\:IfDim#1(#2,#3){ \z@=\z@ \fi
   \:DecOp\:d=#2;  \:DecOp\:dd=#3; \ifdim  \:d#1\:dd }

\def\EqInt{ \:IfInt= }  \def\LtInt{ \:IfInt< }
\def\GtInt{ \:IfInt> }

\def\EqDec{ \:IfDim= }  \def\LtDec{ \:IfDim< }
\def\GtDec{ \:IfDim> } \def\Do(#1,#2)#3{
   \expandafter\let
   \csname :Back\the\:level\endcsname\:Do
\expandafter\edef\csname :DoVars\the\:level\endcsname{
   \DoReg\the\DoReg \:DoB\the\:DoB}
\advance\:level  \@ne
   \DoReg#1  \:DoB#2  \relax
   \ifnum \DoReg<\:DoB
      \def\:Do{\ifnum \DoReg>\:DoB
                  \let\:Do\relax
               \else  #3\advance\DoReg  \@ne \fi
               \:Do}
   \else
      \def\:Do{\ifnum \DoReg<\:DoB
                  \let\:Do\relax
               \else  #3\advance\DoReg \m@ne  \fi
               \:Do}
   \fi  \def\:nextdo{ \:Do \advance\:level  \m@ne
\csname :DoVars\the\:level\endcsname
\def\:temp{\let\:Do}
\expandafter\:temp\csname
   :Back\the\:level\endcsname  } \:nextdo}

\let\:Do\relax

\newcount\DoReg   \let\:DoReg\DoReg       \newcount\:DoB
 \newcount\:level\def\::divide#1{   \:DI\:DK   \:dddd\:DL
   \advance\:DI -\:Cons\:dddd#1
   \:IntFromPt\:J\:dddd  \advance\:dddd -\:J\p@
   \multiply\:J  \@M   \:IntFromPt\:K{\@M\:dddd}
   \advance\:J \:K     \:dddd\@M\p@
   \divide\:dddd \:J   \advance#1 \:Cons\:DI\:dddd  }

\def\:divide#1#2{   \:DK#1   \:DL#2   #1\z@
   \::divide#1  \::divide#1  \::divide#1
   \::divide#1  \::divide#1  }
\def\:Sqrt#1{ \ifdim #1<\:mmp   #1\z@  \else
   \:dd#1   \divide\:dd \tw@
   \def\::Sqrt{  \:ddd#1
      \:divide\:ddd\:dd      \:AbsDif\:d\:dd\:ddd
      \advance\:dd \:ddd   \divide\:dd \tw@
      \ifdim  \:d < \:mmmp
         \let\::Sqrt\relax  \fi
      \::Sqrt}
   \::Sqrt   #1\:dd   \fi }\Define\:length{
   \:dd \:AbsVal \:ddd
   \:abs\:x  \:abs\:y
   \ifdim \:dd<\:x \:dd\:x \fi
   \ifdim \:dd<\:y \:dd\:y \fi
   \ifdim \:dd>\:mp
      \:divide\:x\:dd     \:sqr\:x
      \:divide\:y\:dd     \:sqr\:y
      \:divide\:ddd\:dd   \:sqr\:ddd
      \advance\:x  \:y  \advance\:x \:ddd
      \:y\:dd  \:Sqrt\:x  \:d\:Cons\:x\:y
   \else \:d\:dd \fi}\Define\:distance(2){\MarkLoc(@^)
   \MoveToLoc(#1)  \:x \:X  \:y \:Y  \:ddd\:Z
   \MoveToLoc(#2)
   \advance\:x  -\:X   \advance\:y  -\:Y
   \advance\:ddd  -\if:IIID \:Z \else \:ddd \fi
   \:length  \MoveToLoc(@^)}\def\:NormalizeDeg#1{
   \:DL#1   \:K\:InCons\:DL
   \divide\:K  \:cccvx   \multiply\:K  \:cccvx
   \advance #1 -\:K\p@
   \ifdim #1<\z@ \advance #1  \:CCCVX \fi
   \ifdim #1=\z@
      \ifdim\:DL=\z@ \else
         \advance #1  \:CCCVX \fi \fi }\def\:CosSin#1{ \:DK#1
   \:NormalizeDeg\:DK \def\:tempA{}
\ifdim \:CVXXX<\:DK
   \def\:tempA{\:y-\:y}
   \advance\:DK -\:CCCVX  \:DK-\:DK   \fi
\ifdim \:XC<\:DK
   \edef\:tempA{\:x-\:x \:tempA}
   \advance\:DK -\:CVXXX  \:DK-\:DK   \fi
\ifdim 45\p@<\:DK
   \edef\:tempA{\:d\:x \:x\:y \:y\:d \:tempA}
   \advance\:DK -\:XC   \:DK-\:DK   \fi
   \:x\p@   \:y0.01745\:DK   \:d\:y   \:K\@ne
   \edef\:next{\advance\:K \@ne
      \:sqr\:d  \divide\:d \:K  \advance}
   \:next \:x -\:d   \:next \:y -\:d
   \:next \:x  \:d   \:next \:y  \:d
   \:next \:x -\:d   \:next \:y -\:d
   \:next \:x  \:d   \:next \:y  \:d
   \:tempA   }   \Define\:rInitOval(2){
   \XSaveUnits   \let\du=\:rdu
   \:Zunits\p@   \:dd#1\:Zunits
   \:d\:Cons\:dd\:Xunitsx   \edef\:Xx{\:Cons\:d}
   \:d\:Cons\:dd\:Xunitsy   \edef\:Xy{\:Cons\:d}
   \:dd#2\:Zunits
   \:d\:Cons\:dd\:Yunitsx   \edef\:Yx{\:Cons\:d}
   \:d\:Cons\:dd\:Yunitsy   \edef\:Yy{\:Cons\:d}
   \XRecallUnits  \let\du=\:svdu
}
\Define\:xyInitOval(2){
   \:d#1\:Xunits   \edef\:Xx{\:Cons\:d} \def\:Xy{0}
   \:d#2\:Yunits   \edef\:Yy{\:Cons\:d} \def\:Yx{0} }
 \def\:FigSize#1#2#3{
   \:x\:LBorder   \:y\:RBorder   \:d\:TeXLoc
   {\Object\:temp{#3}
    \setbox\:box\hbox{ \:temp
       \multiply\:x by \tw@  \multiply\:y by \tw@
       \xdef\:FSize{ \noexpand#1=\:Cons\:x;
                       \noexpand#2=\:Cons\:y;}}}
   \global\:LBorder\:x   \global\:RBorder\:y
   \global\:TeXLoc \:d
   \:FSize}
\expandafter\let \csname 0:Ln \endcsname\:Ln
\expandafter\def\csname 1:Ln \endcsname{
      \advance\:x -\:X  \advance\:y -\:Y
      \csname 0:Ln \endcsname(\:x,\:y)  }

\newcount\:ClipLevel   \:ClipLevel\@ne

\Define\Clip{\futurelet\:next\:Clip}

\Define\:Clip{
   \ifx \:next[  \expandafter\:DefClipOut
        \else    \expandafter\:DefClip     \fi }

\def\:DefClipOut[#1]{ \:DefClip(#1) }
\Define\:DefClip(1){  \def\:temp{#1}
   \ifx \:temp\empty
      \:ClipLevel\@ne     \def\:next{\let\:Ln}
      \expandafter\:next\csname 0:Ln \endcsname
   \else  \def\:temp{\::DefClip(#1)}  \fi  \:temp }
\Define\::DefClip(1){  \MarkLoc(^)
   \:x\:X  \:y \:Y  \Move(#1)
   \ifdim\:x>\:X  \:dd\:X  \:X\:x  \:x\:dd  \fi
   \ifdim\:y>\:Y  \:dd\:Y  \:Y\:y  \:y\:dd  \fi
   \advance\:ClipLevel  \@ne
   \expandafter\edef\csname \the
      \:ClipLevel :Ln \endcsname{
          \:xxx \the\:x   \:yyy \the\:y
          \:xxxx\the\:X  \:yyyy\the\:Y
          \ifx \:next[   \noexpand\:ClipOut
          \else          \noexpand\:ClipIn     \fi }
   \let\:Ln\:ClipLn   \MoveToLoc(^)  }\def\:ClipLn(#1,#2){
   \:x#1\:Xunits \:y#2\:Yunits
   {  \let\:Xunits\empty  \let\:Yunits\empty
   \advance\:x \:X   \advance\:y \:Y
   \ifdim \:x<\:X  \:dd\:X \:X\:x \:x\:dd
                   \:dd\:Y \:Y\:y \:y\:dd  \fi
   \:diff\:dd\:X\:x   \:diff\:ddd\:Y\:y
\:Z \:AbsVal \:dd
\advance\:Z  \:AbsVal\:ddd
\ifdim \:Z>\sixt@@n\p@
   \divide\:dd   128
   \divide\:ddd  128  \fi
\:Z\:Cons\:y\:dd
\advance\:Z -\:Cons\:x\:ddd
\ifdim \:dd<\z@
  \:dd-\:dd  \:ddd-\:ddd  \:Z-\:Z
\fi
   \csname \the \:ClipLevel :Ln \endcsname      } 
   \advance\:X \:x   \advance\:Y \:y  }\Define\:ClipIn{
   \def\:next{\let\:next}
   \expandafter\:next\csname \the
      \:ClipLevel :Ln \endcsname
   \advance\:ClipLevel \m@ne
   { \:ClipLeft\:xxxx \:ClipDown\:yyy \:ClipUp\:yyyy \:next }
   { \:ClipRight\:xxx \:ClipDown\:yyy \:ClipUp\:yyyy \:next  }
   { \:ClipDown\:yyyy \:next }
     \:ClipUp\:yyy    \:next }

\Define\:ClipOut{
   \:ClipLeft\:xxx      \:ClipRight\:xxxx
   \:ClipUp  \:yyyy     \:ClipDown \:yyy
   \def\:next{\let\:next}
   \expandafter\:next
      \csname \the\:ClipLevel :Ln \endcsname
   \advance\:ClipLevel \m@ne      \:next  }\def\:ClipLeft#1{
   \ifdim       \:x<#1  \:KilledLine
   \else \ifdim \:X<#1  \:X#1
      \ifdim \:dd>\:mmmp
         \:Y\:Cons\:ddd\:X  \advance\:Y \:Z
         \:divide\:Y\:dd
   \fi \fi \fi     \:CondKilLn  }

\def\:ClipRight#1{
   \ifdim       \:X>#1  \:KilledLine
   \else \ifdim \:x>#1  \:x#1
      \ifdim \:dd>\:mmmp
         \:y\:Cons\:ddd\:x  \advance\:y \:Z
         \:divide\:y\:dd
   \fi \fi \fi    \:CondKilLn }

\Define\:CondKilLn{
   \:d\:x  \advance\:d -\:X
   \ifdim \:d<\z@  \:d-\:d \fi
   \ifdim \:y<\:Y  \advance\:d \:Y \advance\:d -\:y
   \else           \advance\:d \:y \advance\:d -\:Y  \fi
   \ifdim  \:d<\:mmp   \:KilledLine \fi
   \ifdim  \:thickness=\z@ \:KilledLine \fi}

\Define\:KilledLine{
   \let\:ClipLeft\:gobble  \let\:ClipRight\:gobble
   \let\:ClipUp  \:gobble  \let\:ClipDown \:gobble
   \let\:next\relax
   \expandafter\def\csname 1:Ln \endcsname{}}\def\:ClipUp#1{
   \:AbsDif\:d\:y\:Y
   \ifdim  \:d<\:ragged
      \advance\:y  0.5\:thickness
\advance\:Y -0.5\:thickness
\ifdim       \:Y>#1  \:KilledLine
\else \ifdim \:y>#1
   \:thickness#1 \advance\:thickness -\:Y
   \advance\:Y  0.5\:thickness  \:y\:Y
\else
   \advance\:Y  0.5\:thickness
   \advance\:y -0.5\:thickness
\fi  \fi
\:dd\p@   \:ddd\z@
\def\:temp{  \:Z\:Y }  \:temp
   \else \let\:temp\relax
          \ifdim \ifdim\:Y<\:y\:Y\else\:y\fi >#1  \:KilledLine
   \else  \ifdim  \::ClipUp#1\:X\:Y
   \else  \ifdim  \::ClipUp#1\:x\:y
   \fi \fi \fi \fi   \:CondKilLn  }

\def\::ClipUp#1#2#3{
#3>#1   #3#1
\ifdim \:AbsVal\:ddd>\:mmmp
   #2\:Cons\:dd#3  \advance#2 -\:Z
   \:divide#2\:ddd
\fi  }\def\:ClipDown#1{   \:Ex2#1
   \:Flip\:y  \:Flip\:Y  \:ClipUp#1
   \:Flip\:y  \:Flip\:Y  \:temp }

\def\:Flip#1{   #1-#1  \advance#1 \:Ex  }
\Define\:SetDrawWidth{
   \hsize\:RBorder      \advance\hsize -\:LBorder
   \leftskip -\:LBorder \rightskip\z@}
\newdimen\:thickness   \:thickness0.75\p@

\Define\PenSize(1){\:thickness#1\relax}
\Define\:Draw{   \ifvmode \noindent\hfil\fi
   \global\:TeXLoc\z@
   \vbox\bgroup
           \begingroup
   \def\EndDraw{
           \endgroup   \:SetDrawWidth
         \egroup}
   \:DraCatCodes      \parindent\z@    \everypar{}
\leftskip\z@     \rightskip\z@    \boxmaxdepth\maxdimen
\linepenalty10   \let\FigSize\:FigSize
\def\Draw{\:wrn1{}} \:CommonIID   \:InDraw }
\Define\:CommonIID{\def\RotateTo{\:RotateTo}
\def\Rotate{\:Rotate}
\def\MoveF{\:MvF}\def\LineToLoc{\:LnToLoc}}\newdimen\:Xunits   \:Xunits\p@
\newdimen\:Yunits   \:Yunits\p@
\Define\:InDraw{\:Opt()\::InDraw{\:Xunits,\:Yunits}}
\Define\::InDraw(2){   \:X\z@   \:Y\z@
   \:Xunits#1 \relax  \:Yunits#2  \:AdjRunits
   \global\:LBorder\z@    \global\:RBorder\z@   \gdef\:UBorder{0pt}
   \:loadIID   \leavevmode }

\Define\:Resize(2){
   \:Xunits#1\:Xunits  \relax
   \:Yunits#2\:Yunits  \:AdjRunits  }
\Define\:Units(2){
   \:Xunits#1 \relax   \:Yunits#2    \:AdjRunits   }
\Define\:loadIID{
   
   \def\LineTo{\:LnTo}
   \def\MoveTo{\:MvTo}
   \def\Line{\:Ln}
   \def\Move{\:Mv}
   \def\MoveF{\:MvF}
   \:rotatedfalse\def\:InitOval{\:xyInitOval}
   
   \def\Units{\:Units}}\let\SaveAll\begingroup
\let\RecallAll\endgroup
\def\DrawOn{\def\Draw{\:Draw}}                      \DrawOn
\def\DrawOff{\def\Draw{\begingroup \:J\@cclv
                       \:NoDrawSpecials \:NoDraw}}
\catcode`\/0 \catcode`\\11
/def/:NoDraw#1\EndDraw{/endgroup}
/catcode`/\0 /catcode`//12

\def\:NoDrawSpecials{\catcode\:J11
  \ifnum \:J=\z@
     \let \:NoDrawSpecials\relax \fi
  \advance\:J  \m@ne \:NoDrawSpecials}
\let\:XunitsReg\:Xunits   \let\:YunitsReg\:Yunits  \Define\EntryExit(4){
   \edef\:InOut##1{
      \noexpand\ifcase ##1\space
         #1\noexpand\or #2\noexpand\or
         #3\noexpand\or #4\noexpand\fi}}

\EntryExit(0,0,0,0)\Define\:DrawBox{   \:x0.5\wd\:box
   \:y\ht\:box \advance\:y  \dp\:box
   \divide\:y \tw@
   \advance \:X -\:InOut0\:x
   \advance \:Y -\:InOut1\:y
   { \advance \:X -\:x
     {  \advance \:Y \:y
\ifdim\:Y>\:UBorder \xdef\:UBorder{\the\:Y}\fi
\ifdim\:X<\:LBorder \global\:LBorder\:X\fi
\advance\:X \wd\:box
\ifdim\:X>\:RBorder \global\:RBorder\:X\fi }
      \advance \:Y -\:y  \advance \:Y  \dp\:box
     {\:d\:X \advance\:d \wd\:box
 \advance\:X -\:TeXLoc   \global\:TeXLoc\:d
 \vrule width\:X depth\z@ height\z@
 \raise \:Y \box\:box} }
   \edef\MoveToExit(##1,##2){
   \:X\the\:X   \:Y\the\:Y
   \:x\the\:x   \:y\the\:y
   \advance\:X  ##1\:x
   \advance\:Y  ##2\:y}
   \advance \:X  \:InOut2\:x
   \advance \:Y  \:InOut3\:y }
 \Define\ThreeDim{\:Opt[]\:ThreeDim{\p@}}

\def\:ThreeDim[#1](#2){\::ThreeDim[#1](#2,,)}

\def\::ThreeDim[#1](#2,#3,#4,#5){ \bgroup\begingroup
   \def\EndThreeDim{          \endgroup\egroup}
   \:IIIDtrue   \:Zunits#1  \:Z\z@
   \def\:temp{#4}
   \ifx \:temp\empty   \:CosSin{#3\p@}  \:divide\:x\:y       \:Ey\:x
\:CosSin{#2\p@}  \:Ex\:Cons\:Ey\:x  \:Ey\:Cons\:Ey\:y
\let\:project\:projectPar

   \else               \:Ex#2\:Xunits \:Ey#3\:Yunits \:Ez#4\:Zunits
\let\:project\:projectPer
 \fi
   
\def\LineTo{\:tLnTo}
\def\MoveTo{\:tMvTo}
\def\Line{\:tLn}
\def\Move{\:tMv}

\def\Units{\:tUnits}\def\RotateTo{\:tRotateTo}
\def\Rotate{\:tRotate}
\def\MoveF{\:tMvF}
\:Vdirection\z@ \def\LineToLoc{\:tLnToLoc} }  \Define\:projectPer{
   \:diff\:x\:X\:Ex   \:diff\:y\:Y\:Ey
   \:diff\:xxxx\:Z\:Ez
   \:divide\:x\:xxxx      \:divide\:y\:xxxx
   \:x-\:Cons\:Ez\:x     \:y-\:Cons\:Ez\:y
   \advance\:x  \:Ex    \advance\:y \:Ey}
\Define\:projectPar{
   \:x\:Cons\:Ex\:Z  \advance\:x \:X
   \:y\:Cons\:Ey\:Z  \advance\:y \:Y}
\def\:tDLn(#1,#2,#3,{\:tMvTo(#1,#2,#3)
                     \:tLnTo(}
\Define\:tMvTo(3){   \:X#1\:Xunits
   \:Y#2\:Yunits    \:Z#3\:Zunits}
\Define\:tMv(3){     \advance\:X  #1\:Xunits
   \advance\:Y  #2\:Yunits
   \advance\:Z  #3\:Zunits}
\Define\:tResize(3){   \:Xunits#1\:Xunits
   \:Yunits#2\:Yunits \:Zunits#3\:Zunits
   \:AdjRunits}
\Define\:tUnits(3){    \:Xunits#1  \relax
   \:Yunits#2  \relax       \:Zunits#3
   \:AdjRunits}\Define\:tLnTo(3){
   \:project   \edef\:temp{\:x\the\:x \:y\the\:y}
   \:X#1\:Xunits  \:Y#2\:Yunits  \:Z#3\:Zunits
   \:DLN}

\Define\:tLn(3){
   \:project   \edef\:temp{\:x\the\:x \:y\the\:y}
   \advance\:X  #1\:Xunits
   \advance\:Y  #2\:Yunits
   \advance\:Z  #3\:Zunits   \:DLN}

\Define\:DLN{
   \TwoDim      \:temp        \LineTo(\:x\du,\:y\du)
   \EndTwoDim }\Define\TwoDim{\bgroup\begingroup
   \def\EndTwoDim{\endgroup\egroup}
   \:loadIID
   \if:IIID  \:IIIDfalse \:project \:X\:x \:Y\:y
             \:CommonIID \fi
   \Units(\:Xunits,\:Yunits)}  
\newif\if:rotated
\Define\RotatedAxes(2){\begingroup
   
   \if:IIID  \:IIIDfalse
             \:project \:X\:x \:Y\:y
             \:CommonIID    \fi
   \:DK#2\p@ \advance\:DK -#1\p@
\advance\:DK -\:CVXXX   \:NormalizeDeg\:DK
\ifdim \:DK<\:mmp \:wrn4{(#1,#2)}
\else
  \:DK#1\p@ \:NormalizeDeg\:DK    \:CosSin\:DK
  \:Xunitsx\:x  \:Xunitsy\:y
  \:DK#2\p@ \advance\:DK -\:XC \:NormalizeDeg\:DK \:CosSin\:DK
  \edef\Units(##1,##2){\noexpand\:Units(##1,##2)
    \:Xunitsx \:Cons\:Xunitsx\:Xunits
    \:Xunitsy \:Cons\:Xunitsy\:Xunits
    \:Yunitsx-\:Cons\:y\:Yunits
    \:Yunitsy \:Cons\:x\:Yunits  } \fi
 \def\MoveTo{\:rMvTo}
\def\Move{\:rMv}
\def\LineTo{\:rLnTo}
\def\Line{\:rLn}

\def\MoveF{\:rMvF}\def\:InitOval{\:rInitOval}

   \MarkLoc(:org) \Units(\:Xunits,\:Yunits)  \:rotatedtrue }
\newdimen\:Xunitsx  \newdimen\:Xunitsy
\newdimen\:Yunitsx  \newdimen\:Yunitsy
\Define\:rResize(2){
   \:Xunits #1\:Xunits    \:Yunits #2\:Yunits
   \:Xunitsx#1\:Xunitsx   \:Xunitsy#1\:Xunitsy
   \:Yunitsx#2\:Yunitsx   \:Yunitsy#2\:Yunitsy  }
\def\:rdu#1{%
   \ifx #1\:Xunits \units:i\else
   \ifx #1\:Yunits \units:i\else
   \ifx #1\:Zunits \units:i\else #1\fi \fi \fi }
\def\units:i{\Units(\p@,\p@)}
\Define\:rMvTo{\MoveToLoc(:org) \:rMv}
\Define\:rMv(1){ \:rxy(#1)
   \advance\:X  \:x \advance\:Y  \:y}
\Define\:rLnTo(1){ \MarkLoc(:a) \:rMvTo(#1) {\LineToLoc(:a)} }
\Define\:rLn(1){ \:rxy(#1) \:Ln(\:x\du,\:y\du) }
\def\:rDLn(#1,#2,{\:rMvTo(#1,#2)\:rLnTo(}
\Define\:rxy(2){                           \XSaveUnits   \let\du=\:rdu
   \:Zunits\p@   \:Zunits#1\:Zunits
   \:x\:Cons\:Zunits\:Xunitsx
   \:y\:Cons\:Zunits\:Xunitsy            \XRecallUnits  \XSaveUnits
   \:Zunits\p@   \:Zunits#2\:Zunits
   \advance\:x  \:Cons\:Zunits\:Yunitsx
   \advance\:y  \:Cons\:Zunits\:Yunitsy \XRecallUnits  \let\du=\:svdu}
\Define\:rMvF(1){  \:CosSin\:direction   \XSaveUnits   \let\du=\:rdu
   \:Zunits\p@            \:Zunits#1\:Zunits
   \:x\:Cons\:Zunits\:x   \:y\:Cons\:Zunits\:y
   \edef\:temp{(\:Cons\:x,\:Cons\:y)}   \expandafter\Move\:temp
   \XRecallUnits  \let\du=\:svdu}
\Define\:AdjRunits{
   \:Xunitsx\:Xunits   \:Xunitsy\z@
   \:Yunitsx\z@        \:Yunitsy\:Yunits}
\Define\Text{  \setbox\:box
   \vtop\bgroup    \edef\DoReg{\the\DoReg}
      \hyphenpenalty\@M  \exhyphenpenalty\@M
      \catcode`\ 10 \catcode`\^^M13 \catcode`\^^I10
      \catcode`\&4  \let~\space
      \:Text}                       \catcode`\^^M13 %
\def\:Text(--#1--){%
      \:SetLines#1\hbox{}^^M--)^^M %
   \egroup                  %
   \if:IIID \TwoDim  \:DrawBox  \EndTwoDim %
   \else             \:DrawBox  \fi}       %
\def\:SetLines#1^^M{        %
   \def\:TextLine{#1}       %
   \ifx \:TextLine\:LastLine   \let\:temp\relax      %
   \else  \def\:temp{                                 %
             \:IndirectLines#1\relax~~--)~~\:SetLines}%
   \fi  \:temp }                      \catcode`\^^M9
\def\:IndirectLines#1~~{    \def\:TextLine{#1}
  \ifx \:TextLine\:LastLine   \let\:temp\relax
  \else  \def\:temp{\:AddLine{#1}\:IndirectLines}
  \fi  \:temp }

\Define\:LastLine{--)}

\def\:AddLine#1{
   \ifvmode \noindent    \hsize\z@ \else
      \hfil \penalty-500 \hbox{}    \fi
   \hfil#1
   \setbox\:box\hbox{#1}
   \ifdim \wd\:box>\hsize \hsize\wd\:box \fi}\def\TextPar#1#2{
   \def\:TxtPar##1(##2){##1(--##2--)}
   \edef\:temp{\expandafter\noexpand\csname :\string#2\endcsname}
   \edef#2{\noexpand\:TextPar\expandafter\noexpand\:temp}
   \expandafter\let\:temp\:undefined
   \expandafter#1\:temp}
                                   \catcode`\^^M13
\def\:TextPar#1{\begingroup     \catcode`\&4         %
   \catcode`\ 10 \catcode`\^^M13 \catcode`\^^I10 %
   \:TPar{#1}}                     \catcode`\^^M9  %

\def\:TPar#1(--#2--){\endgroup
   #1(--#2--)  }
 \newdimen\:direction  \newdimen\:Vdirection
\Define\:Rotate(1){   \advance \:direction   #1\p@
                      \:NormalizeDeg\:direction  }
\Define\:RotateTo(1){ \:direction  #1\p@
                      \:NormalizeDeg\:direction  }

\Define\:tRotate(2){
   \advance\:Vdirection  #2\p@
   \:NormalizeDeg\:Vdirection  \:Rotate(#1)}
\Define\:tRotateTo(2){
   \:Vdirection #2\p@
   \:NormalizeDeg\:Vdirection  \:RotateTo(#1)}
\Define\LineF(1){\MarkLoc(,)\MoveF(#1) {\LineToLoc(,)}}
\Define\:MvF(1){    \:CosSin\:direction
   \:d#1\:Xunits   \advance\:X  \:Cons\:d\:x
   \:d#1\:Yunits   \advance\:Y  \:Cons\:d\:y   }
\Define\MoveFToOval(2){
   \:NormalizeDeg\:direction \:CosSin\:direction  \:Zunits\p@
   \:d#2\:Zunits    \:x\:Cons\:d\:x
   \ifdim \:d=\z@  \:err5{...,\:Cons\:d} \fi
   \:d#1\:Zunits    \:y\:Cons\:d\:y
   \ifdim \:d=\z@  \:err5{\:Cons\:d,...} \fi     \:SearchDir
   \XSaveUnits
      \if:rotated
        \:Zunits\p@     \:Zunits#1\:Zunits
        \:x\:Cons\:Zunits\:Xunits
        \:Zunits\p@     \:Zunits#2\:Zunits
        \:y\:Cons\:Zunits\:Yunits
      \else
        \:x#1\:Xunits  \:y#2\:Yunits  \fi     \relax
      \Units(\:x,\:y)
      \:Vdirection\:direction  \:direction\:ddd
      \MoveF(\@ne)  \:direction\:Vdirection  \XRecallUnits  }
\Define\:tMvF(1){    \:CosSin\:Vdirection
   \:d #1\:Zunits   \advance\:Z  \:Cons\:y\:d
   \:xx#1\:Xunits     \:yy#1\:Yunits
   \:xx\:Cons\:x\:xx    \:yy\:Cons\:x\:yy
   \:CosSin\:direction
   \:x\:Cons\:xx\:x    \:y\:Cons\:yy\:y
   \advance\:X  \:x   \advance\:Y  \:y } \Define\CSeg{\:Opt[]\:CSeg1}
\def\:CSeg[#1]#2(#3,#4){   \MarkLoc($^)
   \MoveToLoc(#4) \:x\:X \:y\:Y
   \if:IIID  \:d\:Z  \fi   \MoveToLoc(#3)
   \advance\:x -\:X  \:x#1\:x
   \advance\:y -\:Y  \:y#1\:y
   \if:IIID   \advance\:d -\:Z  \:d#1\:d \fi
   \:t{#2}
   \edef\:temp{\the\:t(
      \expandafter\:Frac\the\:x\noexpand\:du,
      \expandafter\:Frac\the\:y\noexpand\:du \if:IIID ,
      \expandafter\:Frac\the\:d\noexpand\:du \fi)}
   \MoveToLoc($^)     \:temp}\Define\LSeg{\:Opt[]\:LSeg1}
\def\:LSeg[#1]#2(#3,#4){   \:distance(#3,#4)
   \:d#1\:d  \:t{#2}
   \edef\:temp{\the\:t(\expandafter\:Frac\the\:d\noexpand\:du)}
   \:temp}\Define\DSeg{\:Opt[]\:DSeg1}

\def\:DSeg[#1]#2(#3,#4){   \MarkLoc(^)
   \MoveToLoc(#4)  \:xxx\:X   \:yyy\:Y \:xxxx\:Z
   \MoveToLoc(#3)
   \advance\:xxx -\:X   \advance\:yyy -\:Y
   \ifdim \:AbsVal\:xxx<\:mmmp
      \ifdim \:AbsVal\:yyy<\:mmmp \:wrn5{#3,#4}
   \fi\fi
   \if:IIID
      \advance \:xxxx  -\:Z
      \:divide\:xxx\:Xunitsx
      \:divide\:yyy\:Yunitsy
      \:divide\:xxxx\:Zunits
      \:x\:xxx  \:y\:yyy   \:ddd\z@  \:length
      \:x\:d    \:y\:xxxx  \:SearchDir
      \:yyyy\:ddd
      \:x\:xxx  \:y\:yyy
   \else
      \:x  \:AbsVal\:Yunitsy
\:y  \:AbsVal\:Xunitsy  \ifdim \:y>\:x  \:x\:y  \fi
\:y  \:AbsVal\:Yunitsx  \ifdim \:y>\:x  \:x\:y  \fi
\:y  \:AbsVal\:Xunitsx  \ifdim \:y>\:x  \:x\:y  \fi
\:K  \:InCons\:x  \relax
      \ifnum \:K<\thr@@     \:K\@ne
\else \ifnum \:K<\sixt@@n   \:K4
\else \ifnum \:K<\:XC       \:K\sixt@@n
\else \ifnum \:K<\@m        \:K\@cclvi
\fi \fi \fi \fi
\divide\:xxx \:K
\divide\:yyy \:K  
      \:x \:Cons\:Yunitsy\:xxx
\advance\:x -\:Cons\:Yunitsx\:yyy
      \:y-\:Cons\:Xunitsy\:xxx
\advance\:y \:Cons\:Xunitsx\:yyy     
   \fi
   \:SearchDir   \:ddd#1\:ddd   \:t{#2}
   \edef\:temp{\the\:t(
      \:Cons\:ddd \if:IIID ,\:Cons\:yyyy \fi)}
   \MoveToLoc(^) \:temp}   \def\:theDoReg{\def\DoReg{\the\:DoReg}}

\Define\MarkLoc{  \:theDoReg
   \expandafter\edef \csname \:MarkLoc}

\Define\MarkGLoc{  \:theDoReg
   \expandafter\xdef \csname \:MarkLoc}

\Define\:MarkLoc(1){ Loc\space#1:\endcsname{
   \:X\the\:X  \:Y\the\:Y
   \if:IIID \:Z\the\:Z \fi  }       \let\DoReg\:DoReg }

\Define\MoveToLoc(1){   \:theDoReg   \expandafter\ifx
   \csname Loc\space#1:\endcsname\relax \:err2{#1}\fi
   \csname Loc\space#1:\endcsname     \let\DoReg\:DoReg  }

\Define\MarkPLoc(1){  \:theDoReg
   \if:IIID   \:project
      \expandafter\edef \csname Loc\space#1:\endcsname{
         \:X\the\:x  \:Y\the\:y  \:Z\z@}
   \else \:err1\MarkPLoc  \fi   \let\DoReg\:DoReg}

\Define\WriteLoc(1){{ \:theDoReg \edef\:temp{#1}
   \ifx \:temp\empty \else \MoveToLoc(#1) \fi
   \immediate\write\sixt@@n{...
      \:temp=(\the\:X,\the\:Y\if:IIID,\the\:Z\fi)}}}

\Define\:LnToLoc(1){
   \:x\:X \:y\:Y   \MoveToLoc(#1)
   { \:LnTo(\:x\du,\:y\du) }}

\Define\:tLnToLoc(1){
   \:xx\:X \:yy\:Y \:dd\:Z     \MoveToLoc(#1)
   { \:tLnTo(\:xx\du,\:yy\du,\:dd\du) }}\def\:GetLine#1#2#3#4#5{
   \MoveToLoc(#1)   \divide\:X \:eight  \divide\:Y \:eight
   #3\:X  #4\:Y
   \MoveToLoc(#2)   \divide\:X \:eight  \divide\:Y \:eight
   \advance #3 -\:X  \advance #4 -\:Y
   #5\:Cons#3\:Y      \advance #5 -\:Cons#4\:X
   \divide #3 \:eight     \divide  #4 \:eight \relax      }
\def\MoveToLL(#1,#2)(#3,#4){
   \:GetLine{#1}{#2}\:x \:y \:xxx
   \:GetLine{#3}{#4}\:xx\:yy\:xxxx
   \:ddd \:Cons\:x \:yy     \advance\:ddd -\:Cons\:xx\:y
   \ifdim  \:AbsVal\:ddd < \:mmmp
      \:X\@cclv\p@  \:Y\:X
      \:wrn3{(\string#1,\string#2)(\string#3,\string#4)}
   \else
      \:divide\:xxx\:ddd    \:divide\:xxxx\:ddd
      \:X\:Cons\:xxx\:xx   \advance\:X -\:Cons\:xxxx\:x
      \:Y\:Cons\:xxx\:yy   \advance\:Y -\:Cons\:xxxx\:y
   \fi   }\Define\MoveToCC{\:Opt[]\:MoveToCC{}}

\def\:MoveToCC[#1](#2,#3)(#4,#5){
  \:UserUnits(#2,#3)(#4,#5)
  \:distance($#2,$#4)
\ifnum \:d<\:mp \MoveToLoc(#3)
   \:wrn3{(#1,#2)(#3,#4)}
\else                 \:xx \:d
   \:distance($#2,$#3)  \:xxx\:d
   \:distance($#5,$#4) 
     \:yy \:xxx  \advance\:yy -\:d
\:yyy\:xxx  \advance\:yyy \:d
\:divide\:yy\:xx     \:yy\:Cons\:yy\:yyy
\advance\:yy \:xx  \divide\:yy \tw@
     \:yyy \ifdim \:AbsVal\:xxx>\:AbsVal\:yy \:xxx \else \:yy \fi
\ifdim \:AbsVal\:yyy<\:mp \:yyy\z@ \else
   \:divide\:xxx\:yyy  \:sqr\:xxx
   \:yyyy\:yy  \:divide\:yyyy\:yyy  \:sqr\:yyyy
   \advance\:xxx -\:yyyy   \:Sqrt\:xxx
   \:yyy\:Cons\:yyy\:xxx
\fi 
     \MoveToLoc($#4)  \:x\:X  \:y\:Y  \:divide\:yy\:xx
\MoveToLoc($#2)
\advance\:x -\:X     \advance\:y -\:Y
\advance\:X \:Cons\:yy\:x
\advance\:Y \:Cons\:yy\:y
     \:divide\:yyy\:xx
\advance\:X  #1\:Cons\:yyy\:y
\advance\:Y -#1\:Cons\:yyy\:x

  \fi  \:SysUnits  }\def\:UserUnits(#1,#2)(#3,#4){
   \:xx\:Xunitsx  \:xxx\:Xunitsy
   \:yy\:Yunitsx  \:yyy\:Yunitsy
   \:xxxx\:Cons\:yy\:xxx  \advance\:xxxx \:Cons\:yyy\:xx
   \ifdim \:AbsVal\:xxxx>\:mmmp
      \:divide\:xx\:xxxx   \:divide\:xxx{-\:xxxx}
      \:divide\:yy{\:xxxx} \:divide\:yyy\:xxxx
   \fi
   \:UnLoc(#1)  \:UnLoc(#2)
   \:UnLoc(#3)  \:UnLoc(#4)}\Define\:SysUnits{
   \MoveTo(\:Cons\:X,\:Cons\:Y)}\Define\:UnLoc(1){
   \MoveToLoc(#1)  \:d\:Cons\:yyy\:X
   \advance\:d \:Cons\:yy\:Y
   \:Y\:Cons\:xx\:Y
   \advance\:Y  \:Cons\:xxx\:X   \:X\:d
   \MarkLoc($#1)}

\def\:MoveToLC[#1](#2,#3)(#4,#5){
   \:UserUnits(#2,#3)(#4,#5)
   \MoveToLoc($#2)  \:x\:X  \:y\:Y
\MoveToLoc($#3)  \advance\:x -\:X
                \advance\:y -\:Y
   \edef\:temp{ \:xxx\the\:x  \:yyy\the\:y }
\MoveToLoc($#4)
\advance\:X \:y  \advance\:Y -\:x
\MarkLoc(^$) \MoveToLL($#4,^$)($#2,$#3)
\MarkLoc(^$)  
   \:distance($#4,$#5)  \:xx\:d  \:distance($#4,^$)
\ifdim      \:d>\:xx        \:wrn3{(#2,#3)(#4,#5)}
\else \ifdim \:d<\:mmp     \:yy\:xx   \else   \:yy\:d
      \:divide\:yy\:xx  \:sqr\:yy
      \:yy-\:yy   \advance\:yy \p@
      \:Sqrt\:yy  \:yy\:Cons\:xx\:yy
\fi \fi
   \:temp   \:x\:xxx  \:y\:yyy  \:length  \:xx\:d
\:divide\:xxx\:xx
\:divide\:yyy\:xx
\advance\:X  #1\:Cons\:yy\:xxx
\advance\:Y  #1\:Cons\:yy\:yyy
   \:SysUnits }  \def\Object#1{\:Opt(){\:DefineSD#1}0}

\def\:DefineSD#1(#2){\begingroup  \:multid#1
   \:DraCatCodes   \:DefSD#1(#2)}

\def\:DefSD#1(#2)#3{
   \expandafter\::Define\csname\string#1.\endcsname(#2){
      \:t{\:SubD{#3}}
      \if:IIID \edef\:temp{\noexpand\TwoDim \the\:t
                           \noexpand\EndTwoDim}
      \else    \def\:temp{\the\:t}   \fi         \:temp}
   \def#1{\def\:SDname{\csname\string#1.\endcsname}
          \:Opt[]\:CallSD{}}}

\def\:CallSD[#1]{ \edef\:Entry{#1} \:SDname }\def\:SubD#1{
   \let\::RecallXLoc\:AddXLoc   \gdef\:AddXLoc{}
   \edef\:RecallBor{ \global\:LBorder  \the\:LBorder
                  \global\:RBorder  \the\:RBorder
                  \xdef\noexpand\:UBorder{\:UBorder}
                  \global\:TeXLoc\the\:TeXLoc }
\global\:TeXLoc\z@
\setbox\:box\vbox{\EntryExit(0,0,0,0)
\begingroup
   \def\MarkXLoc{\:MarkXLoc}
   \:InDraw  #1
\endgroup
\:SetDrawWidth                    \let\:XLoc\relax
\xdef\:AddXLoc{\:dd\the\:LBorder  \:AddXLoc}}
\:RecallBor
   \:ddd\dp\:box
   \ifx \:Entry\empty
   \:DrawBox
\else
   \let\:RecallIn\:InOut
   \:x\:X    \:y\:Y
   \def\:XLoc(##1,##2,##3){
      \def\:temp{##1}
      \ifx \:temp\:Entry \:X\:x  \advance\:X -##2
                         \:Y\:y  \advance\:Y -##3
      \fi}
   \:AddXLoc
   \advance\:X  \:dd      \advance\:Y -\:ddd
   \EntryExit(-1,-1,\:InOut2,\:InOut3)  \:DrawBox
   \let\:InOut\:RecallIn
\fi
   \MarkLoc(^)
      \MoveToExit(-1,-1)
   \:xxx\:X   \:yyy\:Y   \advance\:yyy  \:ddd
   \def\:XLoc(##1,##2,##3){
      \:X\:xxx  \advance\:X  ##2     \advance\:X -\:dd
      \:Y\:yyy  \advance\:Y  ##3
      \MarkLoc(##1)}
   \:AddXLoc
   \MoveToLoc(^)
   \ifx \:Entry\empty \else     \MoveToLoc(\:Entry) \fi
   \global\let\:AddXLoc\::RecallXLoc }
\Define\:MarkXLoc(1){  \:theDoReg
   \let\:XLoc\relax
   \xdef\:AddXLoc{\:AddXLoc \:XLoc(#1,\the\:X,\the\:Y)}
   \let\DoReg\:DoReg} 
   \catcode`\ 10 \catcode`\^^M5 \catcode`\^^I10
   \let\wlog\:wlog  \let\:wlog\:undefined
\:RestoreCatcodes     \tracingstats1  

\expandafter\edef\csname :RestoreCatcodes\endcsname{%
   \catcode`\noexpand\noexpand\noexpand \^=\the\catcode`\^%
}
\catcode`\^=7
\expandafter\edef\csname :RestoreCatcodes\endcsname{%
   \csname :RestoreCatcodes\endcsname
   \catcode`\noexpand \_=\the\catcode`\_%
   \catcode`\noexpand :=\the\catcode`:%
   \catcode`\noexpand &=\the\catcode`&%
   \catcode`\noexpand \#=\the\catcode`\#%
   \catcode`\noexpand \^^M=\the\catcode`\^^M%
   \let\expandafter\noexpand
       \csname:RestoreCatcodes\endcsname=\noexpand\undefined}
\catcode`\:=11
\ifx \AlDraTex\:undefined  \def\AlDraTex{chart,diagram} \fi

\newif\if:circle
\newif\if:option
\def\:CheckOption#1{ \def\:temp{#1}
   \:optionfalse
   \expandafter\:GetOptions\AlDraTex,,// }

\def\:GetOptions#1,#2//{  \def\:next{#1}
   \ifx        \:next\empty
   \else \ifx  \:temp\:next   \:optiontrue   \let\:next=\relax
   \else \def\:next{\:GetOptions#2//}
   \fi  \fi  \:next }
\def\:DefineExt#1{%
   \Define#1{\begingroup  \:DraCatCodes
             \csname .:\string#1\endcsname}%
  \expandafter\Define\csname .:\string#1\endcsname}\IntVar\:C   \IntVar\:Ca
\IntVar\:Cb  \IntVar\:Cc  \newtoks\:tk
\DecVar\Va   \DecVar\Vb
\DecVar\:Vc   \DecVar\:Vd   \DecVar\:Ve

\Define\:GetX(2){\Va=#1;}
\Define\:GetY(2){\Vb=#2;}
\Define\:GetXY(2){\Va=#1;\Vb=#2;}\Define\:GetExit(2){
    \:Vc=#1;  \IF \LtDec(\:Vc,0)  \THEN  \:Vc=-\:Vc;  \FI
    \:Vd=#2;  \IF \LtDec(\:Vd,0)  \THEN  \:Vd=-\:Vd;  \FI
    \IF  \LtDec(\:Vd,\:Vc)  \THEN  \:Vd=\:Vc;  \FI
    \Va=#1;   \Vb=#2;
    \IF \GtDec(\:Vd,0)  \THEN  \Va/\:Vd;  \Vb/\:Vd;   \FI }
\Define\:GetEdgeFront(1){\def\:EdgeFront{#1}}
\Define\:Count(1){\:C+1;}\let\Defend=\noexpand

\def\RemoveSpaces{\let\:Spaces=\relax}
\RemoveSpaces\def\:edef#1#2{\edef\:Temp{\noexpand\Define
   \noexpand#1{#2}}\:Temp }
\:CheckOption{chart}\if:option
   \edef\AlDraTex{\AlDraTex,pie,bar,xy}
\else \:CheckOption{pie}\if:option
   \edef\AlDraTex{\AlDraTex,chart}
\else \:CheckOption{xy}\if:option
   \edef\AlDraTex{\AlDraTex,chart}
\else \:CheckOption{bar}\if:option
   \edef\AlDraTex{\AlDraTex,chart}
\fi\fi\fi\fi
\:CheckOption{chart}\if:option
   \def\Compute#1(#2){\def\NextDefine{\::Comp#1(#2)}\Define\:Comp(#2)}

\def\::Comp#1(#2){\Define\:Temp(#2){}\ifx \:Temp\:Comp
         \let\:Compute=\empty
   \else \Define\:Compute(1){\:Comp(##1)\:Return#1}\fi}

\Compute\I(1){}\Define\:ChooseColor(1){
\IF \EqDec(#1,1) \THEN  \SetBrush(0,0){}
\ELSE\IF \EqDec(#1,2) \THEN
   \SetBrush(2pt\du,1pt\du){  \PenSize(0.2pt)
        \Units(1pt,1pt)  \Move(0,-0.37)  \Line(0,0.75)}
\ELSE\IF \EqDec(#1,3) \THEN
   \SetBrush(1pt\du,2pt\du){  \PenSize(0.2pt)
        \Units(1pt,1pt)  \Move(-0.37,0)  \Line(0.75,0)}
\ELSE\IF \EqDec(#1,4) \THEN
   \SetBrush(2pt\du,2pt\du){  \PenSize(0.2pt)
        \Units(1pt,1pt){\Move(-1,0)\Line(2,0)}
               \Move(0,-1)\Line(0,2)}
\ELSE\IF \EqDec(#1,5) \THEN
   \SetBrush(5pt\du,2pt\du){  \PenSize(0.2pt)
        \Units(1pt,1pt)  \Move(0,-0.37)  \Line(0,0.75)}
\ELSE\IF \EqDec(#1,6) \THEN
   \SetBrush(2pt\du,5pt\du){  \PenSize(0.2pt)
        \Units(1pt,1pt)  \Move(-0.37,0)  \Line(0.75,0)}
\ELSE
   \SetBrush(#1pt\du,#1pt\du){
                \Units(1pt,1pt)\PaintRect(0.37,0.37)}
\FI\FI\FI\FI\FI\FI} 
\fi
\:CheckOption{pie}\if:option
   \Define\PieChart{\begingroup   \:Spaces  \:PieChart}

\Define\:PieChart(1){\endgroup
   \Table\:Temp{#1}
   \:Ve=0;  \:Temp(0,99){\:total}
   \MarkLoc(o:) \MoveF(\:PieR)  \MarkLoc(b:)
   \DSeg\:Vd(o:,b:) \MoveToLoc(o:)
   \:C=0;  \:Temp(0,99){\:InsertSlice} }  \Define\PieChartSpec(1){
   \IF  \EqText(,#1)  \THEN\ELSE   \:PieSpec(#1)
   \FI  \:PieCol}

\Define\:PieSpec(3){
   \:edef\:PieR{#2}      \:edef\:LblDist{#3}
   \IF  \EqInt(1,#1)  \THEN
      \def\:ReturnComp##1{\:slice(\Val##1,\Val##1)}
      \Define\:total(1){
         \ifx  \:Compute\empty          \:Ve+##1;
         \else \let\:Return=\:ReturnAdd \:Compute(##1)\fi}
      \Define\:InsertSlice(1){
         \ifx \:Compute\empty           \:slice(##1,##1)
         \else \let\:Return=\:ReturnComp\:Compute(##1) \fi}
   \ELSE
      \def\:ReturnComp##1{\edef\:tempa{\Val##1}}
      \Define\:total(2){
         \ifx  \:Compute\empty          \:Ve+##1;
         \else \let\:Return=\:ReturnAdd \:Compute(##1)\fi}
      \Define\:InsertSlice(2){
         \ifx \:Compute\empty           \:slice(##1,##2)
         \else \let\:Return=\:ReturnComp\:Compute(##1)
                                       \:slice(\:tempa,##2) \fi}
   \FI  }
\def\:ReturnAdd#1{\:Ve+#1;}\Define\:slice(2){  \:C+1;  \Va=#1;  \Va/\:Ve;
    \Va*360; {
       \Va/2;  \Rotate(\Val\Va)   \:DetachSlice
       \MarkLoc(o:)
       \MoveF(1pt\du)  \MarkLoc(x:)  \CSeg\:SetEnEx(x:,o:)
       \MoveToLoc(o:)  \MoveF(\:LblDist)  \:AdjChatLbl(0,99){\:MvLbl}
       \SliceText(--#2--)  \MoveToLoc(o:)
       \Vb=\:Vd; \Vb+\Va;  \Vb+\Va;
       \DrawOvalArc(\:PieR,\:PieR)(\Val\:Vd,\Val\Vb)
       \:Ca=0; \:PieColors(0,99){\:Ca+1;\:PaintPie}
   {\RotateTo(\Val\:Vd)  \LineF(\:PieR)}
       \RotateTo(\Val\Vb)  \LineF(\:PieR)
    }  \Rotate(\Val\Va)  \:Vd+\Va;  }

\Define\:SetEnEx(2){\:Vc=#1; \Vb=#2; \EntryExit(\Val\:Vc,\Val\Vb,0,0)}
\Define\DetachSlices(2){
   \IF  \EqInt(#2,0)  \THEN
      \Define\:DetachSlice{}
   \ELSE\IF  \EqText(,#1) \THEN
      \def\:DetachSlice{ \MoveF(#2) }
   \ELSE
      \Table\:DSlice{#1}
      \Define\:OptDSlice(1){
         \IF \EqInt(\:C,##1) \THEN
            \MoveF(#2)
         \FI}
      \Define\:DetachSlice{ \:DSlice(0,99){\:OptDSlice} }
   \FI  \FI } \Define\MovePieLabels(1){
   \IF  \EqText(,#1)  \THEN
      \Table\:AdjChatLbl{0,0,0}
   \ELSE     \Table\:AdjChatLbl{0,0,0&#1}  \FI}

\Define\:MvLbl(3){
   \IF \EqInt(\:C,#1)  \THEN  \Move(#2,#3)  \FI }
\:DefineExt\:PieCol(1){\endgroup
   \IF \EqText(,#1)  \THEN
   \ELSE                   \Table\:PieColors{#1} \FI}
\Define\:PaintPie(1){{
   \IF \EqInt(\:C,\:Ca) \THEN
      \IF \GtDec(#1,0) \THEN
         \:ChooseColor(#1)
          \:Ve=\:PieR;
          \IF \GtDec(\Va,90) \THEN
             \Va=\:Vd;  \Va+181;
             \PaintOvalArc(\Val\:Ve,\Val\:Ve)(\Val\:Vd,\Val\Va)
             \Va-1;  \:Vd=\Va;    \FI
          \PaintOvalArc(\Val\:Ve,\Val\:Ve)(\Val\:Vd,\Val\Vb)
          \RotateTo(\Val\Vb)  \MoveF(\:PieR)  \MarkLoc(a:)
\MoveToLoc(o:)  \RotateTo(\Val\:Vd)  \MoveF(\:PieR)
\MarkLoc(b:)    \PaintQuad(o:,o:,a:,b:) 
          \MoveF(-0.5pt\du)  \MarkLoc(b:) \MoveToLoc(o:)
\Rotate(1) \MoveF(\:PieR)  \MarkLoc(b':)
\MoveToLoc(o:) \RotateTo(\Val\Vb)  \MoveF(\:PieR)
\MoveF(-0.5pt\du)  \MarkLoc(a:) \MoveToLoc(o:)
\Rotate(-1) \MoveF(\:PieR) \MarkLoc(a':)
\PaintQuad(a:,a':,b':,b:)
   \FI \FI  }}\Define\Legend(1){   \AlignGrid(-1,\:LgEn)
   \PictNode(2){  \EntryExit(-1,\:LgEn,0,0)
      \DrawRect(\:LgSzX,\:LgSzY)
      \IF  \GtDec(##1,0)  \THEN
         \:ChooseColor(##1)
         \IF  \EqDec(##1,4)  \THEN
            \Move(0.05,0.05)  \:Vd=\:LgSzX;  \:Ve=\:LgSzY;
            \:Vd-0.1;  \:Ve-0.1;
            \PaintRect(\Val\:Vd,\Val\:Ve) \Move(-0.05,-0.05)
         \ELSE  \PaintRect(\:LgSzX,\:LgSzY)  \FI
      \FI
      \:Vd=\:LgSzX;  \:Vd+5;
      \:Ve=\:LgEn;  \:Ve+1;  \:Ve*\:LgSzY;  \:Ve/2;
      \Move(\Val\:Vd,\Val\:Ve)
 \Text(--\ignorespaces##2--)    }
    \GridDiagram(#1)()()}

\:DefineExt\LegendSpec(3){\endgroup
    \:edef\:LgSzX{#1}\:edef\:LgSzY{#2}\:edef\:LgEn{#3}}

\LegendSpec(10,10,0) 
\PieChartSpec(1,40,20)(0)
\TextPar\Define\SliceText(1){\Text(--#1--)}
\DetachSlices(,0) \MovePieLabels()

\fi
\:CheckOption{xy}\if:option
   \Define\Axis(2){     \MoveToLoc(#1)   \CSeg\:GetXY(#1,#2)
   \begingroup   \:Spaces  \:axis}

\Define\:axis(2){ \endgroup\:axs#1//{#2}}

\def\:axs#1#2//#3{\::axs(#1,#2,#3)}
\Define\::axs(3){ \SaveAll  \Units(1pt,1pt)
   \def\:tic{\:ticLn( \Val\:Vc,0)}
\IF       \EqText(W,#1)  \THEN
   \def\:MvTicTxt{\Move(-\:TxtPos,0)}
   \IF  \LtDec(\Va,0)  \THEN \Move(\Val\Va,0)    \:Vc=-\Va;
   \ELSE    \:Vc=\Va;   \FI
   \EntryExit(1,0,0,0)  \Va=0;
\ELSE\IF  \EqText(E,#1)  \THEN
   \IF  \GtDec(\Va,0)  \THEN
      \Move(\Val\Va,0)     \:Vc=-\Va;
   \ELSE     \:Vc=\Va;  \FI
   \def\:MvTicTxt{\Move(\:TxtPos,0)}
   \Va=0;   \EntryExit(-1,0,0,0)
\ELSE
   \def\:tic{\:ticLn(0, \Val\:Vc)}
   \IF  \EqText(N,#1)  \THEN
      \IF  \GtDec(\Vb,0)  \THEN
         \Move(0,\Val\Vb)       \:Vc=-\Vb;
      \ELSE        \:Vc=\Vb;  \FI
      \EntryExit(0,-1,0,0)
      \def\:MvTicTxt{\Move(0,\:TxtPos)}
   \ELSE     \EntryExit(0,1,0,0)
      \IF  \LtDec(\Vb,0)  \THEN
         \Move(0,\Val\Vb)    \:Vc=-\Vb;
      \ELSE    \:Vc=\Vb;       \FI
      \def\:MvTicTxt{\Move(0,-\:TxtPos)}
   \FI
   \Vb=0;
\FI\FI
   \:Ca=#2;   \IF \EqText(#2,-0) \THEN   \:Ca=-1; \FI
   \IF  \LtInt(\:Ca,0)  \THEN
      \:Ca=-#2;  { \Line(\Val\Va,\Val\Vb)
\MarkLoc(bk:)
\IF  \EqDec(\Va,0)  \THEN
    \IF  \GtDec(\Vb,0)  \THEN
        \ifx \:FrArrowHead\empty \else  \Line(0,5)  \fi
        \Vb=\:ArrowLength;
    \ELSE
        \ifx \:FrArrowHead\empty \else   \Line(0,-5) \fi
        \Vb=-\:ArrowLength; \FI
\ELSE
    \IF  \GtDec(\Va,0)  \THEN
        \ifx \:FrArrowHead\empty \else  \Line(5,0)  \fi
        \Va=\:ArrowLength;
    \ELSE
        \ifx \:FrArrowHead\empty \else  \Line(-5,0)  \fi
        \Va=-\:ArrowLength; \FI
\FI
\Move(\Val\Va,\Val\Vb)  \MarkLoc(fr:) \:FrArrowHead  }

   \FI
   \:C=\:Ca;  \:C/10;  \:C*10;  \:Ca-\:C;
   \IF  \GtInt(\:Ca,4)  \THEN   \:Ca-5; \let\:PutTics=\:TicsByLoc  \FI
   \IF \EqInt(\:Ca,4) \THEN
  \Define\:ticLn(2){\DoLine(##1,##2)(7pt\du){
     \Units(1pt,1pt) \MoveF(-4) \LineF(2) }}
  \:Ca-1;
\ELSE  \let\:ticLn=\Line  \FI
   \IF  \LtInt(\:Ca,3)  \THEN
   \IF  \GtDec(\:Vc,0)  \THEN \:Vc= \:TicLn;
   \ELSE                     \:Vc=-\:TicLn;
\FI\FI
\ifcase \Val\:Ca
       \:Vc=   0;   \def\:TxtPos{3}  
   \or \:Vc=-\:Vc;   \:Vd=\:TicLn; \:Vd+3; \edef\:TxtPos{\Val\:Vd} 
\else              \def\:TxtPos{3} \fi
   \:Ca=\:C;  \:Ca/10;
   \Table\:Temp{#3}  \:C=-1;  \:Temp(0,99){\:Count}
   \IF  \GtInt(\:C,0) \THEN  \:PutTics \FI \RecallAll}
\Define\:PutTics{       \Va/\:C;   \Vb/\:C;
  { \:Temp(0,99){  {\PenSize(\:TicTh)\:tic} \:TicText }  }
    \IF  \GtInt(\:Ca,0)  \THEN
       \:C*\:Ca;  \Va/\:Ca;  \Vb/\:Ca;
       \IF  \GtDec(\:Vc,\:TicLn)  \THEN  \:Vc= \:TicLn;  \FI
       \IF  \LtDec(\:Vc,-\:TicLn) \THEN  \:Vc=-\:TicLn;  \FI
       \:Vc*\:TicFc;       { \let\:ticLn=\:Ln
       \Do(0,\Val\:C){  {\PenSize(\:TicTh)\:tic}
                          \Move(\Val\Va,\Val\Vb)  }}
     \FI}\Define\:TicsByLoc{ \MarkLoc(x:)
  \Define\:tempa(2){\MoveTo(0,##1) \MarkLoc(a:)
                    \MoveTo(0,##2) \MarkLoc(b:)
                    \LSeg\:Ve(a:,b:)}
  \:Temp(0,0){\:tempa }
  \Define\:tempa(2){ \MoveTo(0,##1) \MarkLoc(b:)
     \LSeg\:Vd(a:,b:)  \:Vd/\:Ve;
     { \Va*\:Vd;   \Vb*\:Vd;  \MoveToLoc(x:)
       \Move(\Val\Va,\Val\Vb)
       {\PenSize(\:TicTh)\:tic} \:TicText(##2) }}
  \:Temp(1,99){\:tempa  } }\Define\:TicText(1){          \MarkLoc(o:)
   \:MvTicTxt\TicText(--#1--) \MoveToLoc(o:)
   \Move(\Val\Va,\Val\Vb)}

\Define\TicSpec(3){
    \:edef\:TicLn{#1} \:edef\:TicFc{#2} \:edef\:TicTh{#3}}

\TicSpec(6,0.6,0.2 pt)

\TextPar\Define\TicText(1){ \IF \EqText(,#1) \THEN
                     \ELSE \Text(--\strut#1--) \FI}  
\fi
\:CheckOption{bar}\if:option
   \:DefineExt\BarChart(1){\endgroup
   \Table\:Chart{#1}  \MarkLoc(Origin)
   \XSaveUnits    \:Ca=0;  \Vb=0;  \:Vc=0;
\:Chart(0,999){\:Ca+1;\:MaxBar}
\ifx  \:Hchart\empty
   \:Ca*10;  \:Vc+\:BarGrdOver;
\ifx H\:BarDir
   \Move(\Val\:Vc,0)  \MarkLoc(NE)
   \Move(0,-\Val\:Ca) \MarkLoc(SE)
   \Vb-\:Vc;  \Move(\Val\Vb,0) \Move(-\:BarGrdUnder,0)
   \MarkLoc(SW)   \Move(0,\Val\:Ca)   \MarkLoc(NW)
   \edef\:BsLn{{\noexpand\Line(0,-10)}}
\else
   \Move(0,\Val\:Vc)  \MarkLoc(NW)
   \Move(\Val\:Ca,0) \MarkLoc(NE)
   \Vb-\:Vc;  \Move(0,\Val\Vb)
   \Move(0,-\:BarGrdUnder)    \MarkLoc(SE)
   \Move(-\Val\:Ca,0) \MarkLoc(SW)
   \edef\:BsLn{{\noexpand\Line(10,0)}}
\fi  \MoveToLoc(Origin) 
\else
   \Va=\:Vc; \:Vd=\:Vc;  \:Vd-\Vb;  \Va/\:Vd;
\ifx H\:BarDir
   \Va*\:Hchart;   \Move(\Val\Va,0)
   \Move(\:BarGrdOver,0)  \MarkLoc(NE)
   \Move(0,-\:Vchart) \MarkLoc(SE)
   \Move(-\:Hchart,0) \Move(-\:BarGrdOver,0)
       \Move(-\:BarGrdUnder,0) \MarkLoc(SW)
   \Move(0,\:Vchart) \MarkLoc(NW)
   \:Ca*10;  \edef\:BsLn{{\noexpand\Line(0,-10)}}
\else
   \Va*\:Vchart;   \Move(0,\Val\Va)
   \Move(0,\:BarGrdOver)  \MarkLoc(NW)
   \Move(\:Hchart,0) \MarkLoc(NE)
   \Move(0,-\:Vchart) \Move(0,-\:BarGrdOver)
       \Move(0,-\:BarGrdUnder) \MarkLoc(SE)
   \Move(-\:Hchart,0) \MarkLoc(SW)
   \:Ca*10;  \edef\:BsLn{{\noexpand\Line(10,0)}}
\fi  \MoveToLoc(Origin) 
   \Va=\:Hchart;   \Va/\:Ca;
   \:Vc-\Vb;  \:Vd=\:Vchart;  \:Vd/\:Vc;
   \:Hflip\Units(\Val\Va pt,\Val\:Vd pt)
   \:Vc+\Vb;
\fi   
   \def\:GrdLn{}
\ifx  \:BarGrdDist\empty   \def\:BsLn{}    \else
    \let\du=\relax
   \IF  \GtInt(\:BarGrdDist,0)  \THEN
      \:Vd=\:Vc;  \:Vd/\:BarGrdDist;  \:C[\:Vd];
      \IF \GtDec(\:C,\:Vd) \THEN \:C-1; \FI
      \IF \GtInt(\:C,0)  \THEN
         \edef\:GrdLn{{  { \noexpand\noexpand\noexpand\:BsLn }
   \PenSize(\:TicTh)
   \noexpand\noexpand\noexpand\Do(1,\Val\:C){
   \:Hflip{\noexpand\noexpand\noexpand\Move}(0,\:BarGrdDist)
   { \:HflipMY{\noexpand\noexpand\noexpand\Line}(10,0) }  }}}
\:GridLoc(Max)

      \FI
      \Vb+\:Vc;   \Vb-\:BarGrdUnder;
      \Vb/\:BarGrdDist;  \:C[\Vb];
      \IF \LtDec(\:C,\:Vd) \THEN \:C+1; \FI
      \IF \LtInt(\:C,0)  \THEN
             \edef\:GrdLn{
   \ifx\:AddBarDepth\:AddBarDp\:HflipMY{\noexpand\Move}(-\:HBarx,\:VBarx)  \fi
    { \noexpand\:BsLn }
    \:GrdLn { \PenSize(\:TicTh)
        \noexpand\Do(-1,\Val\:C){  \:Hflip{\noexpand\Move}(0,-\:BarGrdDist)
        { \:HflipMY{\noexpand\Line}(10,0) }  }}
}
\:GridLoc(Min)

      \ELSE  \edef\:GrdLn{{\ifx\:AddBarDepth\:AddBarDp
            \:HflipMY{\noexpand\Move}(-\:HBarx,\:VBarx)  \fi

            {\noexpand\:BsLn}   \:GrdLn}}
        \FI
      \MoveToLoc(Origin)
   \FI
   \let\du=\:SvDu
\fi 
   \Va=0;  \:Vc=0;
   \:Chart(0,99){
      \:Hflip\Move(0,\Val\:Vc) 
      \:NoBarClip()      \:DrawBars}
   \ifx\:AddBarDepth\:AddBarDp
   \MoveToLoc(Origin) \:HflipMY\Move(-\:HBarx,\:VBarx)  \MarkLoc(Origin)
   \MoveToLoc(NW) \:HflipMY\Move(-\:HBarx,\:VBarx)  \MarkLoc(NW)
   \MoveToLoc(NE) \:HflipMY\Move(-\:HBarx,\:VBarx)  \MarkLoc(NE)
   \MoveToLoc(SW) \:HflipMY\Move(-\:HBarx,\:VBarx)  \MarkLoc(SW)
   \MoveToLoc(SE) \:HflipMY\Move(-\:HBarx,\:VBarx)  \MarkLoc(SE)
\fi
   \XRecallUnits   \MoveToLoc(Origin)
   \let\:Return=\:undefined   \:NoBarClip()   }  \Define\:DrawBar(2){
   \:Bars(#1,#1){\:GetXY}  \Vb-\Va;
   \ifx \:Compute\empty                       \edef\:tempa{#2}
   \else    \def\:Return##1{\edef\:tempa{\Val##1}}\:Compute(#2)
   \fi
   \ifx  H\:BarDir      \Move(0,-\Val\Va)
      \:DrawBr(\:tempa,-\Val\Vb) \:AddBarPaint(#1,\:tempa)
      \Clip(\:tempa,-\Val\Vb)
      \:AddBarDepth(#1,\:tempa) \:Cluster(\:tempa) \Move(0,\Val\Va)
   \else                \Move(\Val\Va,0)
      \:DrawBr(\Val\Vb,\:tempa) \:AddBarPaint(#1,\:tempa)
      \Clip(\Val\Vb,\:tempa)
      \:AddBarDepth(#1,\:tempa) \:Cluster(\:tempa) \Move(-\Val\Va,0)
   \fi  }

\Define\DrawBar(2){ \DrawRect(#1,#2)}

\Define\:DrawBr(1){\MarkLoc(x:)\DrawBar(#1)\MoveToLoc(x:)}
\:DefineExt\BarChartSpec(1){\endgroup
   \IF  \EqText(,#1) \THEN \ELSE \:BarSpc(#1)    \FI  \:BarCol }

\Define\:BarSpc(2){
   \let\:Temp=\:SeqBars       \let\:BarDir=V  \let\:newStck=\empty
\let\:endnewStck=\empty
\Define\:AddBarDepth(2){}  \Define\:Cluster(1){}
\:OptBars#1...//    \Table\:Bars{#2}
   \:C=-1;  \:Bars(0,999){\:Count}
   \ifx \:Temp\:TxtBar  \IF  \GtInt(\:C,1)  \THEN
      \:tk{#1,#2}  \:aldwarn8{}
   \FI \fi
  \:Temp }\def\:OptBars#1#2#3#4//{\:OpBr(#1)\:OpBr(#2)\:OpBr(#3)}

\Define\:OpBr(1){
        \IF \EqText(H,#1)  \THEN   \let\:BarDir=H
   \ELSE\IF \EqText(T,#1)  \THEN   \let\:Temp=\:TxtBar
   \ELSE\IF \EqText(S,#1)  \THEN   \Define\:Cluster(1){ \:Vc-##1;
     \:Hflip\Move(0,##1) }
\def\:newStck{
   \edef\:svMinBar{\Val\Vb}   \Vb=0;
   \edef\:svMaxBar{\Val\:Vc}   \:Vc=0;  }
\def\:endnewStck  {
   \IF \LtDec(\:svMinBar,\Vb)  \THEN  \Vb=\:svMinBar;  \FI
   \IF \GtDec(\:svMaxBar,\:Vc)  \THEN  \:Vc=\:svMaxBar;  \FI }

   \ELSE\IF \EqText(3,#1)  \THEN   \let\:AddBarDepth=\:AddBarDp
   \FI\FI\FI\FI  }  \def\:BarsCommand#1{\EqInt(\Val\:C,#1)\THEN
       \Define\:DrawBars(}

\Define\:SeqBars{
   \IF \:BarsCommand01){ \:DrawBar(0,##1) \:GrdLn
           \:HflipMY\Move(10,0)  \:Vc=0;}
        \Define\:MaxBar{\:maxBar}

   \ELSE  \IF \:BarsCommand12){  \:DrawBar(0,##1)  \:DrawBar(1,##2)
           \:slotgrdln
 }
        \Define\:MaxBar(2){
   \:newStck   \:maxBar(##1)\:maxBar(##2)   \:endnewStck   }

   \ELSE  \IF \:BarsCommand23){ \:DrawBar(0,##1)
           \:DrawBar(1,##2)   \:DrawBar(2,##3)
           \:slotgrdln
 }
        \Define\:MaxBar(3){
   \:newStck   \:maxBar(##1)\:maxBar(##2)\:maxBar(##3)   \:endnewStck   }

   \ELSE  \IF \:BarsCommand34){ \:DrawBar(0,##1)
           \:DrawBar(1,##2)  \:DrawBar(2,##3)      \:DrawBar(3,##4)
           \:slotgrdln
 }
        \Define\:MaxBar(4){   \:newStck
   \:maxBar(##1)\:maxBar(##2)\:maxBar(##3)\:maxBar(##4)   \:endnewStck  }

   \ELSE  \IF \:BarsCommand45){ \:DrawBar(0,##1)
           \:DrawBar(1,##2)  \:DrawBar(2,##3)
           \:DrawBar(3,##4)  \:DrawBar(4,##5)
           \:slotgrdln
 }
        \Define\:MaxBar(5){   \:newStck   \:maxBar(##1)
   \:maxBar(##2)\:maxBar(##3)\:maxBar(##4)\:maxBar(##5) \:endnewStck  }

   \ELSE \:aldwarn7{}
   \FI  \FI  \FI  \FI  \FI  }\Define\:GOBBLE(1){}
\Define\:AddBarDp(2){       \MarkLoc(a:)
   \:HflipMY\Line(-\:HBarx,\:VBarx)  \MarkLoc(b:)
   \:HflipMY\Line(0,#2)         \MarkLoc(c:)
   \:HflipMY\Line(\:HBarx,-\:VBarx)  \MarkLoc(d:)
   \ifx \:Cluster\:GOBBLE  \:TopIIIdTop
   \else
      \IF \GtInt(\:C,#1) \THEN\ELSE \:TopIIIdTop \FI
   \fi
   \:HflipMY\Move(\:HBarx,0)    \MarkLoc(f:) \MoveToLoc(b:)
   \:HflipMY\Move(0,-\:VBarx) \MarkLoc(b:) \CSeg\Clip(b:,f:)
   \MoveToLoc(a:) }

\Define\BarDepth(1){
   \edef\:HBarx{#1pt\noexpand\du}  \Va=#1;  \Va*1.25;
   \edef\:VBarx{\Val\Va pt\noexpand\du} }

\BarDepth(3)

\def\:Hflip#1(#2,#3){%
   \ifx  H\:BarDir   #1(#3,#2)\else  #1(#2,#3)\fi}
\def\:HflipMY#1(#2,#3){%
   \ifx  H\:BarDir   #1(#3,-#2)\else  #1(#2,#3)\fi}  \Define\:TopIIIdTop{
   \:HflipMY\Line(\Val\Vb,0)  \MarkLoc(e:)
   \:HflipMY\Line(-\:HBarx,\:VBarx) \MarkLoc(f:)
   {\LineToLoc(c:)}}\def\::Compute#1{\ifx \:Compute\empty
         \edef\:tempa{#1}%
   \else \def\:Return##1{\edef\:tempa{\Val##1}}\:Compute(#1)\fi}

\Define\:TxtBar{
   \Define\:DrawBars(2){
      \MarkLoc(a:)  \:Bars(0,0){\:GetXY}  \Va+\Vb;  \Va/2;
\:HflipMY\Move(\Val\Va,0)  \MarkLoc(BarBot)
\::Compute{##1}  \:Hflip\Move(0,\:tempa)         \MarkLoc(BarTop)  
      \BarText(--##2--)  \MoveToLoc(a:)
      \:DrawBar(0,##1)   \:GrdLn
      \:HflipMY\Move(10,0)  \:Vc=0;}
   \Define\:MaxBar{\:maxTxtBar}}

\Define\:maxTxtBar(2){\:maxBar(#1)}

\TextPar\Define\BarText(1){
   \:Hflip\Move(0,4pt\du)  \:Hflip\EntryExit(0,-1,0,0)
   \Text(--\ignorespaces#1--)}\:DefineExt\:BarCol(1){\endgroup
   \IF  \EqText(,#1) \THEN\ELSE
      \Table\:Color{#1&0&0&0&0} \FI}  \Define\:GetColor(1){\:Vc=#1;}

\Define\:AddBarPaint(2){{    \PenSize(0.2pt)
   \MarkLoc(a:)  \MoveTo(0,0)  \MarkLoc(b:)
\ifx  V\:BarDir    \MoveTo(\Val\Vb,#2)
\else              \MoveTo(#2,-\Val\Vb)    \fi
\MarkLoc(c:) \MoveToLoc(a:)  \Units(1pt,1pt) 
   \:Color(#1,#1){\:GetColor}
   \IF         \EqDec(\:Vc,2)  \THEN
      \:VrPnt
   \ELSE\IF    \EqDec(\:Vc,3)  \THEN
      \:HrPnt
   \ELSE\IF    \EqDec(\:Vc,4)  \THEN
      \PenSize(0.2pt)  {\:VrPnt}  \:HrPnt 
   \ELSE \IF   \GtDec(\:Vc,0)  \THEN
      \:ChooseColor(\Val\:Vc)
      \CSeg\PaintRect(b:,c:)
   \FI  \FI  \FI  \FI}}\Define\:VrPnt{   \CSeg\:GetXY(b:,c:)  \:C[\Va];
   \:Vc=\Va;  \:Vc-2;  \:Vc/1.5;  \:C[\:Vc];
   \:Vc=\:C; \:Vc*1.5;  \:Vc-\Va;  \:Vc/-2;    \Move(\Val\:Vc,0)
   \Do(0,\Val\:C){ {\Line(0,\Val\Vb)} \Move(1.5,0) } }
\Define\:HrPnt{   \CSeg\:GetXY(b:,c:)  \:C[\Vb];
   \:Vc=\Vb; \ifx  V\:BarDir \:Vc-2; \else \:Vc+2; \fi
   \:Vc/1.5; \:C[\:Vc];  \:Vc=\:C;  \:Vc*1.5;
   \:Vc-\Vb; \:Vc/-2;  \Move(0,\Val\:Vc)
   \Do(0,\Val\:C){ {\Line(\Val\Va,0)}
       \Move(0,\ifx  H\:BarDir-\fi 1.5)} }
\Define\BarClipOn{\Define\:NoBarClip(1){}}
\Define\BarClipOff{\let\:NoBarClip=\Clip}

\BarClipOff\Define\:maxBar(1){
   \ifx \:Compute\empty                       \edef\:tempa{#1}
   \else    \def\:Return##1{\edef\:tempa{\Val##1}}\:Compute(#1)
   \fi
   \IF  \GtDec(\:tempa,0)  \THEN
       \ifx \:Cluster\:GOBBLE
   \IF \GtDec(\:tempa,\:Vc)  \THEN  \:Vc=\:tempa;  \FI
\else
   \:Vc+\:tempa;
\fi

   \ELSE
       \ifx \:Cluster\:GOBBLE
   \IF \LtDec(\:tempa,\Vb)  \THEN  \Vb=\:tempa;  \FI
\else
   \Vb-\:tempa;
\fi

   \FI}
\Define\ChartSize(1){
   \IF  \EqText(,#1)  \THEN
      \def\:Hchart{}  \def\:Vchart{}
   \ELSE  \:ChartSz(#1)  \FI}

\Define\:ChartSz(2){
   \MarkLoc(o:) \Move(#1   \du,#2 \du)  \MarkLoc(x:)
   \CSeg\:GetXY(o:,x:)     \edef\:Hchart{\Val\Va}
   \edef\:Vchart{\Val\Vb}  \MoveToLoc(o:) }

\ChartSize()\let\:SvDu=\du

\Define\:GridLoc(1){
   \:Vd=\:C; \:Vd*\:BarGrdDist; \MoveToLoc(Origin)
   \ifx V\:BarDir
      \Move(10,0)  \MarkLoc(#1)
      \MoveToLL(Origin,#1)(NE,SE)  \Move(0,\Val\:Vd)
   \else
      \Move(0,-10)  \MarkLoc(#1)
      \MoveToLL(Origin,#1)(SW,SE)  \Move(\Val\:Vd,0)
   \fi
   \let\du=\:SvDu
   \ifx\:AddBarDepth\:AddBarDp  \:HflipMY\Move(-\:HBarx,\:VBarx)  \fi
   \let\du=\relax   \MarkLoc(#1)}
\Define\BarGrid(3){
   \def\:BarGrdDist{#1}  \def\:BarGrdUnder{#2}
   \def\:BarGrdOver{#3}}
\BarChartSpec(V,2.5,7.5)(0)\BarGrid(0,0,0)
\Define\:slotgrdln{
   \:HflipMY\Move(0,\Val\:Vc) \:GrdLn
   \:HflipMY\Move(10,0)  \:Vc=0;         }

\fi
\:CheckOption{diagram}\if:option
   \edef\AlDraTex{spread,grid,tree,\AlDraTex}
\else \:CheckOption{spread}\if:option
   \edef\AlDraTex{diagram,\AlDraTex}
\else \:CheckOption{grid}\if:option
   \edef\AlDraTex{diagram,\AlDraTex}
\else \:CheckOption{tree}\if:option
   \edef\AlDraTex{diagram,\AlDraTex}
\fi\fi\fi\fi
\:CheckOption{xy}     \if:option  \:circletrue  \fi
\:CheckOption{diagram}\if:option  \:circletrue  \fi
\if:circle     \Define\ArrowHeads(1){
   \Define\:FrArrowHead{}
   \Define\:BkArrowHead{}
   \Define\:FrCrvArrowHead{}
   \Define\:BkCrvArrowHead{}
   \IF  \GtInt(#1,0) \THEN
      \Define\:FrArrowHead{\:ArrowHead(bk:,fr:)}
      \Define\:FrCrvArrowHead{
         \DSeg\RotateTo(2:,2':)
         \MoveToLoc(2:)  \MoveF(\:ArrowLength pt\du)
         \MoveF(0.1 pt\du) \MarkLoc(1'':)   \MarkLoc(2'':)
         \MoveToLoc(2':) \MoveF(\:ArrowLength pt\du)
         \MoveF(0.1 pt\du) \MarkLoc(1'':)  \MarkLoc(2':)
         \:ArrowHead(2'':,2:) }
   \FI
   \IF  \GtInt(#1,1) \THEN
      \Define\:BkArrowHead{\:ArrowHead(fr:,bk:)}
      \Define\:BkCrvArrowHead{
         \DSeg\RotateTo(1:,1':)
         \MoveToLoc(1:)  \MoveF(\:ArrowLength pt\du)
         \MoveF(0.1 pt\du) \MarkLoc(1'':)  \MarkLoc(1'':)
         \MoveToLoc(1':) \MoveF(\:ArrowLength pt\du)
         \MoveF(0.1 pt\du) \MarkLoc(1'':)  \MarkLoc(1':)
         \:ArrowHead(1'':,1:) }
   \FI }

\Define\:ArrowHead(2){
   \LSeg\Va(#1,#2)
   \IF  \LtDec(\Va,\:ArrowLength)  \THEN      \:aldwarn0{}
   \ELSE   \MoveToLoc(#2)  \DSeg\RotateTo(#1,#2)
      \MoveF(-\:ArrowLength pt\du)
      { \Ragged(\:ArrowRagged)  \Rotate(90)  \Va=\:ArrowWidth;
        \Va/2;     \MoveF(\Val\Va pt\du)  \MarkLoc(#2')
        \Va*2;  \MoveF(-\Val\Va pt\du)  \MarkLoc(#2'')
        \:head(#2)   }   \MarkLoc(#2)                    \FI   }

\Define\ArrowSpec(1){\:ArSp(#1,)}

\Define\:ArSp(2){
  \IF \EqText(V,#1) \THEN \Define\:head(1){
   { \Va/2;  \MoveF(\Val\Va pt\du) \MarkLoc(a:)
     \CSeg\:EdgeLine(a:,##1)  }
    \MarkLoc(a:) \CSeg\:EdgeLine(a:,##1)
    \MarkLoc(a:) \CSeg\:EdgeLine(a:,##1') }        \FI
  \IF \EqText(H,#1) \THEN \Define\:head(1){
   \MarkLoc(a:) \CSeg\:EdgeLine(a:,##1)
   \MarkLoc(a:) \CSeg\:EdgeLine(a:,##1')
   \MarkLoc(a:) \CSeg\:EdgeLine(a:,##1'')}  \FI
  \IF \EqText(F,#1) \THEN \Define\:head(1){\PaintQuad(##1,##1,##1',##1'')}
  \FI
  \IF \EqText(,#2) \THEN\ELSE  \:ArrowHeadSpec(#2) \FI }

\Define\:ArrowHeadSpec(4){  \:edef\:ArrowLength{#1}
   \:edef\:ArrowWidth{#2}  \:edef\:ArrowRagged{#3}}
\ArrowHeads(0)   \ArrowSpec(F,10,6,5)     \:circlefalse \fi
\:CheckOption{diagram}\if:option
      \Define\NewNode(2){
   \Define#1(1){
      \Define\:FrameType{\csname \string#1:Frame\endcsname}
      \Define\:NodeName{##1}
      \Indirect\Define<##1.mvto>{#2}  \:NodeBody}
   \expandafter\Define\csname \string#1:Frame\endcsname}
\Define\TextNode{
    \TextPar\Define\:NodeBody(1){
       \Object\:tmp{
          \:NodeContent(--##1--)
          \:FrameType }
       \:SetNode }
    \TextPar\Define\:NodeContent}

\Define\PictNode{
    \Define\:NodeBody(1){
       \Object\:tmp{ \Object\:temp{\:NodeContent(##1)}
              \:temp \:FrameType }
       \:SetNode }
    \TextPar\Define\:NodeContent}

\TextNode(1){\Text(--#1--)}\Define\:SetNode{\:SetNodeA }

\Define\:SetNodeA{  \:tmp
   \MarkLoc(a:)  \MoveToExit(0,0)  \MarkLoc(\:NodeName)
                 \MoveToExit(1,1)  \MarkLoc(\:NodeName;:11)
   \MoveToLoc(a:)   \:circlefalse }
\Define\DefNode{
   \Define\:SetNode{
      \Indirect\let<\:NodeName :call> = \:tmp
\def\:temp{\Indirect\let<\:NodeName :bdy> = }
   \expandafter\:temp \csname \string \:tmp .\endcsname
\Indirect\let<\:NodeName :frm> = \:FrameType
\Indirect\let<\:NodeName :cnt> = \:NodeContent
\Indirect\edef<\:NodeName :crc>{
   \if:circle \noexpand\:circletrue
   \else      \noexpand\:circlefalse \fi}

      \FigSize\Va\Vb{ \:SetNodeA }   \Va/2; \Vb/2;
\MarkLoc(\:NodeName)      \Move(\Val\Va,\Val\Vb)
\MarkLoc(\:NodeName;:11)  \MoveToLoc(\:NodeName)

      \:circlefalse
      \Define\:SetNode{\:SetNodeA }  }
   \:DefNode
}
\def\:DefNode#1(#2){#1(#2)}
\Define\PutNode(1){
   \def\:NodeName{#1}
   \Indirect{ \expandafter \let \csname \string \:tmp .\endcsname = }
               <#1:bdy>
   \Indirect<#1:crc>
   \Indirect{ \let\:FrameType   = }<#1:frm>
   \Indirect{ \let\:NodeContent = }<#1:cnt>
   \Indirect{ \let\:tmp        = }<#1:call>
   \ifx\:tmp\relax  \:aldwarn{10}{#1}\fi
   \:SetNodeA
   \Indirect\let<\:NodeName :call> = \:undef
\Indirect\let<\:NodeName :cnt>  = \:undef
\Indirect\let<\:NodeName :frm>  = \:undef
\Indirect\let<\:NodeName :bdy>  = \:undef
\Indirect\let<\:NodeName :crc>  = \:undef
   }
\Define\ZeroNodesDim{
   \def\MaxX{0}    \def\MaxY{0}
   \def\WidthX{0}  \def\WidthY{0} }
\Define\AddNodeDim(1){
   \CSeg\:GetXY(#1,#1;:11)
   \IF \GtDec(\Va,\MaxX) \THEN  \edef\MaxX{\Val\Va} \FI
   \Va*2;      \Va+\WidthX;   \edef\WidthX{\Val\Va}
   \IF \GtDec(\Vb,\MaxY) \THEN  \edef\MaxY{\Val\Vb} \FI
   \Vb*2;     \Vb+\WidthY;   \edef\WidthY{\Val\Vb}  }

\Define\NewCIRCNode(3){
   \ifx \CIRC\:undefined                  \:aldwarn9\relax
      \font\CIRC=lcircle10\space scaled\magstep5
   \fi
   \Define\:Temp(2){
      \IF  \EqText(,#3)  \THEN
         \NewNode(#1,\MoveToOval){
            \Text(--\hbox to##1pt{%
               \vrule height ##2pt depth ##2pt%
                           width0pt%
            \hss\hbox to##2pt{%
               \CIRC \char #2\hss}}--)}
       \ELSE
          \NewNode(#1,\MoveToOval){
             \Text(--\hbox to##1pt{%
                \hss\hbox to##2pt{%
                   \CIRC \char #3\hss}}--)
             \Text(--\hbox to##1pt{%
                \vrule height ##2pt depth ##2pt%
                           width0pt%
             \hss\hbox to##2pt{%
                \CIRC \char #2\hss}}--)}
        \FI  }
   \setbox\:box=\hbox{\CIRC \char #2}
   \IF \EqText(,#3) \ELSE \IF  \GtInt(#3,#2)  \THEN
       \setbox\:box=\hbox{\CIRC \char #3}
   \FI  \FI
   \edef\:tempa{\noexpand\Va=\:Cons\wd\:box;}
   \:tempa  \Vb=\Va; \Vb/2;
   \edef\:tempa{\noexpand\:Temp(\Val\Va,\Val\Vb)}
   \:tempa  }
\Define\NodeMargin(2){\def\:XNodeMargin{#1}\def\:YNodeMargin{#2}}
\NodeMargin(2,2) \Define\GetNodeSize{
    \MoveToExit(1,1)
    \Move(\:XNodeMargin pt\du,\:YNodeMargin  pt\du)
    \MarkLoc(:11)    \MoveToExit(0,0)   \MarkLoc(:00)
    \CSeg\:GetXY(:00,:11)}

\Define\SetMinNodeSize{
    \IF \LtDec(\Va,\:NodeX)  \THEN  \Va=\:NodeX;  \FI
    \IF \LtDec(\Vb,\:NodeY)  \THEN  \Vb=\:NodeY;  \FI }

\Define\MinNodeSize(1){\:MinNodeSize(#1,)}
\Define\:MinNodeSize(2){
   \IF  \EqText(,#2)  \THEN
      \edef\:Temp{\noexpand\:MinNdSz(\csname :#1::\endcsname,)}
   \ELSE
      \def\:Temp{\:MinNdSz(#1,#2)}
   \FI  \:Temp}
\Define\:MinNdSz(3){
    \Va=#1;  \Va/2;  \edef\:NodeX{\Val\Va}
    \Va=#2;  \Va/2;  \edef\:NodeY{\Val\Va}  }

\Define\SaveNodeSize(1){ \CSeg\:GetXY(#1,#1;:11)
   \Va*2;  \Vb*2;
   \expandafter\edef\csname :#1::\endcsname{\Val\Va,\Val\Vb}  }

               \Define\:Temp{
\def\SaveDrawSize(##1)##2{
   \FigSize\Va\Vb{##2}
   \expandafter\edef\csname :##1::\endcsname{\Val\Va,\Val\Vb}  }
                      }\:Temp

\MinNodeSize(0,0)\def\:NodeLine{\Line}   \def\:NodeArc{\DrawOvalArc}
\Define\SRectNodeSpec(1){\:edef\:ShadowSize{#1}}
\Define\VRectNodeSpec(1){\:edef\:Vdepth{#1}}

\VRectNodeSpec(4)  \SRectNodeSpec(3) 
\NewNode(\RectNode,\MoveToRect){   \:RectShape }
\NewNode(\VRectNode,\:MoveToVRect){ \:RectShape
   \:Vc=\:Vdepth;  \:Vc/2;  \:NodeLine(\:Vdepth,-\Val\:Vc)
   \:NodeLine(\Val\Va,0)    {\:NodeLine(-\:Vdepth,\Val\:Vc)}
   \:NodeLine(0,\Val\Vb)     \:NodeLine(-\:Vdepth,\Val\:Vc) }
\NewNode(\SRectNode,\MoveToRect){  \Units(1pt,1pt)
   \GetNodeSize   \SetMinNodeSize
   \Move(-\Val\Va,-\Val\Vb)   \Va*2;   \Vb*2;
   \:NodeLine(0,\Val\Vb)  \:NodeLine(\Val\Va,0)
   \PenSize(\:ShadowSize pt)
   \:Ve=\:ShadowSize;   \:Ve/2; \Va+\:ShadowSize;
   \Move(\Val\:Ve,0)   \Line(0,-\Val\Vb)
   \Move(\Val\:Ve,-\Val\:Ve)  \Line(-\Val\Va,0) }

\Define\:RectShape{  \Units(1pt,1pt)
   \GetNodeSize   \SetMinNodeSize
   \Move(-\Val\Va,-\Val\Vb)   \Va*2;   \Vb*2;
   \:NodeLine(0,\Val\Vb)  \:NodeLine(\Val\Va,0)
   \:NodeLine(0,-\Val\Vb) \:NodeLine(-\Val\Va,0) }

\Define\MoveToRect(3){
    \CSeg[#2]\:GetX(#1,#1;:11)
    \CSeg[#3]\:GetY(#1,#1;:11)    \MoveToLoc(#1)
    \Move(\Val\Va pt\du,\Val\Vb pt\du)}
\NewNode(\RRectNode,\:MoveToRRect){    \Units(1pt,1pt)
   \GetNodeSize   \SetMinNodeSize
   \Va+3;   \Move(-\Val\Va,-\Val\Vb)   \Va*2;   \Vb*2;
   \:Vc=\Va;  \:Vc/4;   \:Vd=\Vb;  \:Vd/4;
\IF \GtDec(\:Vc,9) \THEN   \:Vc=9; \FI
\IF \GtDec(\:Vd,9) \THEN   \:Vd=9; \FI
\Va-\:Vc;   \Va-\:Vc;   \Vb-\:Vd;  \Vb-\:Vd;
   { \Move       (\Val\:Vc,\Val\:Vd)
     \:NodeArc(\Val\:Vc,\Val\:Vd)(180,225)  }
   \Move(0,\Val\:Vd)  \:NodeLine(0,\Val\Vb)    \Move(\Val\:Vc,0)
   \:NodeArc(\Val\:Vc,\Val\:Vd)(90,180)
   \Move(0,\Val\:Vd)  \:NodeLine(\Val\Va,0)    \Move(0,-\Val\:Vd)
   \:NodeArc(\Val\:Vc,\Val\:Vd)(45,90)
      \:NodeArc(\Val\:Vc,\Val\:Vd)(0,45)
   \Move(\Val\:Vc,0)  \:NodeLine(0,-\Val\Vb)  \Move(-\Val\:Vc,0)
   \:NodeArc(\Val\:Vc,\Val\:Vd)(270,360)
   \Move(0,-\Val\:Vd)  \:NodeLine(-\Val\Va,0)   \Move(0,\Val\:Vd)
   \:NodeArc(\Val\:Vc,\Val\:Vd)(225,270)
}\NewNode(\ORectNode,\:MoveToORect){  \Units(1pt,1pt)
   \GetNodeSize   \SetMinNodeSize   \Va-2;
   { \Move(\Val\Va,0)   \:NodeArc(9,\Val\Vb)(-90,90) }
   { \Move(-\Val\Va,0)  \:NodeArc(9,\Val\Vb)(90,-90) }
   \Move(-\Val\Va,\Val\Vb)   \Va*2;   \Vb*2;
   \:NodeLine(\Val\Va,0)
   { \Move(0,-\Val\Vb)   \:NodeLine(-\Val\Va,0)} }
\NewNode(\DRectNode,\:MoveToDRect){ \Units(1pt,1pt)
   \GetNodeSize   \SetMinNodeSize
   { \Move(\Val\Va,0) \Move(9,0)
     {\:NodeLine(-9,\Val\Vb)}
     {\:NodeLine(-9,-\Val\Vb)} }
   { \Move(-\Val\Va,0)  \Move(-9,0) {\:NodeLine(9,\Val\Vb)}
     {\:NodeLine(9,-\Val\Vb)} }
   \Move(\Val\Va,\Val\Vb)   \Va*2;
   \:NodeLine(-\Val\Va,0) { \Move(0,-\Val\Vb)
      \Move(0,-\Val\Vb)
   \:NodeLine(\Val\Va,0)} }\Define\CircleNode{ \:circletrue  \OvalNode}

\NewNode(\OvalNode,\MoveToOval){
   \Units(1pt,1pt)    \GetNodeSize     \Va/0.707;  \Vb/0.707;
   \SetMinNodeSize
   \if:circle
      \IF        \LtDec(\Va,\Vb)  \THEN \Va=\Vb;
      \ELSE  \IF \LtDec(\Vb,\Va)  \THEN \Vb=\Va;
      \FI \FI
   \fi        \:NodeArc(\Val\Va,\Val\Vb) (0,360) }
\NewNode(\DiamondNode,\:MoveToDiamond){
    \Units(1pt,1pt)      \GetNodeSize
    \Vb*1.66;  \Va+\Vb;  \Vb=\Va; \Vb*0.66; \SetMinNodeSize
    \Move(0,\Val\Vb)  \:NodeLine(-\Val\Va,-\Val\Vb)
    \:NodeLine(\Val\Va,-\Val\Vb)
    \:NodeLine(\Val\Va,\Val\Vb)
    \:NodeLine(-\Val\Va,\Val\Vb)}\NewNode(\Node,\MoveToRect){ \:NoShape  }

\Define\:NoShape{
   \GetNodeSize   \SetMinNodeSize
   \MoveToExit(0,0) \EntryExit(0,0,0,0)   \Va*2; \Vb*2;
   \Text(--\vbox to\Val\Vb pt{\hsize=\Val\Va pt
           \leavevmode\hfil\vfil}--)}\Define\FcNode(1){
   \edef\:Temp{\noexpand\Indirect\let<#1.mvto>=\noexpand\MoveToRect}
   \:Temp  \MarkLoc(#1)  \Move(0.01pt\du,0.01pt\du)
   \MarkLoc(#1;:11) \MoveToLoc(#1) }
\Define\MoveToNodeDir(1){\edef\:Temp{\noexpand\::MvToExit(#1)}\:Temp}

\Define\::MvToExit(2){
   \DSeg\:Vc(#1,#2)  \let\:Temp=\relax
   \IF          \EqDec(\:Vc,0)  \THEN  \MoveToNode(#1,1,0)
   \ELSE   \IF  \EqDec(\:Vc,90) \THEN  \MoveToNode(#1,0,1)
   \ELSE   \IF  \EqDec(\:Vc,180)\THEN  \MoveToNode(#1,-1,0)
   \ELSE   \IF  \EqDec(\:Vc,270)\THEN  \MoveToNode(#1,0,-1)
   \ELSE   \DSeg\Va(#1,#1;:11)
\IF  \LtDec(\:Vc,\Va)  \THEN         \:PreSrchExit(#1,1,\Val\:Ve,0,1)
\ELSE  \Vb=180;  \Vb-\Va;
   \IF  \LtDec(\:Vc,\Vb)  \THEN      \:PreSrchExit(#1,\Val\:Ve,1,1,-1)
   \ELSE  \Vb=180;  \Vb+\Va;
      \IF  \LtDec(\:Vc,\Vb)  \THEN   \:PreSrchExit(#1,-1,\Val\:Ve,1,-1)
      \ELSE  \Vb=360;  \Vb-\Va;
         \IF  \LtDec(\:Vc,\Vb)  \THEN\:PreSrchExit(#1,\Val\:Ve,-1,-1,1)
         \ELSE                      \:PreSrchExit(#1,1,\Val\:Ve,-1,0)
\FI\FI\FI\FI  \FI\FI\FI\FI   \:Temp }\Define\:PreSrchExit(5){
   \def\:Temp{\def\:tempa{(#1,#2,#3}\:SearchExit(#1,#4,#5)} }
\Define\:SearchExit(3){
   \:Ve=#3;  \:Ve-#2;  \IF \LtDec(\:Ve,0)  \THEN \:Ve=-\:Ve;  \FI
   \IF  \LtDec(\:Ve,0.01)  \THEN \let\:Temp=\relax
   \ELSE
      \:Ve=#2;  \:Ve+#3;  \:Ve/2;  \edef\:tempb{\Val\:Ve}
      \expandafter\MoveToNode\:tempa)   \MarkLoc(x:)
      \DSeg\:Vd(#1,x:)
      \IF  \GtDec(\:Vd,\:Vc) \THEN
         \edef\:Temp{\noexpand\:SearchExit(#1,#2,\:tempb)}
      \ELSE
         \edef\:Temp{\noexpand\:SearchExit(#1,\:tempb,#3)}
    \FI \FI  \:Temp  }\Define\TraceDiagramOn{\def\:TrcDiag##1{##1}}
\Define\TraceDiagramOff{\def\:TrcDiag##1{}}
\TraceDiagramOff \Define\TagNode(1){{\:TagNode(#1)}}

\Define\::TagNode(1){ \let&=\relax
   \xdef\:tgnd{\:tgnd&#1,\Val\:C..\Val\:Ca}}

\Define\:bgDiTags{
   \Define\:TagNode(1){}
   \Define\:EdgeNode(1){}
   \let\:svTgNd=\:tgnd  \gdef\:tgnd{}
   \let\:svEdNd=\:ednd  \gdef\:ednd{} }

\Define\:endDiTags{
   \MarkLoc(o:)
   \ifx\:tgnd\empty \else
      \Table\:Temp{\:tgnd}
      \:Temp(1,999){\:SetNdTg} \fi
   \global\let\:tgnd=\:svTgNd
      \ifx\:ednd\empty \else
      \Table\:Temp{\:ednd}
      \:C=0;
      \:Temp(1,999){\::CurrEdge} \fi
   \global\let\:ednd=\:svEdNd

   \MoveToLoc(o:)}

\Define\:SetNdTg(2){
   \Indirect{ \Indirect\let<#1.mvto>= }<#2.mvto>
   \MoveToLoc(#2)  \MarkLoc(#1)
   \MoveToLoc(#2;:11)  \MarkLoc(#1;:11)  }

   \fi
   \:CheckOption{diagram}\if:option
      \Define\MoveToNode(3){
   \XSaveUnits \Units(1pt,1pt)
   \edef\:Temp{\noexpand\Indirect<#1.mvto>(#1,#2,#3)}
   \:Temp  \XRecallUnits }\Define\:MoveToVRect(3){
   \MoveToRect(#1,#2,#3) \:Ve=\:Vdepth;  \:Ve/2;
   \:Vd=\:Ve;  \:Vd/2;
   \IF \GtDec(#2,0.98) \THEN      \Move(-\Val\:Ve,0)
      \IF \GtDec(#3,0.98) \THEN   \Move(0,-\Val\:Vd)
   \FI \FI
   \IF \LtDec(#3,-0.98) \THEN     \Move(0,\Val\:Vd)
      \IF \LtDec(#2,-0.98) \THEN  \Move(\Val\:Ve,0)
   \FI \FI}  \Define\:MoveToRRect(3){   \MoveToRect(#1,#2,#3)
   \Va=#2;  \IF \LtDec(\Va,0)  \THEN \Va=-\Va;  \FI
   \Vb=#3;  \IF \LtDec(\Vb,0)  \THEN \Vb=-\Vb;  \FI
   \IF  \GtDec(\Va,1)  \THEN \Va=1;  \FI
   \IF  \GtDec(\Vb,1)  \THEN \Vb=1;  \FI
   \IF        \LtDec(\Va,0.9) \THEN \Va=0; \Vb=0;
   \ELSE \IF  \LtDec(\Vb,0.9) \THEN \Va=0; \Vb=0;
   \ELSE \Va-1;  \Va=-\Va;  \Vb-1;  \Vb=-\Vb;  \FI \FI
   \IF \GtDec(#2,0) \THEN  \Va=-\Va;  \FI
   \IF \GtDec(#3,0) \THEN  \Vb=-\Vb;  \FI
   \Va*30;  \Vb*30;
   \Move(\Val\Va pt\du,\Val\Vb pt\du) }\Define\:MoveToORect(3){
   \MoveToRect(#1,0,0)   \MarkLoc(o:mv)
   \MoveToRect(#1,#2,#3) \MarkLoc(x:mv)  \MarkLoc(y:mv)
   \CSeg\:GetX(o:mv,x:mv)
   \MoveToRect(#1,1,0)  \MarkLoc(a:mv)   \LSeg\:Ve(a:mv,x:mv)
\IF \LtDec(\:Ve,0.1) \THEN  \Move(9pt\du,0)
\ELSE   \DSeg\:Ve(a:mv,x:mv)
   \IF \LtDec(\:Ve,90) \THEN
      \MoveToRect(#1,1,1)  \MarkLoc(1:mv)
      \DSeg\:Vd(o:mv,x:mv)      \DSeg\:Ve(o:mv,1:mv)
      \:Vd/\:Ve; \:Vd*90;  \RotateTo(\Val\:Vd)  \MoveToLoc(a:mv)
      \MoveFToOval(9 pt\du,\Val\Vb pt\du)   \MarkLoc(y:mv)
   \ELSE \IF \GtDec(\:Ve,270)\THEN
      \MoveToRect(#1,1,-1)  \MarkLoc(1:mv)
      \DSeg\:Vd(o:mv,x:mv)      \DSeg\:Ve(o:mv,1:mv)
      \:Vd-360;   \:Ve-360;   \:Vd/\:Ve;
      \:Vd*90;  \RotateTo(-\Val\:Vd)  \MoveToLoc(a:mv)
      \MoveFToOval(9 pt\du,\Val\Vb pt\du)   \MarkLoc(y:mv)
\FI \FI \FI
   \MoveToRect(#1,-1,0)  \MarkLoc(a:mv)    \LSeg\:Ve(a:mv,x:mv)
\IF \LtDec(\:Ve,0.1) \THEN  \Move(-9pt\du,0)
\ELSE    \DSeg\:Ve(a:mv,x:mv)
   \IF \GtDec(\:Ve,90) \THEN
      \IF \LtDec(\:Ve,180) \THEN
         \MoveToRect(#1,-1,1)  \MarkLoc(1:mv)
         \DSeg\:Vd(o:mv,x:mv)      \DSeg\:Ve(o:mv,1:mv)
         \:Vd-180;  \:Ve-180;  \:Vd/\:Ve;
         \:Vd*90;  \:Vd-180;  \:Vd=-\:Vd;
         \RotateTo(\Val\:Vd)  \MoveToLoc(a:mv)
         \MoveFToOval(9 pt\du,\Val\Vb pt\du)   \MarkLoc(y:mv)
      \ELSE\IF \LtDec(\:Ve,180) \THEN
         \MoveToRect(#1,-1,-1)  \MarkLoc(1:mv)
         \DSeg\:Vd(o:mv,x:mv)      \DSeg\:Ve(o:mv,1:mv)
         \:Vd-180;  \:Ve-180;  \:Vd/\:Ve;
         \:Vd*90;  \:Vd+180;  \:Vd=-\:Vd;
         \RotateTo(\Val\:Vd)  \MoveToLoc(a:mv)
         \MoveFToOval(9 pt\du,\Val\Vb pt\du)   \MarkLoc(y:mv)
\FI \FI \FI \FI  \MoveToLoc(y:mv) }\Define\:MoveToDRect(3){ \XSaveUnits
   \MoveToRect(#1,0,0)   \MarkLoc(o:mv)
   \MoveToRect(#1,#2,#3) \MarkLoc(x:mv)  \MarkLoc(y:mv)
   \CSeg\:GetX(o:mv,x:mv)
   \MoveToRect(#1,1,0)  \MarkLoc(a:mv)     \LSeg\:Ve(a:mv,x:mv)
\IF \LtDec(\:Ve,0.1) \THEN  \Move(9pt\du,0)
\ELSE     \DSeg\:Ve(a:mv,x:mv)
   \IF \LtDec(\:Ve,90) \THEN
      \Move(0,\Val\Vb)  \MarkLoc(.:)
      \Move(9,-\Val\Vb) \MarkLoc(..:)
      \MoveToLL(.:,..:) \MarkLoc(y:mv)
   \ELSE \IF \GtDec(\:Ve,270) \THEN
      \Move(0,-\Val\Vb)  \MarkLoc(.:)
      \Move(9,\Val\Vb) \MarkLoc(..:)
      \MoveToLL(.:,..:) \MarkLoc(y:mv)
\FI \FI \FI
   \MoveToRect(#1,-1,0)  \MarkLoc(a:mv)    \LSeg\:Ve(a:mv,x:mv)
\IF \LtDec(\:Ve,0.1) \THEN  \Move(-9pt\du,0)
\ELSE      \DSeg\:Ve(a:mv,x:mv)
   \IF    \GtDec(\:Ve,90) \THEN
      \IF \LtDec(\:Ve,180)\THEN
         \Move(0,\Val\Vb)  \MarkLoc(.:)
         \Move(-9,-\Val\Vb) \MarkLoc(..:)
         \MoveToLL(.:,..:) \MarkLoc(y:mv)
      \ELSE  \IF \LtDec(\:Ve,270)\THEN
         \Move(0,-\Val\Vb)  \MarkLoc(.:)
         \Move(-9,\Val\Vb) \MarkLoc(..:)
         \MoveToLL(.:,..:) \MarkLoc(y:mv)
\FI \FI \FI \FI
\XRecallUnits \MoveToLoc(y:mv) }\Define\MoveToOval(3){
   \CSeg[#2]\:GetX(#1,#1;:11)
   \CSeg[#3]\:GetY(#1,#1;:11)
   \MoveTo(0,0)   \MarkLoc(00:)
   \MoveTo(\Val\Va pt\du,\Val\Vb pt\du)  \MarkLoc(11:)
    \MoveToLoc(#1)   \LSeg\Va(00:,11:)
   \IF \GtDec(\Va,0.001) \THEN
      \DSeg\RotateTo(00:,11:)
      \CSeg\MoveFToOval(#1,#1;:11)
   \FI}\Define\:MoveToDiamond(3){
   \CSeg\:GetXY(#1,#1;:11)
   \IF \LtDec(#2,0)  \THEN  \Va=-\Va;  \FI
   \IF \LtDec(#3,0)  \THEN  \Vb=-\Vb;  \FI
   \MoveToLoc(#1)  \Move(\Val\Va pt\du,0)  \MarkLoc(.:)
   \MoveToLoc(#1)  \Move(0,\Val\Vb pt\du)  \MarkLoc(..:)
   \CSeg\:GetXY(#1,#1;:11)  \Va*#2;  \Vb*#3;
   \MoveToLoc(#1)  \Move(\Val\Va pt\du,\Val\Vb pt\du)
   \MarkLoc(..:..)  \LSeg\Va(#1,..:..)
   \IF \GtDec(\Va,0) \THEN  \MoveToLL(.:,..:)(#1,..:..)  \FI }

   \fi
   \:CheckOption{diagram}\if:option
      \Define\EdgeSpec(1){\:EdgeSpec#1....//}
\def\:EdgeSpec#1#2#3#4.//{\:EdSp(#1)\:EdSp(#2)\:EdSp(#3)}

\Define\:EdSp(1){
   \IF  \EqText(L,#1) \THEN
      \def\:EdgeLine{\Line}     \def\:EdgeCurve{\Curve}
   \FI
   \IF  \EqText(D,#1) \THEN
      \let\:EdgeLine=\:DotLine  \let\:EdgeCurve=\:DotCurve
   \FI
   \IF  \EqText(C,#1) \THEN
      \def\:LocByAddr##1(##2){ \MoveTo(##2)
                         \MarkLoc(;:) ##1(;:)}
   \FI
   \IF  \EqText(T,#1) \THEN
      \let\:LocByAddr=\relax
   \FI
   \IF  \EqText(R,#1) \THEN
      \Define\:EdgeCorner{
   \MoveToLoc(fr:)   \MarkLoc(nxt':)   \MarkLoc(fr':)
   \:SetEndCorner(fr:,bk:)      \:SetEndCorner(nxt':,nxt:)
   \:EdgeCurve(fr:,fr:1,nxt':1,nxt':) }
   \FI
   \IF  \EqText(S,#1) \THEN
      \Define\:EdgeCorner{\MoveToLoc(fr:)    \MarkLoc(nxt':)}
   \FI}

\Define\:DotLine(1){
   {  \MarkLoc(a:)  \Move(#1) \MarkLoc(b:)
      \DSeg\RotateTo(a:,b:)  \LineF(-\:thickness \du)
      \MoveToLoc(a:)  \LineF(\:thickness\du)}
    \DoLine(#1)(3pt\du){
       \Units(1pt,1pt) \MoveF(-0.5) \LineF(1) }}

\Define\:DotCurve(4){
   { \DSeg\RotateTo(#1,#2)  \MoveToLoc(#1)
     \divide\:thickness by2  \LineF(\:thickness \du)  }
   { \DSeg\RotateTo(#4,#3)  \MoveToLoc(#4)
     \divide\:thickness by2  \LineF(\:thickness \du)  }
     \DoCurve(#1,#2,#3,#4)(3pt\du){
        \Units(1pt,1pt) \MoveF(-0.5) \LineF(1) } }

\Define\EdgeAt(6){
   \MoveToNode(#1,#2,#3)  \MarkLoc(EdgeBack)
   \Table\:EdgePath{#6}
   \:EdgePath(0,0){\:GetEdgeFront}
   \MoveToNode(#4,#5,\:EdgeFront)  \:InsertEdge(#6)}\Define\Edge(2){
   \Table\:EdgePath{#1&#2}
   \:C=-1; \:EdgePath(0,99){\:Count}
   \IF \EqDec(\:C,1)  \THEN
      \MoveToNodeDir(#1,#2)    \MarkLoc(EdgeBack)
      \MoveToNodeDir(#2,#1)
   \ELSE  \:EdgePath(2,2){\:LocByAddr\MoveToLoc}  \MarkLoc(:)
\MoveToNodeDir(#1,:)  \MarkLoc(EdgeBack)
\Define\:Temp(1){\def\:Temp{##1}}
\:EdgePath(1,1){\:Temp}
\:EdgePath(\Val\:C,\Val\:C){\:LocByAddr\MoveToLoc}
\MarkLoc(EdgeFront)   \MoveToNodeDir(\:Temp,EdgeFront)     \FI
   \Table\:EdgePath{#2}    \:InsertEdge(#2) }\Define\EdgeTo(4){
   \Table\:EdgePath{#4}
   \:C=-1; \:EdgePath(0,1){\:Count}
   \IF \EqDec(\:C,0)  \THEN
      \MoveToNode(#2,#3,#4)  \MarkLoc(EdgeFront)
      \MoveToNodeDir(#1,EdgeFront)    \MarkLoc(EdgeBack)
      \MoveToLoc(EdgeFront)
   \ELSE   \:EdgePath(1,1){\:LocByAddr\MoveToLoc}  \MarkLoc(:)
\MoveToNodeDir(#1,:)    \MarkLoc(EdgeBack)
\:EdgePath(0,0){\:GetEdgeFront}
\MoveToNode(#2,#3,\:EdgeFront)    \FI
   \:InsertEdge(#4) }\Define\EdgeFrom(4){
   \Table\:EdgePath{#1&#4}
   \:C=-1; \:EdgePath(0,99){\:Count}
   \MoveToNode(#1,#2,#3)  \MarkLoc(EdgeBack)
   \IF \EqDec(\:C,1)  \THEN
      \MoveToNodeDir(#4,EdgeBack)
   \ELSE  \Define\:Temp(1){\def\:Temp{##1}}
\:EdgePath(1,1){\:Temp}
\:EdgePath(\Val\:C,\Val\:C){\:LocByAddr\MoveToLoc}
\MarkLoc(EdgeFront)      \MoveToNodeDir(\:Temp,EdgeFront)   \FI
   \Table\:EdgePath{#4}    \:InsertEdge(#4) }\Define\:InsertEdge(1){
    \MarkLoc(EdgeFront)       \:C=-1;    \:EdgePath(0,2){\:Count}
    \IF  \EqInt(\:C,0)  \THEN
       \:OneSegmentEdge
    \ELSE
       \Define\:Temp(1){
   \:LocByAddr\MoveToLoc(##1)\MarkLoc(temp:)}
\:EdgePath(1,1){\:Temp}
\MoveToLoc(EdgeBack)  \MarkLoc(bk:)
\MoveToLoc(temp:)     \MarkLoc(fr:)   \:BkArrowHead
       \IF  \GtInt(\:C,1)  \THEN
          \def\:temp##1&##2(&){\Table\:EdgePath{##2}}
         \:temp#1(&)  \:EdgePath(1,99){\:InterEdge}
       \FI
       \MoveToLoc(EdgeFront) \MarkLoc(nxt:)
\:EdgeCorner \:EdgeSegment
\MoveToLoc(nxt':) \MarkLoc(bk:)
\MoveToLoc(nxt:)  \MarkLoc(fr:)
\:FrArrowHead  \:EdgeSegment
    \FI  \def\EdgeLabel{\:PutLabel\:LineLbl}  }\Define\:OneSegmentEdge{
   \MoveToLoc(EdgeBack)   \MarkLoc(bk:)
   \MoveToLoc(EdgeFront)     \MarkLoc(fr:)
   \:BkArrowHead  \:FrArrowHead  \:EdgeSegment }\Define\:InterEdge(1){
   \:LocByAddr\MoveToLoc(#1) \MarkLoc(nxt:)
   \:EdgeCorner \:EdgeSegment
   \MoveToLoc(nxt':) \MarkLoc(bk:)
   \MoveToLoc(nxt:)  \MarkLoc(fr:)}

\Define\:EdgeSegment{\MoveToLoc(bk:) \CSeg\:EdgeLine(bk:,fr:)}
\Define\:SetEndCorner(2){ \MoveToLoc(fr':)
   \CSeg[0.5]\Move(fr':,#2)  \MarkLoc(#1)  \LSeg\Va(fr':,#1)
   \IF  \GtDec(\Val\Va,10)  \THEN
      \MoveToLoc(fr':)  \DSeg\RotateTo(fr':,#2)
      \MoveF(10)  \MarkLoc(#1)
   \FI
   \MoveToLoc(fr':) \CSeg[0.5]\Move(fr':,#1)  \MarkLoc(#1 1)}
\Define\CurvedEdgeAt(6){
    \MoveToNode(#1,#2,#3)  \MarkLoc(1:)   \MarkLoc(EdgeBack)
    \MoveToNode(#4,#5,#6)  \MarkLoc(2:)   \MarkLoc(EdgeFront)
    \:CurvedEdgeAtDir}

\Define\:CurvedEdgeAtDir(4){
    \RotateTo(#3)   \LSeg[#4]\MoveF(1:,2:)
    \MarkLoc(2':)   \MoveToLoc(1:)
    \RotateTo(#1)   \LSeg[#2]\MoveF(1:,2:)
    \MarkLoc(1':)   \:InCrvEd    }

\Define\CurvedEdge(2){
   \DSeg\RotateTo(#1,#2)   \Rotate(\:StartEdgeDir)
   \MoveToLoc(#1)   \MoveF(50)  \MarkLoc(1:)
   \MoveToNodeDir(#1,1:)    \MarkLoc(EdgeBack)
   \MarkLoc(1:)  \LSeg[\:StartEdgeTension]\MoveF(#2,#1)
   \MarkLoc(1':)
   \DSeg\RotateTo(#2,#1)   \Rotate(\:EndEdgeDir)
   \MoveToLoc(#2)   \MoveF(50)  \MarkLoc(2:)
   \MoveToNodeDir(#2,2:)    \MarkLoc(EdgeFront)
   \MarkLoc(2:)  \LSeg[\:EndEdgeTension]\MoveF(#2,#1)
   \MarkLoc(2':) \:InCrvEd }\Define\:InCrvEd{
   \:BkCrvArrowHead\:FrCrvArrowHead
    \MoveToLoc(1:)     \MarkLoc(EdgeBack')
    \MoveToLoc(1':)    \MarkLoc(EdgeBack'')
    \MoveToLoc(2':)    \MarkLoc(EdgeFront'')
    \MoveToLoc(2:)     \MarkLoc(EdgeFront')
    \:EdgeCurve(1:,1':,2':,2:)   \def\EdgeLabel{\:PutLabel\:CurveLbl} }

\Define\CurvedEdgeSpec(4){
   \:edef\:StartEdgeDir{#1}  \:edef\:StartEdgeTension{#2}
   \:edef\:EndEdgeDir{#3}    \:edef\:EndEdgeTension{#4}   }

\CurvedEdgeSpec(10,0.2,-10,0.2)\Define\HVEdge(2){  \:XYEdge(#1,#2,1,)  }
\Define\VHEdge(2){  \:XYEdge(#1,#2,,1)  }   \Define\:XYEdge(4){
   \MoveToNode(#1,0,0) \MarkLoc(1:)  \Move(#30,#40)  \MarkLoc(1':)
\MoveToNode(#2,0,0) \MarkLoc(2:)  \Move(#40,#30)  \MarkLoc(2':)
\MoveToLL(1:,1':)(2:,2':) 
   \MarkLoc(x:)  \MarkLoc(bk:)
   \CSeg\:GetExit(1:,x:)  \MoveToNode(#1,\Val\Va,\Val\Vb)
   \MarkLoc(EdgeBack)
   \CSeg\:GetExit(2:,x:)  \MoveToNode(#2,\Val\Va,\Val\Vb)
   \MarkLoc(EdgeFront)
   \LSeg\Va(#1,EdgeBack)   \LSeg\Vb(#1,x:)
\IF \GtDec(\Vb,\Va) \THEN
   \LSeg\Va(#2,EdgeFront)   \LSeg\Vb(#2,x:)
   \IF \GtDec(\Vb,\Va) \THEN   \:Vc=1;  \ELSE  \:Vc=-1;  \FI
\ELSE  \:Vc=-1;  \FI
   \IF  \GtDec(\:Vc,0)  \THEN
      \MarkLoc(fr:)  \:FrArrowHead
      \CSeg\:EdgeLine(fr:,bk:)  \MarkLoc(fr:)
      \MoveToLoc(EdgeBack) \MarkLoc(bk:) \:BkArrowHead
      \CSeg\:EdgeLine(bk:,fr:)
      \def\EdgeLabel{\:PutLabel\:XYLbl}
   \ELSE  \CSeg\:GetExit(#1,#2)
\IF  \LtDec(\Va,0)  \THEN
   \Va=-\Va;  \:Vc=-1;  \ELSE \:Vc=1;
\FI
\IF  \LtDec(\Vb,0)  \THEN
   \Vb=-\Vb;  \:Vd=-1;  \ELSE \:Vd=1;
\FI
\IF  \GtDec(\Va,\Vb)  \THEN
      \:Vd=0;  \:MarkMid(#1,#2)  \HHEdge(#1,#2,A:)
\ELSE \:Vc=0;  \:MarkMid(#1,#2)  \VVEdge(#1,#2,A:)   \FI
 \FI }\Define\:MarkMid(2){
   \MoveToNode(#1,\Val\:Vc,\Val\:Vd)   \MarkLoc(1:)
   \MoveToNode(#2,-\Val\:Vc,-\Val\:Vd) \MarkLoc(2:)
   \CSeg[0.5]\Move(2:,1:)  \MarkLoc(A:)}
\Define\VVEdge(3){  \:XXedge(#1,#2,#3,,1) }
\Define\HHEdge(3){  \:XXedge(#1,#2,#3,1,) }\Define\:XXedge(5){
   \MoveToNode(#1,0,0) \MarkLoc(1:)  \Move(#40,#50)  \MarkLoc(1':)
   \MoveToNode(#2,0,0) \MarkLoc(2:)  \Move(#40,#50)  \MarkLoc(2':)
   \MoveToLoc(#3)      \MarkLoc(3:)  \Move(#50,#40)  \MarkLoc(3':)
   \DSeg\Va(1:,3:)  \IF \GtDec(\Va,179) \THEN  \Va-180;  \FI
\DSeg\Vb(2:,3:)  \IF \GtDec(\Vb,179) \THEN  \Vb-180;  \FI
\IF  \EqText(,#4)  \THEN
   \IF       \EqDec(\Va,0) \THEN  \HVEdge(#1,#2)
   \ELSE \IF \EqDec(\Vb,0) \THEN  \VHEdge(#1,#2)
   \ELSE \::XXEdge(#1,#2)  \FI\FI
\ELSE
   \IF       \EqDec(\Va,90) \THEN  \VHEdge(#1,#2)
   \ELSE \IF \EqDec(\Vb,90) \THEN  \HVEdge(#1,#2)
   \ELSE \::XXEdge(#1,#2)  \FI\FI
\FI}

\Define\::XXEdge(2){
   \MoveToLL(2:,2':)(3:,3':)  \MarkLoc(y:)
   \MoveToLL(1:,1':)(3:,3':)  \MarkLoc(x:)  \MarkLoc(fr:)
   \CSeg\:GetExit(1:,x:)  \MoveToNode(#1,\Val\Va,\Val\Vb)
   \MarkLoc(EdgeBack) \MarkLoc(bk:) \:BkArrowHead
   \CSeg\:EdgeLine(bk:,fr:)
   \CSeg\:EdgeLine(x:,y:) \MarkLoc(bk:)
   \CSeg\:GetExit(2:,y:)  \MoveToNode(#2,\Val\Va,\Val\Vb)
   \MarkLoc(EdgeFront)  \MarkLoc(fr:)  \:FrArrowHead
   \CSeg\:EdgeLine(fr:,bk:)
   \def\EdgeLabel{\:PutLabel\:XXLbl}  } \Define\CycleEdge(1){
   \MoveToNode(#1,0,0)  \MarkLoc(o:)  \MoveF(10)
   \MarkLoc(oo:)  \Rotate(\:CrvEdgeDir)   \MoveF(10)
   \MarkLoc(1':)  \CSeg\:GetExit(oo:,1':)
   \MoveToNode(#1,\Val\Va,\Val\Vb)    \MarkLoc(EdgeBack)
   \MarkLoc(1:) \MoveF(\:CrvEdgeTension)
   \MarkLoc(1':) \DSeg\RotateTo(o:,oo:)
   \MoveToLoc(o:)  \Rotate(-\:CrvEdgeDir)   \MoveF(10)
   \MarkLoc(2':)  \CSeg\:GetExit(o:,2':)
   \MoveToNode(#1,\Val\Va,\Val\Vb) \MarkLoc(EdgeFront)
   \MarkLoc(2:)  \MoveF(\:CrvEdgeTension)
   \MarkLoc(2':) \:InCrvEd   \:EdgeCurve(1:,1':,2':,2:) }

\Define\CycleEdgeSpec(2){
   \:edef\:CrvEdgeDir{#1}  \:edef\:CrvEdgeTension{#2}}

\CycleEdgeSpec(30,20)\EdgeSpec(LTS) \Define\TransEdge(4){
   \IF  \EqText(,#3)  \THEN
      \IF  \EqText(,#4)  \THEN
             \Edge(#1,#2)
      \ELSE  \Edge(#2,#1)
         \EdgeLabel(--#4--) \relax \FI
   \ELSE \IF  \EqText(,#4)  \THEN
      \Edge(#1,#2)    \EdgeLabel(--#3--)\relax
   \ELSE \IF  \EqText(#1,#2)  \THEN
      \RotateTo(#3)   \CycleEdge(#1)
      \EdgeLabel(--#4--)\relax
   \ELSE              \CurvedEdge(#1,#2)
      \EdgeLabel(--#3--) \relax \CurvedEdge(#2,#1)
      \EdgeLabel(--#4--) \relax   \FI  \FI  \FI }
\:DefineExt\:DiagEg(1){\endgroup
   \IF  \EqText(,#1)  \THEN  \ELSE
      \Define\:DiagramEdge{#1}
   \FI }\Define\GridEdge(2){
   \LSeg\Va(#1,#2)
   \IF  \EqDec(\Va,0)  \THEN       \:aldwarn2{#2}    \FI
   \MoveToNode(#1,0,1)  \MarkLoc(1:)
\MoveToNode(#2,0,-1) \MarkLoc(2:)   \CSeg\:GetXY(1:,2:)
\IF  \GtDec(\Vb,-0.1)  \THEN
   \MoveToNode(#1,1,0)  \MarkLoc(1:)
   \MoveToNode(#2,-1,0) \MarkLoc(2:)   \CSeg\:GetXY(1:,2:)
   \IF  \GtDec(\Va,-0.1)  \THEN \EdgeAt(#1,1,1,#2,-1,-1)
   \ELSE
      \MoveToNode(#1,-1,0)  \MarkLoc(1:)
      \MoveToNode(#2,1,0) \MarkLoc(2:)   \CSeg\:GetXY(1:,2:)
      \IF  \LtDec(\Va,0.1)  \THEN  \EdgeAt(#1,-1,1,#2,1,-1)
      \ELSE  \EdgeAt(#1,0,1,#2,0,-1)
    \FI \FI  
   \ELSE \MoveToNode(#1,0,-1)  \MarkLoc(1:)
\MoveToNode(#2,0,1) \MarkLoc(2:)   \CSeg\:GetXY(1:,2:)
\IF  \LtDec(\Vb,0.1)  \THEN
   \MoveToNode(#1,1,0)  \MarkLoc(1:)
   \MoveToNode(#2,-1,0) \MarkLoc(2:)   \CSeg\:GetXY(1:,2:)
   \IF  \GtDec(\Va,-0.1)  \THEN \EdgeAt(#1,1,-1,#2,-1,1)
   \ELSE
      \MoveToNode(#1,-1,0)  \MarkLoc(1:)
      \MoveToNode(#2,1,0) \MarkLoc(2:)   \CSeg\:GetXY(1:,2:)
      \IF  \LtDec(\Va,0.1)  \THEN \EdgeAt(#1,-1,-1,#2,1,1)
      \ELSE   \EdgeAt(#1,0,-1,#2,0,1)
   \FI\FI 
   \ELSE \IF  \GtDec(\Va,0)  \THEN   \EdgeAt(#1,1,0,#2,-1,0)
\ELSE                       \EdgeAt(#1,-1,0,#2,1,0)  \FI
  \FI\FI}    \Define\EdgeToNode(1){{\:EdgeNode(#1)}}

\Define\::EdgeNode(1){ \let&=\relax
   \xdef\:ednd{\:ednd&\Val\:C..\Val\:Ca,#1} }

\Define\::CurrEdge(2){ \:C+1;
   \IF \EqInt(\:C,2) \THEN
      \:CurrEdge(#1,#2)  \:C=0;
   \FI }

   \fi
   \:CheckOption{diagram}\if:option
                \Define\:Temp{
\def\:PutLabel##1{
   \Define\:Lbl{\def\:Temp{\noexpand##1}
      \begingroup   \:DraCatCodes  \:GetLbl}
   \edef\:Temp{\noexpand\:Opt[]
      \noexpand\:Lbl\noexpand{\:LblAdj\noexpand}}
   \:Temp}

\def\:GetLbl[##1]{\endgroup
   \IF  \EqText(+,##1) \THEN
      \edef\:Temp{\:Temp[+\:LblAdj]}
   \ELSE  \edef\:Temp{\:Temp[##1]}  \FI
   \:Temp }
                }\:Temp\Define\:PutLbl(1){
   \MarkLoc(1:)  \XSaveUnits  \Units(1pt,1pt)
   \MoveF(\:LblDis)  \XRecallUnits
   \IF  \EqText(,#1)  \THEN \ELSE    \:AdjLbl(#1)
   \FI    \MarkLoc(2:)    \CSeg\:GetExit(1:,2:)
   \EntryExit(-\Val\Va,-\Val\Vb,0,0) \:EdgeLabel }

\def\:AdjLbl(#1;#2){\Move(#1,#2)}\Define\LabelPos(2){
   \def\:LblAdj{#1}  \def\:LblDis{#2} }
\LabelPos(0.5;,3)          \Define\:Temp{
\def\:LineLbl[##1##2;##3]{
   \MoveToLoc(EdgeBack)   \MarkLoc(a:)
   \MoveToLoc(EdgeFront)  \MarkLoc(b:)
   \DSeg\RotateTo(a:,b:)
   \:C=-1; \:EdgePath(0,99){\:Count}
\IF  \EqInt(\:C,1)  \THEN
    \:EdgePath(1,1){\:LocByAddr\MoveToLoc}
    \MarkLoc(a:)  \MarkLoc(b:)
\ELSE  \IF  \GtInt(\:C,1)  \THEN
   \:C/2;  \:EdgePath(\Val\:C,\Val\:C){\:LocByAddr\MoveToLoc}
   \MarkLoc(a:)
   \:C+1;  \:EdgePath(\Val\:C,\Val\:C){\:LocByAddr\MoveToLoc}
   \MarkLoc(b:)
   \DSeg\RotateTo(a:,b:)
\FI \FI
   \MoveToLoc(a:)
   \IF  \EqText(+,##1) \THEN  \Rotate(-90)
      \CSeg[##2]\Move(a:,b:)
   \ELSE                     \Rotate(90)
      \CSeg[##1##2]\Move(a:,b:)
   \FI
   \:PutLbl(##3) }
          }\:Temp           \Define\:Temp{
\def\:CurveLbl[##1##2;##3]{
   \MoveToLoc(EdgeBack')  \MarkLoc(a:)
   \MoveToLoc(EdgeBack'') \MarkLoc(A:)
   \MoveToLoc(EdgeFront'')   \MarkLoc(B:)
   \MoveToLoc(EdgeFront')    \MarkLoc(b:)
   \DSeg\RotateTo(a:,b:)
   \IF  \EqText(+,##1)  \THEN
       \Rotate(-90)
   \ELSE
       \Rotate(90)
   \FI
   \MoveToCurve[##2](a:,A:,B:,b:) \:PutLbl(##3)  }      }\:Temp
           \Define\:Temp{
\def\:XYLbl[##1##2;##3]{
   \IF  \EqText(+,##1)  \THEN  \Vb=-1; \Va=##2;
                        \ELSE  \Vb= 1; \Va=##1##2;
   \FI
   \MoveToLoc(x:)
   \IF  \LtDec(\Va,0.5)  \THEN
      \MarkLoc(b:)  \MoveToLoc(EdgeBack)  \MarkLoc(a:)
      \DSeg[\Val\Vb]\RotateTo(EdgeFront,x:)
   \ELSE
      \MarkLoc(a:)  \Va-0.5;
      \MoveToLoc(EdgeFront)  \MarkLoc(b:)
      \DSeg[\Val\Vb]\RotateTo(EdgeBack,x:)
   \FI
   \Va*2;
   \MoveToLoc(a:)   \CSeg[\Val\Va]\Move(a:,b:)
   \:PutLbl(##3)  }
            }\:Temp           \Define\:Temp{
\def\:XXLbl[##1##2;##3]{
   \LSeg\Va(x:,y:)
\IF  \LtDec(\Va,0.01) \THEN
   \MoveToLoc(x:)
   \DSeg\RotateTo(EdgeBack,EdgeFront)
   \MoveF(0.1pt\du)  \MarkLoc(y:)
\FI  
   \IF  \EqText(+,##1)  \THEN
       \Vb=-90;  \Va=##2;
   \ELSE
       \Vb=90;   \Va=##1##2;
   \FI
   \IF  \LtDec(\Va,0.333)  \THEN
      \MoveToLoc(x:)    \MarkLoc(b:)
      \MoveToLoc(EdgeBack)  \MarkLoc(a:)
   \ELSE   \IF  \LtDec(\Va,0.666)  \THEN
      \MoveToLoc(x:)   \MarkLoc(a:)
      \MoveToLoc(y:)   \MarkLoc(b:)  \Va-0.3333;
   \ELSE
      \MoveToLoc(y:)   \MarkLoc(a:)
      \MoveToLoc(EdgeFront)   \MarkLoc(b:)  \Va-0.666;
   \FI  \FI
   \Va*3;
   \DSeg\RotateTo(a:,b:)  \Rotate(\Val\Vb)
   \MoveToLoc(a:)   \CSeg[\Val\Va]\Move(a:,b:)
   \:PutLbl(##3)  }
            }\:Temp    \Define\LabelSpec{
   \TextPar\Define\:EdgeLabel }

\LabelSpec(1){\Text(--#1--)}
\fi
\:CheckOption{spread}\if:option
   \:DefineExt\DiagramSpec(1){\endgroup
    \:C=0;  \let\:Next=\:DiagNodes  \:Next#1&&
    \Indirect\let<DiagNod:\Val\:C>=\relax }

             \Define\:Temp{
\def\:DiagNodes##1&{
   \def\:Temp{##1}
   \ifx \:Temp\empty  \let\:Next=\relax
   \else
      \:tk={##1} \:Ca=\:C;  \:C+1;
      \edef\:Temp{\noexpand\:AddDiagNodes
         (\Val\:C,\Val\:Ca,\the\:tk)}   \:Temp
   \fi  \:Next }
              }\:Temp

\Define\:AddDiagNodes(3){
   \Indirect\def<DiagNod:#2>{\begingroup\:Spaces
      \:DiagNds(#1,#3)}}

       \Define\:Temp{
\def\:DiagNds(##1,##2)(##3){\endgroup
   \Table\:Temp{##3} \:SetDiaNodes##2
   \Indirect<DiagNod:##1>}
        }\:Temp \Define\Diagram{ \Indirect<DiagNod:0> }

\def\:SetDiaNodes#1{\:Temp(0,999){#1}}
\DiagramSpec(\DiagramNode&\Edge)
\Define\DiagramNode(4){ \MoveTo(#2,#3) \Node(#1)(--#4--) }  
\fi
\:CheckOption{grid}\if:option
   \:DefineExt\GridDiagramSpec(1){\endgroup
   \IF  \EqText(,#1)  \THEN  \ELSE
      \:CheckSpecTyp(#1,&)
      \IF  \:Temp  \THEN  \Table\:PolyShapes{u,\noexpand#1}
\Define\:ShapeTyp{u,}
      \ELSE               \Define\:ShapeTyp{}
\:SetPolyNodes#1&,&      \let\:Temp=\:PolyShapes
\Table\:PolyShapes{\:Temp}      \FI
   \FI
   \:DiagEg }\Define\:CheckSpecTyp(2){
   \IF  \EqText(#2,&)  \THEN  \def\:Temp{\EqInt(1,1)}
   \ELSE                      \def\:Temp{\EqInt(1,2)}
   \FI  }\GridDiagramSpec(\Node)(\GridEdge)
                                         \Define\:Temp{
\def\:SetPolyNodes##1,##2&{
   \Define\:PolyShapes{##1,\noexpand##2}
   \:MorePolyNodes   }                   }\:Temp

\Define\:MorePolyNodes{
    \def\:Temp{\Define\:PolyShapes}
    \Define\:Next{\:GetPolyNodes}    \:Next  }

                                               \Define\:Temp{
\def\:GetPolyNodes##1,##2&{
   \IF \EqText(,##1) \THEN     \Define\:Next{}
   \ELSE
      \expandafter\expandafter\expandafter\:Temp\expandafter{
              \:PolyShapes & ##1,\noexpand##2 }
   \FI    \:Next}}                              \:Temp
\Define\GridSpace(2){
   \:edef\:VGridSp{#2}   \:edef\:HGridSp{#1}}
\GridSpace(20,20)\Define\AlignGrid(2){
   \def\:XGridEntry{#1}     \def\:YGridEntry{#2}}
\AlignGrid(0,0)

\Define\GridDiagram(2){  \MarkLoc(GridOrg)
   \IF  \LtInt(#1,2)  \:aldwarn{11}{#1,#2} \THEN  \FI
   \IF  \LtInt(#2,2)  \:aldwarn{11}{#1,#2} \THEN  \FI
   \EntryExit(\:XGridEntry,\:YGridEntry,0,0)
   \Do(0,#1){ \:C=\DoReg;  \Indirect\DecVar<y.\Val\:C>  }
   \Do(0,#2){ \:C=\DoReg;  \Indirect\DecVar<x.\Val\:C> }
   \Define\:Nrows{#1}       \Define\:Ncols{#2}
   \:VGridSpace }
                   \Define\:VGridSpace(1){
   \IF  \EqText(,#1)  \THEN
          \Define\:Temp{\:VGridSp,0,1000}
   \ELSE  \Define\:Temp{#1}              \FI
   \Table\:VerAdj{\:Temp}
   \:C=0;  \:VerAdj(0,99){\:Count}
   \IF  \LtInt(\:C,\:Nrows) \THEN
      \:C=-\:C;   \:C+\:Nrows;    \:C-1;
\IF \LtInt(\Val\:C,\:Nrows)  \THEN
   \def\:tempA{\def\:Temp}
   \IF  \GtInt(\Val\:C,0)  \THEN
      \Do(1,\Val\:C){  \expandafter\expandafter
         \expandafter\:tempA\expandafter{
                   \:Temp & \:VGridSp,0,1000}}
   \FI
   \Table\:VerAdj{\:Temp}
\FI
   \ELSE   \IF \GtInt(\:Nrows,1) \THEN
      { \def&{\string&}\:aldwarn4{#1} }
   \FI   \FI
   \:HGridSpace}
                   \Define\:HGridSpace(1){
   \IF  \EqText(,#1)  \THEN
          \Define\:Temp{\:HGridSp,0,1000}
   \ELSE  \Define\:Temp{#1}              \FI
   \Table\:HorAdj{\:Temp}
   \:C=0; \:HorAdj(0,99){\:Count}
   \IF  \LtInt(\:C,\:Ncols) \THEN
      \:C=-\:C;   \:C+\:Ncols;   \:C-1;
\IF \LtInt(\:C,\:Ncols)  \THEN
   \def\:tempA{\def\:Temp}
   \IF  \GtInt(\:C,0)  \THEN
      \Do(1,\Val\:C){  \expandafter\expandafter
         \expandafter\:tempA\expandafter{
                   \:Temp & \:HGridSp,0,1000}}
   \FI
   \Table\:HorAdj{\:Temp}
\FI
   \ELSE \IF \GtInt(\:Ncols,1) \THEN
       { \def&{\string&}\:aldwarn4{#1} }
   \FI   \FI
   \:GridEntries}

                   \Define\:GridEntries{\begingroup \:Spaces  \::GridEntries}

\Define\::GridEntries(1){\endgroup
   \:bgDiTags
   \:C=0;  \Define\:Next{\:C+1;\:GetGridRows}  \:Next#1//
   \Do(1,\:Nrows){\:C=\DoReg;  \:Ca=0;
   \Indirect<row.\Val\:C>(0,99){\:Ca+1;
      \IF  \GtInt(\Val\:C,\:Nrows)  \THEN  \:aldwarn5{}  \FI
      \IF  \GtInt(\Val\:Ca,\:Ncols) \THEN  \:aldwarn6{}  \FI
      \:GetEntrySize}}
   \IF  \GtInt(\:Nrows,1)  \THEN
   \Do(2,\:Nrows){                 \:C=\DoReg;
   \Indirect{\Va=-}<y.\Val\:C>;  \:C-1;
   \Indirect{\Va+}<y.\Val\:C>;
   \Va*\:YGridEntry;            \:C+1;
   \Indirect{\Va+}<y.\Val\:C>;  \:C-1;
   \Indirect<y.\Val\:C>+\Va;         }
\:C=0;  \:VerAdj(0,99){\:C+1; \:AdjGridHeight}
\Do(2,\:Nrows){     \:C=\DoReg;  \:C-1;
   \Indirect{\Va=}<y.\Val\:C>;  \:C+1;
   \Indirect<y.\Val\:C>+\Va;         }    \FI
\IF  \GtInt(\:Ncols,1)  \THEN
   \Do(2,\:Ncols){                   \:C=\DoReg;
   \Indirect{\Va=+}<x.\Val\:C>;  \:C-1;
   \Indirect{\Va-}<x.\Val\:C>;
   \Va*\:XGridEntry;             \:C+1;
   \Indirect{\Va+}<x.\Val\:C>;  \:C-1;
   \Indirect<x.\Val\:C>+\Va;         }
\:C=0;  \:HorAdj(0,99){\:C+1; \:AdjGridWidth}
\Do(2,\:Ncols){     \:C=\DoReg;  \:C-1;
   \Indirect{\Va=}<x.\Val\:C>; \:C+1;
   \Indirect<x.\Val\:C>+\Va;         }    \FI
\Do(1,\:Nrows){  \:C=\DoReg;
   \Indirect{\Indirect<y.\Val\:C>=-}<y.\Val\:C>; }
   \let\:TagNode=\::TagNode
\let\:EdgeNode=\::EdgeNode

\Do(1,\:Nrows){  \:C=\DoReg;
   \Do(1,\:Ncols){   \:Ca=\DoReg;
      \Indirect<\Val\:C..\Val\:Ca>}}
    \Do(1,\:Nrows){  \:C=\DoReg;
   \Do(1,\:Ncols){   \:Ca=\DoReg;
      \Indirect\let<\Val\:C..\Val\:Ca>=\:undefined}}
   \Table\:PolyBaseEdges{\:BaseEdges}
\:C=-1;  \:PolyBaseEdges(0,1){\:Count}
\IF  \EqInt(\:C,1)  \THEN
   \:PolyBaseEdges(1,999){\:InsertBaseEdge}
\FI
   \let\:CurrEdge=\:DiagramEdge
\:endDiTags }
                              \Define\:Temp{
\def\:GetGridRows##1//{
   \IF  \EqText(,##1)  \THEN  \let\:Next=\relax
   \ELSE   {\:C-1; \:TrcDiag{\immediate\write16{...read\space
                                        row\space \Val\:C}}}
      \def\Defend{\noexpand\noexpand\noexpand}
\def\:Comma{,}  \let\:CsCom=\,
\def\,{\noexpand\:Comma}
\let\:Temp=&  \def&{\:Temp\:ShapeTyp}

      \Indirect\Table<row.\Val\:C>{\:ShapeTyp##1}
      \let\Defend=\noexpand\let&=\:Temp    
   \FI  \:Next }              }\:Temp\Define\:GetEntrySize(1){
   \IF \EqText(,#1)  \THEN \ELSE
     \:TrcDiag{{  \:C-1;  \:Ca-1;  \:tk={#1}
    \immediate\write16{\Val\:C..\Val\:Ca\space\space\the\:tk}}}
     \:GetEntrySizeA(#1,)  \FI  }\Define\:GetEntrySizeA(3){
   \def\:Temp{\:tk=}    \:AddBaseEdges(#3)
   \:PolyShapes(0,99){\:GetEntrySizeB(#1,#2)}}               \Define\:Temp{
\def\:GetEntrySizeB(##1,##2)(##3,##4){
   \IF  \EqText(,##2)  \THEN \ELSE
      \IF  \EqText(##1,##3)  \THEN
         \Indirect\Define<\Val\:C..\Val\:Ca>{
   \:AddGridNode(##4,##2)}
         \let\:SvWg=~ \let~=\relax
         \edef\:Temp{\noexpand##4(.;)(--##2--)}
         \let~=\:SvWg
         \FigSize\Va\Vb{\:Temp}   \Va/2; \Vb/2;
         \Indirect{\:Ve=}<x.\Val\:Ca>;
         \IF  \GtDec(\Va,\:Ve)  \THEN
            \Indirect<x.\Val\:Ca>=\Va;      \FI
         \Indirect{\:Ve=}<y.\Val\:C>;
         \IF  \GtDec(\Vb,\:Ve)  \THEN
            \Indirect<y.\Val\:C>=\Vb;
   \FI \FI \FI }}  \:Temp  \Define\:AddGridNode(2){
   \:C-1; \:Ca-1;    \MoveToGrid(\Val\:C,\Val\:Ca)
   \let\:SvWg=~\let~=\relax
   \edef\:Temp{ \noexpand #1(\Val\:C..\Val\:Ca)(--#2--) }
   \let~=\:SvWg
   \:Temp    \:C+1; }\Define\MoveToGrid(2){
    \IF \LtInt(#1,0) \THEN  \:aldwarn3{(#1,#2)}  \FI
    \IF \LtInt(#2,0) \THEN  \:aldwarn3{(#1,#2)}  \FI
    \IF \GtInt(#1,\:Nrows) \THEN  \:aldwarn3{(#1,#2)}  \FI
    \IF \GtInt(#2,\:Ncols) \THEN  \:aldwarn3{(#1,#2)}  \FI
    \Indirect{\Va=}<x.#2>;
    \Indirect{\Vb=}<y.#1>;
    \MoveToLoc(GridOrg)
    \Move(\Val\Va pt\du,\Val\Vb pt\du)             }
\Define\:AddBaseEdges(1){
   \IF  \EqText(,#1)  \THEN  \def\:Next{}
   \ELSE
      \:AddBaseEdgesB(#1)
   \FI \:Next}

\Define\:AddBaseEdgesB(2){
   \expandafter\expandafter\expandafter\:Temp
      \expandafter{\:BaseEdges &}
   \edef\:BaseEdges{  \the\:tk \Val\:C,\Val\:Ca,#1 }
   \def\:Next{\:AddBaseEdges(#2)}}

\Define\:BaseEdges{0,0,0,0}\Define\:AdjGridWidth(3){
   \Indirect{\Va=}<x.\Val\:C>;    \Va+#1;
   \IF \LtDec(\Va,#2) \THEN \Va=#2;  \ELSE
   \IF \GtDec(\Va,#3) \THEN \Va=#3;  \FI  \FI
   \Indirect<x.\Val\:C>=\Va;  }\Define\:AdjGridHeight(3){
   \Indirect{\Va=}<y.\Val\:C>;    \Va+#1;
   \IF \LtDec(\Va,#2) \THEN \Va=#2;  \ELSE
   \IF \GtDec(\Va,#3) \THEN \Va=#3;  \FI  \FI
   \Indirect<y.\Val\:C>=\Va;  }\Define\:InsertBaseEdge(3){
   \:C=#1;  \:Ca=#2;  \:TargetNode(#3,)
\IF  \EqInt(\:C,1) \THEN  \:GetGridLbl(#3) \FI
\:C=#1;    \:Ca=#2;
\expandafter\:GetTarget\:Temp..
   \:C=#1;   \:C-1;   \:Ca=#2;  \:Ca-1;
   \edef\:Temp{\noexpand\:DiagramEdge(\Val\:C..\Val\:Ca,\:Temp\the\:tk)}
   \:Temp} \Define\:TargetNode(2){
   \def\:Temp{#1}
   \IF \EqText(,#2)  \THEN \:C=0;  \:tk={}
   \ELSE \:C=1;  \FI }

\Define\:GetGridLbl(2){ \:tk={,#2} }        \Define\:Temp{
\def\:GetTarget##1##2..##3##4..{
   \:C-1; \:Ca-1;  \edef\:tempa{\Val\:C..\Val\:Ca}
   \IF      \EqText(+,##1)  \THEN
       \:C+##2;  \:Ca##3##4;  \edef\:Temp{\Val\:C..\Val\:Ca}
   \ELSE\IF \EqText(-,##1)  \THEN
       \:C-##2;  \:Ca##3##4;  \edef\:Temp{\Val\:C..\Val\:Ca}
   \ELSE                     \def\:Temp{##1##2..##3##4}  \FI\FI
   \:TrcDiag{\immediate\write16{\:tempa-->\:Temp}}    }
         }\:Temp
\fi
\:CheckOption{tree}\if:option
   \Define\:MoveH(2){
   \ifx H\:TreeDir   \Move(0,-#2)
   \else             \Move(#1,0)  \fi  }

\Define\:MoveV(2){
   \ifx H\:TreeDir   \Move(#1,0)
   \else             \Move(0,#2)  \fi  }

\Define\:MoveHV(2){
   \ifx H\:TreeDir   \Move(0,#2)
   \else             \Move(#1,0)  \fi  }

\Define\:DistH(2){
   \ifx H\:TreeDir   \Va=#2;  \else  \Va=#1;  \fi}
\def\:oochild{0} \Define\:MoveLFr(2){
  \ifx H\:TreeDir \MoveToNode(#2,0,-1) \else \MoveToNode(#2,1,0)  \fi}
\Define\:FrDs(2){  \MoveToLoc(:a) \MarkLoc(:a')
   \if H\:TreeDir\relax  \MoveToNode(#2,0,-1)
   \else \MoveToNode(#2,1,0) \fi
   \MarkLoc(:a)
   \:tempb(\DoReg,\DoReg){\:MoveRFr}
   \CSeg\:AddX(:a,:a')  \CSeg\:AddX(:b',:b)  
   \ifx H\:TreeDir
   \IF \GtDec(\:va,0) \THEN   \:vb-\:va; \:va=0;   \FI
\else
   \IF \LtDec(\:vb,0) \THEN   \:va-\:vb; \:vb=0;   \FI
\fi
          }

\Define\:MoveRFr(2){
   \MoveToLoc(:b)       \MarkLoc(:b')
   \ifx H\:TreeDir \MoveToNode(#1,0,1) \else \MoveToNode(#1,-1,0)  \fi
   \MarkLoc(:b)  }
\Define\:AddX(2){ \ifx H\:TreeDir  \:va+#2;  \else  \:vb+#1; \fi }
\Define\:Temp{  \let\:aprsnd=& }
\:Temp \Define\:left(2){ \def\:frL{#1}}
\Define\:right(2){         \let&=\relax
   \edef\:frL{\:frL,#2}    \let&=\:aprsnd  }
\Define\:copy(1){          \let&=\relax
   \edef\:frL{\:frL & #1}  \let&=\:aprsnd  }

\Define\:addleft(2){   \let&=\relax
   \edef\:frL{\:frL&#1}   \let&=\:aprsnd}

\Define\:both(1){
   \IF \GtInt(#1,0) \THEN
      \Do(1,#1){
         \:tempa(\DoReg,\DoReg){\:addleft}
         \:tempb(\DoReg,\DoReg){\:right}  }
    \FI  }           
   \:DefineExt\TreeSpec(1){\endgroup
   \IF  \EqText(,#1)  \THEN   \ELSE
      \:TrNodes#1&&\let  \:LeafNode=\:TreeNode
   \FI  \:DecideLeafNode  }

\def\:Temp{\catcode`\ =10 \catcode`\^^M=5 \catcode`\&=4 }
\catcode`\ =9 \catcode`\^^M=9  \catcode`\&=13
\def\:TrNodes#1&#2&{
   \IF  \EqText(,#2)  \THEN
      \Define\:TreeNode(2){#1(##1)(--##2--)}
      \let\:Next=\relax
   \ELSE
      \Define\:TreeNode{
         \expandafter\:CondTreeNode\:pars(#1)
         \expandafter\:CondTreeNode\:pars(#2)  }
       \def\:tempa{ \Define\:TreeNode }
       \let\:Next=\:DecideTreeNodes
   \FI \:Next}

\def\:DecideTreeNodes#1&{
   \IF  \EqText(,#1)  \THEN \let\:Next=\relax
      \def\:tempa{\Define\:TrNd}
      \expandafter\expandafter\expandafter
         \:tempa\expandafter{  \:TreeNode  \:Temp  }
      \def\:TreeNode(##1,##2){
           \def\:pars{(##1)(##2)}  \:TrNd}
   \ELSE
      \expandafter\expandafter\expandafter
         \:tempa\expandafter{  \:TreeNode
             \expandafter\:CondTreeNode\:pars(#1)}
   \FI  \:Next}
                \:Temp

                \Define\:Temp{
\def\:CondTreeNode(##1)(##2,##3)(##4,##5){
   \IF \EqText(##2,##4)  \THEN
      \def\:Temp{##5(##1)(--##3--)}
   \FI}
                }\:Temp\:DefineExt\:DecideLeafNode(1){\endgroup
   \IF  \EqText(,#1)  \THEN  \ELSE
      \Define\:LeafNode(2){#1(##1)(--##2--)}
   \FI
   \:DecideTreeEdge}\:DefineExt\:DecideTreeEdge(1){\endgroup
   \IF \EqText(,#1) \THEN  \ELSE
      \Define\:TreeEdge(2){
   \Indirect{\let\:Temp=} <##2.mvto>
   \ifx \:Temp\:ChkInv \else
      \:TrcDiag{\immediate\write16{##1,##2}}
      \ifx H\:TreeDir   \MoveToNode(##1,\:EntryX,0)
      \else             \MoveToNode(##1,0,\:EntryY)  \fi
      \:MvToDivider##1..   \MarkLoc(EdgeGuide)
      #1(##1,##2)
   \fi } 
   \FI} \Define\Tree(1){
   \Define\:PolyVerTreeSp{#1}
   \Define\:BareTr(1){,##1}
   \begingroup   \:Spaces  \:Tree }

\Define\LTree(1){
   \Define\:PolyVerTreeSp{#1} \:LTree  }

\catcode`\#=13
\Define\:LTree{\def#{\Defend&}  \let\-=#  \let\:BareTr=\relax
  \begingroup \catcode`\#=13   \:Spaces  \:Tree }
\catcode`\#=6

\Define\:Tree(1){\endgroup
   \:bgDiTags  
   \let\:MergTree=\relax
\ifx \:Cspace\:InInx  \let\:MergTree=\::MergTree  \fi
\:TrSp(0,99){\:copt}
   \ifx H\:TreeDir      \let\:tempc=\Va
   \def\:TreeEnEx{\EntryExit(\:EntryX,\:EntryY,0,0)}
\else                \let\:tempc=\Vb
   \def\:TreeEnEx{\EntryExit(\:EntryX,\:EntryY,0,0)}  \fi
\:C=-1;
   \let\:Next=\:GetTreeRows  \:Next#1//
   \Do(1,\Val\:C){ \:Ca=\DoReg;
   \Indirect{\Va=-}<y.\Val\:Ca>;  \:Ca-1;
   \Indirect\DecVar<x.\Val\:Ca>
\ifx H\:TreeDir \Indirect<x.\Val\:Ca>=-\:EntryX;
\else           \Indirect<x.\Val\:Ca>=\:EntryY;  \fi
\Indirect<x.\Val\:Ca>+1;
\Indirect{\Indirect<x.\Val\:Ca>*}<y.\Val\:Ca>;
\Indirect<x.\Val\:Ca>/2;  
   \Indirect{\Va+}<y.\Val\:Ca>;
   \ifx H\:TreeDir   \Va*-\:EntryX;
   \else             \Va*\:EntryY;  \fi
   \:Ca+1;  \Indirect{\Va+}<y.\Val\:Ca>;  \:Ca-1;
   \Indirect<y.\Val\:Ca>+\Va;
   \Indirect<y.\Val\:Ca>/2;}\:Ca=-1;
\ifx  \:PolyVerTreeSp\empty  \else
   \let\:tempa=\:PolyVerTreeSp
   \Table\:PolyVerTreeSp{\:tempa}
  \:PolyVerTreeSp(0,99){\:AddVerTreeSp}
\fi
\:Cb=\:Ca;  \:Cb+1;
\IF \GtInt(\:C,\:Cb)  \THEN
  \:Cb+1;
   \Do(\Val\:Cb,\Val\:C){\:AddVerTreeSp(\:VerTreeSp)}
\ELSE \IF \GtInt(\:Cb,\:C) \THEN
   \:aldwarn1{}\FI\FI

   \Indirect{\let\:va=}<x.0>
\Indirect{\let\:vb=}<y.0> 
   \let\:TagNode=\::TagNode
\let\:EdgeNode=\::EdgeNode

\Do(0,\Val\:C){ \:Ca=-1;  \:Cb=-1;
  \Indirect<row.\DoReg>(0,99){
      \:Ca+1;  \:DefSubTree}  }
\gdef\:GetNodesAddr{}    \:AddrCount=1;
\Indirect<0..0>  \global\let\:GetNodesAddr=\:undefined
   \ifx \:BareTr\relax     \let\:EnterTreeEdge=\:::EnterTreeEdge   
\else                \let\:EnterTreeEdge=\::EnterTreeEdge    \fi
\:C-1;
\Do(0,\Val\:C){ \:C=\DoReg; \:Ca=-1;  \:Cb=\DoReg;    \:Cb+1;  \:Cc=-1;
  \Indirect<row.\DoReg>(0,99){
      \:Ca+1;  \:EnterTreeEdge}  } 
   \let\:CurrEdge=\:TreeEdge
\let\:tempa=\:MvToDivider
\def\:MvToDivider##1..##2..{\Move(0.001pt\du,0.001pt\du)}
\:endDiTags  \let\:MvToDivider=\:tempa  
   \let\-=\:svneg  }

\let\:svneg=\-                     \Define\:Temp{
\def\:GetTreeRows##1//{
   \def\:Temp{##1}
   \ifx \:Temp\empty   \let\:Next=\relax  \else
      \:C+1;
      \ifx  \:BareTr\relax  \else
         \TableData\:Temp\:BareTr{##1}
      \fi
      \Indirect\Table<row.\Val\:C>{\:Temp}
      \def\:tempa{}   \:Vc=0;
\def\:tempb{\:tk=}     \:Ca=-1;
\Indirect<row.\Val\:C>(0,99){\:Ca+1; \:MeasureRow}
\Indirect\DecVar<y.\Val\:C>
\Indirect<y.\Val\:C>=\:Vc;
   \fi  \:Next}
                    }\:Temp\Define\::EnterTreeEdge(3){
   \IF \GtInt(#2,0) \THEN
      \Do(1,#2){    \:Cc+1;
         \edef\:temp{\noexpand\:TreeEdge
            (\Val\:C..\Val\:Ca,\Val\:Cb..\Val\:Cc)
            \noexpand\:C=\Val\:C;  \noexpand\:Ca=\Val\:Ca;}
                       \:temp  }
   \FI }\Define\:::EnterTreeEdge(3){
   \edef\:temp{\:tk={\noexpand\:C=\Val\:C;
                     \noexpand\:Ca=\Val\:Ca;}} \:temp
   \IF \GtInt(#2,0) \THEN
      \Table\:BareTr{#1}
      \:BareTr(0,99){    \:Cc+1;  \the\:tk
         \edef\:temp{\noexpand\:TreeEdgePos
            (\Val\:C..\Val\:Ca,\Val\:Cb..\Val\:Cc)  }
            \:temp  \:TreeLabel }
   \FI     \the\:tk  }

\Define\:TreeEdgePos(2){  \:TreeEdge(#1,#2)  \LTreePos(#1,#2)  }

\Define\:TreeLabel(1){ \EdgeLabel(--#1--) }
\Define\LTreePos(2){   \DSeg\Va(#1,#2)
  \IF  \LtDec(\Va,270) \THEN \LabelPos(0.65;-3;3,)
  \ELSE \LabelPos(+0.65;3;3,)     \FI  }\Define\:DefSubTree(3){ \:tk={#3}
   \edef\:Temp{
      \noexpand\::DefSubTree(#2,\the\DoReg,\Val\:Ca,\Val\:Cb,\the\:tk)}
   \:Cb+#2;   \:Temp}\Define\::DefSubTree(5){
   \Indirect\Object<#2..#3>{
      \let\:ochild=\:oochild \let\:currFrL=\:frL  \let\:currDpL=\:dpL
\def\:dpL{-1} \def\:frL{} 
      \:Ca=#3; \:Cc=#2;
      \:TrSp(0,99){\:unl}
      \:Cc+1;  \:Cb=#4;  \:Cb+1;
      \IF  \GtInt(#1,0) \THEN
         \edef\:parmA{#1} \edef\:parmB{#2}
         \:tk={#5}  \def\:parmC{\the\:tk}
         \def\:Temp{\:Cind }
      \ELSE   \:tk={#5}
         \edef\:Temp{\noexpand  \:Clnd  (#2,\noexpand\the\:tk) }
      \FI
      \edef\:Temp{\expandafter\noexpand\:Temp
         \noexpand\:endSubTr (#1,#2,#3)}
      \:Temp }}\Define\:endSubTr(3){
   \Indirect{\Indirect
   { \global\let}<#2..\Val\:Ca.mvto>= }<#2..\Val\:Ca.mvto>
\MoveToLoc(#2..\Val\:Ca)     \MarkXLoc(#2..\Val\:Ca)
\MoveToLoc(#2..\Val\:Ca;:11)  \MarkXLoc(#2..\Val\:Ca;:11)
\edef\:temp{\noexpand\:AddNodesAddr(#2..\Val\:Ca) }  \:temp
   \:Ca=\:dpL;  \:Ca+1; \xdef\:dpR{\Val\:Ca} \let&=\relax
\xdef\:frR{ #2..#3,#2..#3 & \:frL }
\ifx  \:Cspace\:InInx
   \let&=\:aprsnd
   \:Ca=#2;  \:Ca*#1;  \:Ca*\:ochild;  \:Ca*\:lchild;
\IF  \GtInt(\:Ca,0)  \THEN
    \IF \GtInt(\:currDpL,\:dpR) \THEN \def\:tempc{\:dpR}
\ELSE                         \def\:tempc{\:currDpL}  \FI
\Table\:tempa{\:currFrL}
\Table\:tempb{\:frR} \if H\:TreeDir  \MoveToNode(#2..#3,0,1)
\else \MoveToNode(#2..#3,-1,0) \fi            \MarkLoc(:b)
\:tempa(0,0){\:MoveLFr}   \MarkLoc(:a)
\:va=0;  \:vb=0;  \IF  \GtInt(\:tempc,0)  \THEN
    \let\:MoveToInv=\::MoveToInv  
    \DoReg=0\relax   \:tempa(1,\:tempc){\advance\DoReg by 1 \:FrDs}
\FI  
    \MoveToLoc(a:)    \:MoveH(-\Val\:va,\Val\:vb)
    \MarkXLoc(INc:)
\FI   \fi  } \Define\:Cind{
   \edef\:lchild{\Val\:Cb}  
\def\:oochild{0}  
\:trcTree
\Indirect<\Val\:Cc..\Val\:Cb>[\:InInx]   \:SubTreeAddrs
\let\:dpL=\:dpR  \let\:frL=\:frR  
\MoveToLoc(INx:)  \MarkLoc(-12:)
\MoveToLoc(IN:)   \MarkLoc(-11:)
\MoveToLoc(OUTx:) \MarkLoc(12:)
\MoveToLoc(OUT:)  \MarkLoc(11:) 
\IF  \GtInt(\:parmA,1)  \THEN
   \def\:oochild{1}  
   \Do(2,\:parmA){
      \MoveToLoc(\:OutOutx)
      \:MoveH(\:HorTreeSp,\:HorTreeSp)   
      \:TrSp(0,99){\:spbn}    \:Cb+1;
      \:trcTree
      \Indirect<\Val\:Cc..\Val\:Cb>[\:InInx]  \:SubTreeAddrs
      \MarkLoc(o:)
\MoveToLoc(INx:)  \MarkLoc(-\DoReg2:)
\MoveToLoc(IN:)   \MarkLoc(-\DoReg1:)
\MoveToLoc(OUTx:) \MarkLoc(\DoReg2:)
\MoveToLoc(OUT:)  \MarkLoc(\DoReg1:)
\MoveToLoc(o:) 
      \:MergTree    }
\FI \:C=\:parmB;    \:TrSp(0,99){\:lcalg}
\:ChildLoc(\:AlignA,-12:,-11:,11:,12:)   \MarkLoc(a:)
\:ChildLoc(\:AlignB,INx:,IN:,OUT:,OUTx:) \MarkLoc(b:)
\:C=\:parmB;
\IF  \LtDec(\:AlignC,-1)  \THEN
   \MoveToLoc(a:)  \:MoveH(\:AlignC,\:AlignC)
\ELSE  \IF  \GtDec(\:AlignC,1)  \THEN
   \MoveToLoc(b:)  \:MoveH(\:AlignC,\:AlignC)
\ELSE
   \CSeg[0.5]\Move(b:,a:) \MarkLoc(o:)
   \CSeg[\:AlignC]\Move(o:,b:)
\FI  \FI
\:TrSp(0,99){\:mvpr}
\ifx  H\:TreeDir   \Indirect{\Move(-\Val}<y.\Val\:C>,0)
\else              \Indirect{\Move(0,\Val}<y.\Val\:C>)   \fi
 \MarkLoc(o:)  \:TreeEnEx
\edef\:Temp{\noexpand\:TreeNode(\:parmB..\Val\:Ca,\:parmC)}\:Temp
\ifx H\:TreeDir
   \MoveToNode(\:parmB..\Val\:Ca,0,1)   \MarkLoc(a:)
   \MoveToNode(\:parmB..\Val\:Ca,0,-1)  \MarkLoc(b:)
\else
   \MoveToNode(\:parmB..\Val\:Ca,-1,0) \MarkLoc(a:)
   \MoveToNode(\:parmB..\Val\:Ca,1,0)  \MarkLoc(b:)
\fi
\MoveToLoc(o:)   \CSeg\:MoveHV(o:,a:)
\MarkLoc(a:)  \MarkXLoc(IN:)
\MarkXLoc(INc:) 
\MoveToLoc(o:)   \CSeg\:MoveHV(o:,b:)
\MarkXLoc(OUT:)  \MarkLoc(b:)
\MoveToLoc(o:)   \CSeg\:MoveHV(o:,-12:)   \MarkLoc(-12:)
\MoveToLoc(a:)   \CSeg\:DistH(a:,-12:)
\IF  \LtDec(\Va,0) \THEN
   \ifx  H\:TreeDir
   \else             \Move(\Val\Va,0)   \fi
\ELSE
   \ifx  H\:TreeDir  \Move(0,\Val\Va)
   \else            \fi
  \FI
\MarkXLoc(INx:)
\MoveToLoc(o:)   \CSeg\:MoveHV(o:,OUTx:)   \MarkLoc(OUTx:)
\MoveToLoc(b:)   \CSeg\:DistH(b:,OUTx:)
\IF  \GtDec(\Va,0) \THEN
   \ifx  H\:TreeDir
   \else             \Move(\Val\Va,0)   \fi
\ELSE
   \ifx  H\:TreeDir  \Move(0,\Val\Va)
   \else               \fi
 \FI
\MarkXLoc(OUTx:)  }
\Define\:Clnd(2){
   \MarkLoc(o:)  \:TreeEnEx  \:C=#1;  \:LeafNode(#1..\Val\:Ca,#2)
\ifx H\:TreeDir
   \MoveToNode(#1..\Val\:Ca,0,1) \MarkLoc(IN:)
   \MoveToNode(#1..\Val\:Ca,0,-1)
\else
   \MoveToNode(#1..\Val\:Ca,-1,0) \MarkLoc(IN:)
   \MoveToNode(#1..\Val\:Ca,1,0)               \fi
\MarkLoc(OUT:)   \MoveToLoc(o:)   \CSeg\:MoveHV(o:,IN:)
\MarkXLoc(IN:)   \MarkXLoc(INc:)    \MarkXLoc(INx:)
\MoveToLoc(o:)   \CSeg\:MoveHV(o:,OUT:)
\MarkXLoc(OUT:)  \MarkXLoc(OUTx:)   }
\Define\TreeAlign(3){ \let\:TreeDir=#1
   \def\:EntryX{#2}  \def\:EntryY{#3}    \:AlignTree}

\Define\:AlignTree(3){
   \def\:AlignA{#1}   \def\:AlignB{#2}   \def\:AlignC{#3}}

\TreeAlign(V,0,0)(0,0,0)
\Define\:trcTree{\:TrcDiag{\immediate\write16{\space
   \:parmB..\Val\:Ca,\Val\:Cc..\Val\:Cb}}}
\Define\:ChildLoc(2){
\IF  \LtDec(#1,0)  \THEN \Va=-#1;
                   \ELSE \Va= #1;  \FI
\IF \LtDec(\Va,100) \THEN
   \::ChildLoc(#1,#2)
\ELSE   \Va/100; \:C[\Va]; \Vb=#1;
   \IF  \GtDec(#1,0)   \THEN  \Vb-\:C 00;
   \ELSE  \Vb+\:C 00; \FI
   \edef\:Temp{\noexpand\::ChildLoc
      (\Val\Vb,-\Val\:C2:,-\Val\:C1:,\Val\:C1:,\Val\:C2:)}\:Temp
\FI } \Define\::ChildLoc(5){
   \IF  \LtDec(#1,-2)  \THEN
      \MoveToLoc(#2)  \:MoveH(#1,#1)
   \ELSE \IF  \GtDec(#1,2)  \THEN
      \MoveToLoc(#4)  \:MoveH(#1,#1)
   \ELSE \IF  \LtDec(#1,-1)  \THEN
      \MoveToLoc(#3)   \Va=#1;  \Va+1;
      \CSeg[\Val\Va]\Move(#2,#3)
   \ELSE \IF  \GtDec(#1,1)  \THEN
      \MoveToLoc(#4)   \Va=#1;  \Va-1;
      \CSeg[\Val\Va]\Move(#4,#5)
   \ELSE   \MoveToLoc(#3)
      \CSeg[0.5]\Move(#3,#4)  \:MoveH(#1,#1)
   \FI \FI \FI \FI}\def\:GetIx#1..#2//{\:Ca=#1;\:Cb=#2;}

\Define\::MoveToInv(3){
   \Indirect{\let\:temp=}<#1.mvto>
   \expandafter\ifx\:temp\MoveToInv
      \:GetIx#1//        \ifx H\:TreeDir \:Cb-#3; \else \:Cb-#2; \fi
      \edef\:Temp{\noexpand\MoveToInv(\Val\:Ca..\Val\:Cb,
                     \ifx H\:TreeDir 0,#3 \else -#2,0 \fi)}
      \expandafter\:Temp
   \ELSE
      \ifx H\:TreeDir  \Indirect<#1.mvto>(#1,0,#3)
      \else            \Indirect<#1.mvto>(#1,-#2,0)  \fi
   \FI    }

\def\MoveToInv{\:MoveToInv}
\let\:MoveToInv=\MoveToRect \Define\::MergTree{
   \edef\:Temp{\DoReg=\the\DoReg\relax}
   \Table\:tempa{\:frL}
\Table\:tempb{\:frR}
\:tempa(0,0){\:left}
\:tempb(0,0){\:right}  
\IF \GtInt(\:dpR,\:dpL) \THEN
   \:both(\:dpL)  \:DoReg=\:dpL  \advance\:DoReg by 1\relax
\:tempb(\:DoReg,\:dpR){\:copy}        
   \let\:dpL=\:dpR
\ELSE
   \:both(\:dpR)
\IF  \GtInt(\:dpL,\:dpR) \THEN
   \:DoReg=\:dpR    \advance\:DoReg by 1\relax
   \:tempa(\:DoReg,\:dpL){\:copy}
\FI  
\FI      \:Temp } \Define\:copt(3){\::copt(#3,)}
\Define\::copt(2){
   \IF \EqText(#1,C) \THEN
       \let\:MergTree=\::MergTree  \FI }
\Define\:MeasureRow(3){
   \:TrcDiag{{      \:tk={#1,#2,#3}
       \immediate\write16{\Val\:C..\Val\:Ca\space\space\the\:tk}}}
   \IF  \EqInt(#2,0) \THEN
         \FigSize\Va\Vb{\:LeafNode(x,#3)}
   \ELSE \FigSize\Va\Vb{\:TreeNode(x,#3)} \FI
   \expandafter\expandafter\expandafter\:tempb
      \expandafter{\:tempa &}
   \edef\:tempa{\the\:tk #2,\Val\Va,\Val\Vb}
   \IF \LtDec(\:Vc,\:tempc)  \THEN  \:Vc=\:tempc; \FI  }
\Define\:AddVerTreeSp(1){
   \:Ca+1; \Indirect<y.\Val\:Ca>+#1;
   \Va=#1; \Va*\:TreeEdgeSpec; \Indirect<x.\Val\:Ca>+\Va;   }

\:DefineExt\TreeSpace(3){\endgroup
   \IF \EqText(D,#1) \THEN
      \def\:OutOutx{OUT:}   \def\:InInx{IN:}  \FI
   \IF \EqText(C,#1) \THEN
   \def\:OutOutx{OUT:}   \def\:InInx{INc:}
   \let\:Cspace=\:InInx  \FI 
   \IF \EqText(S,#1) \THEN
      \def\:OutOutx{OUTx:}  \def\:InInx{INx:}   \FI
   \:edef\:HorTreeSp{#2}
   \:edef\:VerTreeSp{#3}  }

\TreeSpace(S,10,20)\:DefineExt\AdjTree(1){\endgroup
   \IF \EqText(#1,) \THEN \Table\:TrSp{,,}
   \ELSE                  \Table\:TrSp{#1} \FI }

\AdjTree()\Define\:unl(3){
   \IF \EqText(#1,L) \THEN
      \IF \EqInt(\:Cc,#2) \THEN \TreeSpace(#3,0) \FI
   \ELSE \IF \EqText(#1,N) \THEN
      \:edef\:Temp{\noexpand\EqText(\Val\:Cc..\Val\:Ca,#2)}
      \IF \:Temp \THEN \TreeSpace(#3,0) \FI
   \FI \FI}\Define\:spbn(3){
   \IF \EqText(#1,B) \THEN
      \:edef\:Temp{\noexpand\EqText(\Val\:Cc..\Val\:Cb,#2)}
      \IF \:Temp \THEN  \:MoveH(#3,#3)
   \FI\FI}\Define\:lcalg(3){
   \IF \EqText(#1,A) \THEN
      \:edef\:Temp{\noexpand\EqText(\Val\:C..\Val\:Ca,#2)}
      \IF \:Temp \THEN
         \:xalg(#3)
   \FI\FI}
\Define\:xalg(2){ \def\:EntryX{#1} \:AlignTree(#2)  }
\Define\:mvpr(3){
   \IF \EqText(#1,M) \THEN
      \:edef\:Temp{\noexpand\EqText(\Val\:C..\Val\:Ca,#2)}
      \IF \:Temp \THEN  \:MoveH(#3,-#3)
   \FI\FI}\Define\TreeEdgeSpec(1){\:edef\:TreeEdgeSpec{#1}}
\TreeEdgeSpec(0.5) \def\:ChkInv{\MoveToInv}

                              \Define\:Temp{
\def\:MvToDivider##1..##2..{
   \ifx H\:TreeDir   \Indirect{\Move(\Val}<x.##1>,0)
   \else             \Indirect{\Move(0,-\Val}<x.##1>)  \fi  }
                              }\:Temp \Define\TreeEdge(2){
   \ifx H\:TreeDir   \HHEdge(#1,#2,EdgeGuide)
   \else             \VVEdge(#1,#2,EdgeGuide)  \fi }
\Define\:AddNodesAddr(1){
   \Define\:temp{\gdef\:GetNodesAddr}
   \expandafter\expandafter\expandafter
      \:temp\expandafter{\:GetNodesAddr & #1  }}

\Define\:SubTreeAddrs{
   \Table\:AddrPoly{ \:GetNodesAddr}
   \:AddrPoly(\Val\:AddrCount,999){\:AddrCount+1;\:RecordAddr}  }

\Define\:RecordAddr(1){\MoveToLoc(#1)     \MarkXLoc(#1)
      \MoveToLoc(#1;:11)  \MarkXLoc(#1;:11)}

\IntVar\:AddrCount\TreeSpec(\Node)()(\TreeEdge) 
\fi
\let\if:option=\:undefined
\let\:CheckOption=\:undefined
\let\:GetOptions=\:undefined
\let\AlDraTex=\:undefined 
\def\:aldwarn#1#2{\immediate
   \write16{---AlDraTeX warning--- \ifcase #1
      insufficient space for arrow head   
   \or \string\Tree(too many values)(...) 
   \or     edge from node #2 to itself    
   \or     #2 not a grid point            
   \or     too many entries in (#2)       
   \or     too many rows: \Val\:C         
   \or     row \Val\:C\space has too many columns: \Val\:Ca    
   \or     \Val\:C\space(>5) entries in \string\BarChartSpec   
   \or     more than one bar in \string\BarChartSpec(\the\:tk) 
   \or     no font for \string\NewCIRCNode, trying `\string
           \CIRC=lcircle10 scaled\string\magstep5'   
   \or     \string\PutNode(#2)?                                  
   \or     \string\GridDiagram(#2)? less than 2 rows/columns     
                                   \fi}}

\:RestoreCatcodes

 \input amsppt.sty

 \newif\ifproofing \proofingfalse 
 \newcount\refno         \refno=0
 \def\SetRef#1{\item{[\csname#1\endcsname]}}
 \def\MakeRef#1{\advance\refno by 1 \expandafter\xdef
        \csname#1\endcsname{{
        \ifproofing#1\else\number\refno\fi}}}
 \let\ref\MakeRef
 \newbox\keybox \setbox\keybox=\hbox{[18]\enspace}
 \newdimen\keyindent \keyindent=\wd\keybox
 \def\references{\vskip-\smallskipamount
  \bgroup
   \let\ref\SetRef
   \eightpoint   \frenchspacing
   \parindent=\keyindent
   \parskip=\smallskipamount
  }
 \def\endreferences{\egroup}
 \def\serial#1#2{\expandafter\def\csname#1\endcsname ##1 ##2 ##3
        {\unskip, \egroup #2 {\bf##1} (##2), ##3.}}
  \def\paper{\unskip, \bgroup}
  
  \def\inbook#1\bookinfo#2\publ#3\yr#4\pages#5
   {\unskip, \egroup in ``#1\unskip,'' #2\unskip, #3\unskip,
   #4\unskip, pp.~#5.}
 \def\at.{.\spacefactor3000}

 \serial{acta}{Acta Math.}
 \serial{asens}{Ann. Sci. \'Ecole Norm. Sup.}
 \serial{ca}{Comm. Alg.}
 \serial{claquagra}{Class. Quantum fravity}
 \serial{CR}{C. R. Acad. Sci. Paris}
 \serial{cjm}{Can. J. Math.}
 \serial{comp}{Compositio Math.}
\serial{commatphys}{Commun. Math. Phys.}
 \serial{ja}{J. Algebra}
\serial{jag}{J. Algebraic Geometry}
 \serial{jlms}{J. London Math. Soc.}
 \serial{jmku}{J. Math. Kyoto Univ.}
 \serial{ma}{Math. Ann.}
 \serial{mathz}{Math. Z.}
\serial{matsca}{Math. Scand.}
\serial{nucphy}{Nucl. Phys.}
\serial{phrele}{Phys. Rev. Lett.}
\serial{phylett}{Phys. Lett.}
 \serial{plms}{Proc. London Math. Soc.}
\serial{prma}{Proc. Math}
\serial{prs}{Proc. Roy. Soc.}
 \serial{slnm}{Springer Lecture Notes in Math.}
\serial{tams}{Trans. Amer. Math. Soc.}
 \serial{top}{Topology}


 \ref{AlbanoKatz}
 \ref{AsGrMo}
 \ref{AsGrMo2}
 \ref{Atiyah}
 \ref{Batyrev}
 \ref{BatyrevBorisov}
 \ref{BatyrevStraten}
 \ref{Borisov}
 \ref{CandDaLuSc}
 \ref{CandGrHu}
 \ref{CandLuSc}
 \ref{DifrKiMi}
 \ref{ElSt}
 \ref{Fulton}
 \ref{Fultoni}
 \ref{GZK}
 \ref{GrH}
 \ref{GrHu}
 \ref{GrHuLu}
 \ref{Gross}
 \ref{Harris}
 \ref{Hartshorne}
 \ref{HoKlThYa}
 \ref{HoKlThYato}
 \ref{JohKlei}
 \ref{Katz}
 \ref{Katztre}
 \ref{Katzfire}
 \ref{KlSp}
 \ref{LiTian}
 \ref{Morrison}
 \ref{MumSuo}
 \ref{Piene}
 \ref{Reid}
 \ref{Reidto}
 \ref{Sernesi}
 \ref{Sommervoll}
 \ref{Sommervollpen}
 \ref{Sommervollpto}
 \ref{Stromme}
 \ref{TiYa}
 \ref{Yau}

\newif\ifdates 

\let\:=\colon
\def\fr{{\bf f}}        

\def\O{{\cal O}}

\let\texttilde=\~
\def\~{\ifmmode\widetilde \else\texttilde \fi}
 \def\onto{\to\mathrel{\mkern-15mu}\to}
\def\Onto{\mathrel{\setbox0=\hbox{$\longrightarrow$}%
        \hbox to\wd0{$\relbar\hss\onto$}}}

\def\today{\ifcase\month\or
 January\or February\or March\or April\or May\or June\or
 July\or August\or September\or October\or November\or December\fi
 \space\number\day, \number\year}

   \newcount\date \date=\year \advance\date by -1900
        \multiply\date by 100 \advance\date by \month
        \multiply\date by 100 \advance\date by \day


\define\pen{\Bbb P^1}

\define\p{\Bbb P^3}

\define\pnp#1{\Bbb P^{n_#1}}
\define\kk#1{\varOmega^1_{\Bbb P^{n_#1}}}
\define\pr#1{\text{pr}_#1^\ast}
\define\oo#1#2{\Cal O_{#1}({#2})}
\define\open#1{\Cal O_{\pen}({#1})}
\define\oto#1#2{\Cal O_{}({#1})\oplus\Cal O_{}({#2})}
\define\otopen#1#2{\Cal O_{\pen}({#1})\oplus\Cal O_{\pen}({#2})}
\define\mul{\Bbb P^{n_1}\times \dots \times \Bbb P^{n_k}}
\define\gd#1{g_{#1}}

\define\pa#1{\pi_2(#1)}
\define\pb#1{\pi_2^{-1}(#1)}
\define\ii#1#2{\Bbb I_#1^#2}

\define\kae#1{a_{#1}\epsilon}
\define\rff#1{\lbrack {#1} \rbrack}



\Define\Nu{\hfil 0\hfil}
\Define\En{\hfil 1\hfil}
\Define\Tu{\hfil 2\hfil}
\Define\Tr{\hfil 3\hfil}
\Define\Fr{\hfil 4\hfil}
\Define\Fe{\hfil 5\hfil}
\Define\Se{\hfil 6\hfil}
\Define\Sy{\hfil 7\hfil}
\Define\Te{\hfil 8\hfil}
\Define\Ni{\hfil 9\hfil}

\Define\nu{\hfil 0\hfil \ }
\Define\en{\hfil 1\hfil \ }
\Define\tu{\hfil 2\hfil \ }
\Define\tr{\hfil 3\hfil \ }
\Define\fr{\hfil 4\hfil \ }
\Define\fe{\hfil 5\hfil \ }
\Define\se{\hfil 6\hfil \ }
\Define\sy{\hfil 7\hfil \ }
\Define\te{\hfil 8\hfil \ }
\Define\ni{\hfil 9\hfil \ }

\NewNode(\RectDiNode,\MoveToRect){   \Units(1pt,1pt)
   \GetNodeSize   \SetMinNodeSize
   \Move(-\Val\Va,-\Val\Vb)\Move(-1.0,-1.0)   \Va*2;   \Vb*2;
   \Line(0,\Val\Vb) \Line(0,2.0) \Line(\Val\Va,0)\Line(2.0,0)
   \Line(0,-\Val\Vb)\Line(0,-2.0) \Line(-\Val\Va,0)\Line(-2.0,0)
   \Move(10.0,0) \Line(0,\Val\Vb) \Line(0,2.0)
   \Move(2.0,0) \Line(0,-\Val\Vb)\Line(0,-2.0)}
\NewNode(\RectDiNodex,\MoveToRect){   \Units(1pt,1pt)
   \GetNodeSize   \SetMinNodeSize
   \Move(-\Val\Va,-\Val\Vb)\Move(-1.0,-1.0)   \Va*2;   \Vb*2;
   \Line(0,\Val\Vb) \Line(0,2.0) \Line(\Val\Va,0)\Line(2.0,0)
   \Line(0,-\Val\Vb)\Line(0,-2.0) \Line(-\Val\Va,0)\Line(-2.0,0)
   \Move(14.0,0) \Line(0,\Val\Vb) \Line(0,2.0)
   \Move(2.0,0) \Line(0,-\Val\Vb)\Line(0,-2.0)}

\NewNode(\RectDiNodea,\MoveToRect){   \Units(1pt,1pt)
   \GetNodeSize   \SetMinNodeSize
   \Move(-\Val\Va,-\Val\Vb)\Move(-1.0,-1.0) \Q=\Va; \R=\Vb;  \Va*2;   \Vb*2;
   \Line(0,\Val\Vb) \Line(0,2.0) \Line(\Val\Va,0)\Line(2.0,0)
   \Line(0,-\Val\Vb)\Line(0,-2.0)\Line(-1.0,0)
   \Line(-\Val\Q,0)\Move(-3.5,0)\Line(0,-8.0)\Move(0,8)\Move(3.5,0)
   \Line(-\Val\Q,0)\Line(-1.0,0)
   \Move(10.0,0) \Line(0,\Val\Vb) \Line(0,2.0)
   \Move(2.0,0) \Line(0,-\Val\Vb)\Line(0,-2.0)
   \Move(-12,0)\Move(-3,-3)\DrawCircle(4)}

        \NewNode(\RectDiNoded,\MoveToRect){   \Units(1pt,1pt)
   \GetNodeSize   \SetMinNodeSize
   \Move(-\Val\Va,-\Val\Vb)\Move(-1.0,-1.0) \Q=\Va; \R=\Vb;  \Va*2;   \Vb*2;
   \Line(0,\Val\Vb) \Line(0,2.0) \Line(\Val\Va,0)\Line(2.0,0)
   \Line(0,-\Val\Vb)\Line(0,-2.0)\Line(-1.0,0)
   \Line(-\Val\Q,0)
   \Line(-\Val\Q,0)
   \Move(10.0,0) \Line(0,\Val\Vb) \Line(0,2.0)
   \Move(2.0,0) \Line(0,-\Val\Vb)\Line(0,-2.0)
   \Move(-12,0)\Move(-3,-3)\DrawCircle(4)}

        \NewNode(\RectDiNodee,\MoveToRect){   \Units(1pt,1pt)
   \GetNodeSize   \SetMinNodeSize
   \Move(-\Val\Va,-\Val\Vb)\Move(-1.0,-1.0) \Q=\Va; \R=\Vb;  \Va*2;   \Vb*2;
   \Line(0,\Val\R)\Line(0,1)\Move(0,-3.0)\Line(-8,0)\Move(8,0)\Move(0,3.0)
   \Line(0,\Val\R)\Line(0,1)  \Line(1,0)\Line(\Val\Va,0)\Line(1.0,0)
   \Line(0,-1)\Line(0,-\Val\Vb)
   \Line(-1,0)\Line(-\Val\Q,0)
   \Line(-\Val\Q,0)\Line(-1,0)
   \Move(10.0,0) \Line(0,\Val\Vb) \Line(0,1.0)
   \Move(2.0,0) \Line(0,-\Val\Vb)\Line(0,-1.0)
   \Move(-12,0)\Move(-3,-3)\DrawCircle(4)}

   \NewNode(\RectDiNodef,\MoveToRect){   \Units(1pt,1pt)
   \GetNodeSize   \SetMinNodeSize
    \Move(-\Val\Va,-\Val\Vb)\Move(-1.0,-1.0) \Q=\Va; \R=\Vb;  \Va*2;   \Vb*2;
   \Line(0,\Val\R)\Line(0,1)\Move(0,-3.0)\Line(-8,0)\Move(8,0)\Move(0,3.0)
         \Line(0,\Val\R) \Line(0,2.0)
         \Line(\Val\Va,0)\Line(2.0,0)
   \Line(0,-\Val\Vb)\Line(0,-2.0)\Line(-1.0,0)
   \Line(-\Val\Q,0)\Move(-4,0)\Line(0,-8.0)\Move(0,8)\Move(4,0)
   \Line(-\Val\Q,0)\Line(-1.0,0)
   \Move(10.0,0) \Line(0,\Val\Vb) \Line(0,2.0)
   \Move(2.0,0) \Line(0,-\Val\Vb)\Line(0,-2.0)
   \Move(-12,0)\Move(-3,-3)\DrawCircle(4)}

\headline={\hss\hbox{\rm ---\ \folio \ ---}\hss}

 \newif\ifproofing \proofingfalse 
 \newcount\refno         \refno=0
 \def\SetRef#1{\item{[\csname#1\endcsname]}}
 \def\MakeRef#1{\advance\refno by 1 \expandafter\xdef
        \csname#1\endcsname{{
        \ifproofing#1\else\number\refno\fi}}}
 \let\ref\MakeRef
 \newbox\keybox \setbox\keybox=\hbox{[18]\enspace}
 \newdimen\keyindent \keyindent=\wd\keybox
 \def\references{\vskip-\smallskipamount
  \bgroup
   \let\ref\SetRef
   \eightpoint   \frenchspacing
   \parindent=\keyindent
   \parskip=\smallskipamount
  }                                             
 \def\endreferences{\egroup}
 \def\serial#1#2{\expandafter\def\csname#1\endcsname ##1 ##2 ##3
 
       {\unskip, \egroup #2 {\bf##1} (##2), ##3.}}
  \def\p {paper{\unskip, \bgroup}}
   
  \def\inbook#1\bookinfo#2\publ#3\yr#4\pages#5
   {\unskip, \egroup in ``#1\unskip,'' #2\unskip, #3\unskip,
   #4\unskip, pp.~#5.}
 \def\at.{.\spacefactor3000}

\ref{Atiyah} 
\ref{CandDaLuSc} 
\ref{CandGrHu} 
\ref{CandLuSc} 
\ref{Cl}     
\ref{GrH}  
\ref{GrHu}  
\ref{GrHa}   
\ref{JohKlei} 
\ref{Kleim} 
\ref{Katz} 
\ref{Kley}
\ref{Kn} 
\ref{Ko} 
\ref{Mo}
\ref{Nijsse} 
\ref{O} 
\ref{So}

 \title{\hbox{\vtop{\hbox{Existence of isolated curves on  $K3$-fibered 
 Calabi--Yau varieties }
\vtop{\hbox{With an application to Calabi--Yau varieties in projective space}
}}}
}\endtitle

\centerline{ISOLATED RATIONAL CURVES ON $K3$-FIBERED CALABI--YAU THREEFOLDS}
 
\centerline{by Torsten Ekedahl, Trygve Johnsen, and Dag Einar
Sommervoll} 
 
\subheading{Addresses}

Torsten Ekedahl, Department of Mathematics, Stockholm University,
S-106 91  

Stockholm, SWEDEN. E-mail: {\it teke\@matematik.su.se}

Trygve Johnsen, Department of Mathematics, University of Bergen,
N-5008 

Bergen, NORWAY. E-mail: {\it johnsen\@mi.uib.no}

Dag Einar Sommervoll, Statistics Norway, Dep N-0033 Oslo, NORWAY.

E-mail: {\it des\@ssb.no}


\subheading{Abstract}
 In this paper we study 16 complete intersection $K3$-fibered Calabi--Yau variety types  
 in biprojective
 space $\Bbb P^{n_1}\times \Bbb P^1$. 
 These are all the CICY-types that are $K3$ fibered by the projection on the
 second factor.
 We prove existence of isolated rational curves of bidegree $(d,0)$ for every
 positive integer $d$ on a general
 Calabi--Yau variety of these types. 
 The proof depends heavily on existence theorems for curves on
 $K3$-surfaces proved by S. Mori and K. Oguiso. Some of these varieties are related to
 Calabi--Yau varieties in projective space by a determinantal contraction, and we
 use this to prove existence of rational curves of every degree for a general 
 complete intersection Calabi--Yau variety in projective space.

\subheading{Subject Classification}
Primary: 14J30. Secondary: 14J28, 14H45.

 \subheading{ 1. Introduction}
 \vskip.3cm
 A complete intersection $F$ in a multiprojective space $\mul$ 
 is given by
polynomials $f_1,\cdots, f_m$. 
Each $f_j$ has a  multidegree $(g_{1j}, \cdots, g_{kj})$.
We can visualize these numerical data of $F$ by the following:

$$
\Draw
\NewNode(\RectDiNodeo,\MoveToRect){   \Units(1pt,1pt)
   \GetNodeSize   \SetMinNodeSize
   \Move(-\Val\Va,-\Val\Vb)\Move(-1.0,-1.0)   \Va*2;   \Vb*2;
   \Line(0,\Val\Vb) \Line(0,2.0) \Line(\Val\Va,0)\Line(2.0,0)
   \Line(0,-\Val\Vb)\Line(0,-2.0) \Line(-\Val\Va,0)\Line(-2.0,0)
   \Move(18.0,0) \Line(0,\Val\Vb) \Line(0,2.0)
   \Move(2.0,0) \Line(0,-\Val\Vb)\Line(0,-2.0)}
\RectDiNodeo(a)(--$\matrix\format\c&\quad\c&\quad\c&\quad\c&\quad\r\\ n_1&& g_{11}&\cdots& g_{1m}\\ \vdots&& \vdots&&\vdots \\ n_k&& g_{k1}&\cdots& g_{km}\endmatrix$--) 
\EndDraw
\tag0$$
\vskip.3cm
We use $\lbrack n\vert\vert g_{ij}\rbrack $ as a short form for this.
 We say that $F$ is of type $\lbrack n\vert\vert g_{ij}\rbrack $. 
 
 By a
{\it general} $F$ of type 
$\lbrack n\vert\vert g_{ij}\rbrack $
we understand a generic choice of the defining polynomials
$f_j$ of multidegree $(g_{1j},\cdots g_{kj})$. 
 A nonsingular complete
 intersection is Calabi-Yau, abbreviated {\it CY}, if and only if 
 $$n_i +1=\sum_{j=1}^m g_{ij}$$
 for every $i\in\{1, \dots, k\}$.
 We then say that the threefold is {\it of CICY-type}, where {\it CICY}
  is an abbreviation for complete intersection Calabi-Yau threefold. We shall
 use the notion {\it of CICY-type} for singular varieties as well.
A permutation
of the polynomials $f_j$ will in general change the matrix 
$\lbrack g_{ij}\rbrack $, but define the same variety.
    We are not interested in
distingushing between the different possibilities of indexing
of polynomials and projective spaces, and say that
$\lbrack n\vert\vert g_{ij}\rbrack $ and 
$\lbrack n^\star\vert\vert g^\star_{ij}\rbrack $
 represent the same type if one can go from one to the
other by row interchanges on $\lbrack n\vert\vert g_{ij}\rbrack $
 and column interchanges on $\lbrack g_{ij}\rbrack $.

 Let $C\subset \Bbb P^{n_1} \times \dots \times \Bbb P^{n_k} $
 be an irreducible rational curve, and let  $f\co \pen\longrightarrow C$ be a
parametrisation. We have a natural notion of multidegree
$(d_1,\dots,d_k) $ of the curve $C$ by defining $d_i$ to be the
unique number that satisfies $ \oo{\Bbb P^1}{d_i}\cong(\pi_i\circ f)^\star
\oo{\Bbb P^{n_i}}{1}$ where $\pi_i$ denotes the projection on factor $i$.

Let $Z\subseteq F$ be a closed subscheme.
  The tangent space of the Hilbert scheme $\text{\rm Hilb}_{F}$ at
  the point $\lbrack Z\rbrack $ is isomorphic to $H^0(Z,\Cal N_{Z/F})$.

  \definition{ Definition 1.1 }
Let $Z$ be a point in  $\text{Hilb}_{F}$. If     $H^0(Z,\Cal N_{Z/F})=0$,
  we say that $Z$
  is {\it isolated} in $F$.
  \enddefinition

        Note that a subscheme $Z$ is isolated if the tangent space
        at the point of the Hilbert scheme corresponding to $Z$ is
        zero dimensional.
  Recall also that $Z$ on $F$ is isolated if $Z$ does not deform to first order
on $F$.
   \remark{Remark 1.2}
Let $C$ be a nonsingular rational curve on a smooth CICY threefold $F$. 
The normal sheaf $\Cal N_{C/F}$ is
  a rank $2$ locally free sheaf on $C\cong\pen$. 
 The pullback of this sheaf to $\pen$ is isomorphic to $\oto{a}{b}$,
 where $\Cal O(a)$ is short for  $\Cal O_{\pen}(a)$.
 Since
$F$ is CY we get $a+b=-2$ from the exact sequence
$$\/0 \longrightarrow T_C \longrightarrow T_{F} \vert_C
\longrightarrow\Cal N_{C/ F }   \longrightarrow  0, $$
where $T_C$ is the tangent bundle of $C$ 
and  $T_F\vert_C$ is the tangent bundle of $X$
 restricted to $C$.
Hence C is isolated if and only if $\Cal N_{C/F}\cong\oto{-1}{-1}$.
\endremark
  \vskip.2cm

The following theorem is the main result of the paper.

\proclaim{Theorem 1.3}
There are exactly $16$ types of CICY's in biprojective spaces of
the form $\Bbb P^n\times\Bbb P^1$, for positive integers $n$.
For a general (and hence smooth) CICY of each of these $16$
types, and for all positive integers $d$, there exists an isolated smooth rational
curve $C$ of bidegree $(d,0)$.
\endproclaim


Theorem 1.3 has a nice corollary.

\proclaim{Corollary 1.4}
Let $F$ be a general member  of one of the $5$ well-known  CICY-types in projective spaces
$\Bbb P^n$, for positive integers $n$. Then for all positive
integers $d$, there exists a smooth isolated rational curve $C$ of degree $d$
 in $F$.
\endproclaim

\remark{Remark 1.5}
  
While we were working with the material in this paper, we obtained the
manuscript ([12]). It is clear that our Corollary 1.4 can be obtained
from the results of that paper.  Theorems 5.1 and 5.2 of ([12]), in
combination with Oguiso/Mori's existence theorem on curves on
K3-surfaces, give the conclusion of our Corollary 1.4. Our approach
is independent of Kley's, but share the idea, going back to Clemens,
of using K3-surfaces to study curves on CICY's.

\endremark

Our paper is organized as follows. In paragraph 2 we recall some basic facts about $K3$-surfaces that will
be essential in proving our result. The paragraph ends with Proposition 2.4,
which will be an important tool in paragraph 4. 
In paragraph 3 we prove the first part of Theorem 1.3. We give a list
of the $16$ types of CICY-threefolds in biprojective spaces 
$\Bbb P^n\times \Bbb P^1$ for positive integers $n$. Furthermore,
 we show how one can
reduce the last part of Theorem 1.3 to showing the existence of one 
smooth isolated rational bidegree $(d,0)$ 
curve $C$ on one special CICY of each type.
 
In paragraph 4 we complete the proof of Theorem 1.3. We use the results
from paragraph 3, and we show that for each positive integer $d$, and each
of the $16$ CICY-types there exists an isolated smooth rational curve of bidegree $(d,0)$ contained in a CICY of that type.
The main point is to construct suitable $F$ that are $\Bbb P^1$ fibrations
of $K3$-surfaces that are complete intersections in projective spaces. We
use a result by S. Mori and K. Oguiso (for a more general statement,
see Theorem 8.1 of \cite\Kn )that gives the existence of  a smooth
rational curve $C$ of each degree, such that each $C$ is isolated in
a $K3$-surface $S$. 
For each of the $16$ types we start with such a pair $(C,S)$, and we
show that it is possible to choose the $\Bbb P^1$-fibration such that
$F$ is smooth, and $C$ is isolated, not only in $S$, but also in $F$.
In this analysis the material from paragraph 2 will be essential.
In paragraph 5 we describe  determinantal contractions of CICY`s.
We use this description to show how Corollary 1.4 follows from
Theorem 1.3. 

Corollary 1.4. was proved by H. Clemens for infinitely many $d$
for the quintic in $\Bbb P^4$. S. Katz (\cite\Katz) gave a proof for the
quintic in $\Bbb P^4$ for all $d$ and hence proved Corallary 1.4
in this case. K. Oguiso proved Corollary 1.4 in case of the
complete intersection $(4,2)$ in $\Bbb P^5$. Our proof is inspired by (\cite\O).
Results like Theorem 1.3 and Corollary 1.4 are interesting when one
counts rational curves on Calabi--Yau threefolds. Mirror symmetry
gives predictions on the number of rational curves  of a given
(multi)degree on a generic Calabi--Yau threefold. Nevertheless,
it remains to prove rigorously  that for a generic CICY of a
given type, and for a given (multi)degree, there exists a nonempty
set of isolated curves of this (multi)degree.  For the quintic
in $\Bbb P^4$ this is essentially Clemens' conjecture. The finiteness
has only been proved for $d\leq 9$ (\cite\JohKlei,\cite\Katz,\cite\Nijsse), while the nonemptyness is
proved by S. Katz (\cite\Katz). 
The analogue of Clemens' conjecture is wrong
for general CICY's in multiprojective spaces in general (\cite\So).
Our results are  a step towards proving the analogue of Clemens'
conjecture for complete intersections in $\Bbb P^n\times \Bbb P^1$
(and in $\Bbb P^n$).
There is also a more philosophical side of our approach.
The  moduli space of Calabi--Yau varietes may be connected. This
is popularly known as Reid's fantasy. Candelas et al. proved
connectedness for complete intersections in multiprojective spaces (\cite\CandGrHu).
The proof relies on existence of determinantal contractions that
give a map from one Calabi--Yau variety in one multiprojective space
to another Calabi--Yau in another multiprojective space. This
map is an isomorphism outside strata of codimension $2$ that are
contracted to points. Curves that do not intersect these strata
are mapped isomorphically to curves on another Calabi--Yau in 
another multiprojective space. Our philosophical point of view is
therefore that studying rational curves on a Calabi--Yau variety
in a multiprojective space is too restrictive. This paper is an
illustration of this point. Contructing special $K3$-fibrations
in biprojective spaces does not only prove existence of isolated  
rational curves of every bidegree $(d,0)$, but also proves
existence of isolated rational curves of every degree $d$ for the five 
complete intersection Calabi--Yau threefolds in projective spaces.
 
\vskip.1cm
 We thank the Institut Mittag-Leffler for providing a stimulating  atmosphere
while working on this paper, and the referee for making many good points,
including the discovery of one {\it  CICY}-type. The second author thanks S.L. Kleiman
for introducing him to ([12]), and H. Kley for useful
discussions about this work.

\subheading{ 2.     Some basic theory of $K3$-surfaces}
\vskip.3cm
Let $S$ be a $K3$-surface.
Study the following part of the long exact sequence obtained from the
well-known exponential sequence:
                                                                       
$$ 0 = H^1(S, \Cal O ) \longrightarrow   \text{Pic}(S) = H^1(S, \Cal O^\star ) \overset{c}\to\longrightarrow  H^2(S,\Bbb Z) \overset{\kappa}\to\longrightarrow  H^2(S, \Cal O) = \Bbb C.$$

All $K3$-surfaces $S$ are diffeomorphic, and  $H^2(S, \Bbb Z)$ is torsion free of
rank $22$ for all such $S$. As usual, fix once and for all an isomorphism between the lattice
  $H^2(S,\Bbb Z)$ and a standard rank 22 lattice with a known intersection form.
 This gives an intersection matrix $A = \lbrack a_{ij}\rbrack \ i,j=1,..,22$
 compatible with the cup-product form on  $H^2(S,\Bbb Z)$. Then we complexify 
  $H^2(S,\Bbb Z)$ and obtain  $H^2(S,\Bbb C)$. We identify  $\Bbb P^{21}$ with 
 $\Bbb P(H^{2}(S,\Bbb C))$.
 
Let $P$ be the well-known period map from the set of isomorphism
classes of analytical $K3$-surfaces into  $\Bbb P^{21}$. In short
each isomorphism class determines a unique holomorphic 2-form on $S$
up to multiplicative constant. Integration of this form then gives a well-defined 
element of $\Bbb P^{21}$ = $\Bbb P(H^{2}(S,\Bbb C))$, and this element is the
image of the isomorphism class by $P$. 
  We have:

\definition{Definition 2.1}
 Let the subset $D$ of  $\Bbb P^{21}$ be defined by
$$D=\{(\gamma_1,...,\gamma_{22}) | \sum_{ij} a_{ij}\gamma_i \gamma_j = 0, 
\sum_{ij}  
a_{ij}\gamma_i \bar \gamma_j > 0\}.$$
 \enddefinition

$D$ can then be viewed as a moduli space of analytical $K3$-surfaces
(Locally this holds, even for a fixed identification between $H^2(S,\Bbb Z)$ 
and the standard rank 22 lattice as  described above. Taking into account the 
possibility of different such identifications, one has to divide $D$ with 
a discrete group to obtain a moduli space of analytical $K3$-surfaces, but we shall 
only be interested in the local theory).
We see that $D$ may be viewed as ``half of a hyperquadric in $\Bbb P^{21}$''.

The following result is a consequence of Theorem 14 in \lbrack 13\rbrack :

\proclaim{Proposition 2.2}
Let $C_1, .....C_s$, for some $s$  in the range $1,...,20$, be independent elements
in $\text{Pic}(S)$, for a given fixed (analytical structure on a ) $K3$-surface $S$.
We view $\text{Pic}(S)$ as a subgroup of $H^2(S, \Bbb Z)$ via the map  $c$  above.

Let $H_i$ be the embedded tangent space at $P(S)$ in $\Bbb P^{21}$ of the 
subset of $D$ consisting of period points of $K3$-surfaces with $C_i$ in  
the images of the respective $c$-maps, for $i=1,...,s$. Then the intersection of the
$H_i$ is a  space of dimension $20-s$.
\endproclaim

\remark{Remark 2.3}
This means that the set of $K3$-surfaces that "keep $C_1, \dots,C_s$"
has dimension at most $20-s$, even infinitesimally. See also $\lbrack 8\rbrack$, page 594.
\endremark 
\vskip.5cm
\flushpar
{\bf Algebraic $K3$-surfaces in a fixed projective space}

Let us study the space $R_n$ of embedded smooth $K3$-surfaces in $\Bbb P^n$,
for $n=3,4$. For $n=3$ this is an open subset of $\Bbb P^{34}$, parametrizing
quartic surfaces in $\Bbb P^3$. For $n=4$,  $R_n$  is an open subset of the set
of complete intersections of type $(2,3)$. 
For $n=5$ we define  $R_n$ to be the space of all complete intersection
K3-surfaces (of type (2,2,2)). Then $R_5$ is dense in the set of all embedded smooth $K3$-surfaces 
in $\Bbb P^5$, and also in the space of nets of quadric hypersurfaces in $\Bbb P^5$.
In each case it is well known, and easy to prove, by a simple dimension
count, that there is an $19$-dimensional set of projective equivalence
classes of $K3$-surfaces as described. The spaces $R_n$ are smooth of
dimensions $34, 43, 54$, for $n= 3,4,5$ respectively.
See  $\lbrack 8\rbrack$, p. 591-594.
It is well known that projective equivalence  in these cases amount to the
same thing as abstract isomorphism, in the sense that two surfaces are mapped
to the same point by the period map $P$ described above.
Let $P_n$ be the pull-back of the period map $P$ to $R_n$,
and set $D_n=\text{im}(P_n)$.
\vskip.5cm
\flushpar
{\bf Maps of tangent spaces}

Let $S$ be a fixed $K3$-surface, corresponding to a point $\lbrack S\rbrack$  in $R_n$, for $n=3,4,$
or $5$. Assume that $S$ contains a smooth rational curve $C$ of degree $d$, for
some positive integer $d$.  Let $F$ inside $D_n$ be the locus of those $K3$-surfaces
that ``keep $C$" in the sense described above. Look at the following
maps:                             
  $$          T_{R_n,\lbrack S\rbrack}  \overset{dP_n}\to\longrightarrow  T_{D_n,\lbrack S\rbrack}  \overset{f}\to\longrightarrow  
T_{D_n,\lbrack S\rbrack}  /T_{F,\lbrack S\rbrack}  = \Bbb C.$$
Here $f$ is the natural quotient map. The quotient to the right is isomorphic
to $\Bbb C$, since $\text{dim} T_{D_n,\lbrack S\rbrack}  = 19$, and $\text{dim} T_{F,\lbrack S\rbrack}  = 18$, by Proposistion 2.2
(Put $C_1$ = hyperplane section of $S$ in $\Bbb P^n$, and $C_2 = C$). Clearly the composite
map  $f \circ dP_n$ is surjective and linear. Let us define and study
 natural surjective linear maps  $W_n \overset{\phi_n}\to \longrightarrow  T_{R_n,\lbrack S\rbrack} $ , where $W_n$ will be
$\Bbb A^{35},  \Bbb A^{15} \times A^{35}, (\Bbb A^{21})^3 $, for $n = 3, 4, 5$, respectively.
$W_3$ will  parametrize homeogeneous quartic polynomials in $4$ variables,
$W_4$ will  parametrize pairs of homogeneous quadric and cubic polynomials in
$5$ variables, and $W_5$ will  parametrize triples of homeogeneous quadric
polynomials in $6$ variables.
Let $S$ be a $K3$-surface in $\Bbb P^3$ with equation $q_4(x_0,...,x_3) = 0$. We  construct
a surjective linear map  $W_3 \longrightarrow  T_{R_3,\lbrack S\rbrack} $ in the following manner,

          $$
           \align \phi_3\co   p_4(x_0,...,x_3) \longrightarrow 
\text{``the K3-surface with equation":}&\\
   \  \      q_4(x_0,...,x_3) + tp_4(x_0,...,x_3)= 0&.\\
           \endalign
           $$
        
Here the polynomial $q_4(x_0,...,x_3)$ is fixed. We think of $t$ as an 
infinitesimal variable with $t^2 = 0$, thus the
points of the tangent space $T_{R_3,\lbrack S\rbrack}$  can be given by 
such equations $q_4(x_0,...,x_3) + tp_4(x_0,...,x_3)= 0$
over  $C\lbrack t\rbrack /t^2$.  Clearly  $\phi_3$  is surjective and essentially maps
$p_4(x_0,...,x_3)$ to its image in the vector space of all quartic polynomials in 4 variables modulo $q_4(x_0,...,x_3)$.
Let $S$ be a $K3$-surface in $\Bbb P^4$ with equation  $f_2(x_0,...,x_4) = f_3(x_0,...,x_4) = 0$.
In an analogous way we  construct a linear map  $W_4 \longrightarrow  T_{R_4,\lbrack S\rbrack} $ :

$$
\align 
  \phi_4\co   (p_2(x_0,...,x_4), p_3(x_0,...,x_4))  \longrightarrow   \text{ ``the K3-surface with equation":}&\\
   f_2(x_0,...,x_4) +tp_2(x_0,...,x_4) =0&\\
   f_3(x_0,...,x_4) + tp_3(x_0,...,x_4)= 0&.
\endalign
$$
Clearly  $\phi_4$  is surjective and maps $(p_2(x_0,...,x_4), p_3(x_0,...,x_4))$  to the
pair of images  modulo $f_2(x_0,...,x_4)$ and $(f_2(x_0,...,x_4),f_3(x_0,...,x_4))$,
respectively .
Clearly $T_{R_4,\lbrack S\rbrack}$ can be viewed as the set of such pairs of all such cosets
of quadric and cubic polynomials, respectively.
Let $S$ be a $K3$-surface in $\Bbb P^5$ with equation
   $$ g_2(x_0,...,x_5) = h_2(x_0,...,x_5) = i_2(x_0,...,x_5) = 0.$$
In an analogous manner we  construct a linear map  $W_5 \longrightarrow  T_{R_5,\lbrack S\rbrack}  $:

$$
\align
\phi_5\co   (r_2(x_0,...,x_5), s_2(x_0,...,x_5), u_2(x_0,...,x_5))  \longrightarrow   \text{ ``the K3-surface ":}&\\
g_2(x_0,...,x_5) + tr_2(x_0,...,x_5) =0&\\
 h_2(x_0,...,x_5) +ts_2(x_0,...,x_5) =0&\\
i_2(x_0,...,x_5) +tu_2(x_0,...,x_5) =0&.
\endalign
$$
Clearly  $\phi_5$  is surjective and maps $(r_2(x_0,...,x_5), s_2(x_0,...,x_5)$,
$u_2(x_0,...,x_5))$   to the triple of images  modulo
$(g_2(x_0,...,x_5),h_2(x_0,...,x_5), i_2(x_0,...,x_5) )$.
We see that $T_{R_5, \lbrack S\rbrack} $ can be viewed as the set of such triples of cosets of
quadric polynomials.  Study the composite surjective, linear map
        $$ \eta_n = f \circ dP_n \circ \phi_n\co    W_n \longrightarrow  \Bbb C, \ 
\text{for } n = 3,4,5.$$
The surjectivity of the linear map $\eta_n$ gives:

\proclaim{ Proposition 2.4}
The kernel of the linear map $\eta_n$ is a proper subspace of codimension
$1$ in $W_n$.
\endproclaim

In paragraph 5 we will produce some natural subcones of the $W_n$, arising from
 CICY's that are  $K3$-fibrations. We will prove that general points of these cones are not
contained in  $\text{ker}(\eta_n)$. It will be sufficient to prove that the
linear spans of these cones are all of $W_n$.

\subheading{3. CICY threefolds in
multiprojective spaces}
\vskip.3cm
In the first part of this paragraph we find all CICY-types in
a  biprojective space of the type $\Bbb P^n\times\Bbb P^1$.
We want to discard those CICYs that are projected to a zero dimensional
set by the second projection. Such threefolds are easily identified
with a CICY in $\Bbb P^n$ by the first projection.

\definition{Definition 3.1}
Let  $F$ be a complete intersection CY-variety in $\Bbb P^{n}\times \Bbb P^1$. If 
the projection on the second factor is surjective, then $F$ is {\it $\Bbb P^1$-surjective}. 
We also use this notion for CICY-types, i.e. a CICY-type is $\Bbb P^1$-surjective
if the general member is $\Bbb P^1$-surjective.
\enddefinition

\proclaim{Lemma 3.2}
There are $16$ $\Bbb P^1$-surjective CICY-types and they are  divided
into the following three groups:

Group 1.

 $$
 \Draw
\RectDiNode(a)(--\tr\ \Fr~~\en\ \Tu--) 
 \EndDraw
 $$

   $$
 \Draw
\RectDiNode(b)(--\fr\ \fr\En~~\en\ \nu\Tu--)
\Move(50,0)
\RectDiNode(c)(--\fr\ \fr\En~~\en\ \en\En--)
 \EndDraw
 $$

$$
 \Draw
\RectDiNode(a)(--\fe\ \fr\en\En~~\en\ \nu\en\En--)
\EndDraw
$$

Group 2.

 $$
 \Draw
\RectDiNode(a)(--\fr\ \tr\Tu~~\en\ \tu\Nu--)
\Move(50,0)
\RectDiNode(b)(--\fr\ \tr\Tu~~\en\ \nu\Tu--)
\Move(50,0)
\RectDiNode(c)(--\fr\ \tr\Tu~~\en\ \en\En--)
 \EndDraw
$$

$$
\Draw
\RectDiNode(b)(--\fe\ \tr\tu\En~~\en\ \nu\nu\Tu--)
\Move(60,0)
\RectDiNode(c)(--\fe\ \tr\tu\En~~\en\ \nu\en\En--)
 \Move(60,0)
\RectDiNode(d)(--\fe\ \tr\tu\En~~\en\ \en\nu\En--)
 \EndDraw
 $$
 
   $$
\Draw
\RectDiNode(b)(--\se\ \tr\tu\en\En~~\en\ \nu\nu\en\En--)
\EndDraw
$$

Group 3.
 $$
\Draw
\RectDiNode(b)(--\fe\ \tu\tu\Tu~~\en\ \tu\nu\Nu--)
\Move(60,0)
\RectDiNode(c)(--\fe\ \tu\tu\Tu~~\en\ \en\en\Nu--)
\EndDraw
$$

  $$
\Draw
\RectDiNode(b)(--\se\ \tu\tu\tu\En~~\en\ \nu\nu\nu\Tu--)
\Move(70,0)
\RectDiNode(c)(--\se\ \tu\tu\tu\En~~\en\ \En\nu\nu\En--)
\EndDraw
$$
      $$
\Draw
\RectDiNode(b)(--\sy\ \tu\tu\tu\en\En~~\en\ \nu\nu\nu\en\En--)
\EndDraw
$$
\endproclaim

 \demo{Proof}
We want to determine:

$$
\Draw
\NewNode(\RectDiNodeo,\MoveToRect){   \Units(1pt,1pt)
   \GetNodeSize   \SetMinNodeSize
   \Move(-\Val\Va,-\Val\Vb)\Move(-1.0,-1.0)   \Va*2;   \Vb*2;
   \Line(0,\Val\Vb) \Line(0,2.0) \Line(\Val\Va,0)\Line(2.0,0)
   \Line(0,-\Val\Vb)\Line(0,-2.0) \Line(-\Val\Va,0)\Line(-2.0,0)
   \Move(18.0,0) \Line(0,\Val\Vb) \Line(0,2.0)
   \Move(2.0,0) \Line(0,-\Val\Vb)\Line(0,-2.0)}
\RectDiNodeo(a)(--$\matrix\format\c&\quad\c&\quad\c&\quad\c&\quad\r\\ n_1&& g_{11}&\cdots& g_{1m}\\ \vdots&& \vdots&&\vdots \\ n_k&& g_{k1}&\cdots& g_{km}\endmatrix$--) 
\EndDraw
\tag0$$
\vskip.3cm

\noindent 
where $k=2$, and $\sum g_{1j}= n_1+1$ and $\sum g_{2j}= 2$, and $n_2 = 1$.
Note that a column with only one 1 and rest zeros gives a hyperplane in one of
the factors and thus gives a natural identification with a CICY-type in a lower
dimensional embedding space.  Discarding those, gives
 $n_1\leq 7$. To see this, take  first $n_1=8$. Since we have
$6$ polynomials in $\Bbb P^8\times \Bbb P^1$, at least
three of them must have degree $1$ in the first factor. Furthermore,
since we have at most  two nonzero entries on the second factor,
at least one of the 1's in the first row must have a zero in its
column. Thus the variety can be idenfied with a variety of a
CICY-type a lower dimensional embedding space. If $n_1>8$ then the
first row must contain more than three 1's and the same argument
applies.

Since we are looking at $3$-dimensional CY-varieties, we also have $n_1\geq 3$.

Computation yields the list given above.
 \enddemo

 \remark{ Remark 3.3}
The reason for the partition of the $16$ CICY-types  into 
three groups will become apparent in paragraph 4.  Note also
that all these CICY-types are $K3$-fibrations over $\Bbb P^1$.
In the following we will also adopt  a nonfigurative notation for the 
different CICY-types.  For example we will write $X_{4\vert 41,1\vert 11}$
for 
$$\Draw
\RectDiNode(c)(--\fr\ \fr\En~~\en\ \en\En--)
 \EndDraw
 $$
\endremark  
 \vskip.3cm
  We close this paragraph by proving the following theorem, which
  will be crucial in paragraphs 4 and 5.
 
  \flushpar\proclaim{Theorem 3.4}
  Let $(d_1,\dots,d_k)\in \Bbb N^k$
 and not all  $d_i=0$. Suppose there exists a CI threefold $F $
 of type $\lbrack n\vert\vert g_{ij}\rbrack$ in 
 ${\Bbb P^{n_1}\times \dots \times \Bbb P^{n_k}} $
 containing an isolated nonsingular rational curve $C $ with
multidegree $(d_1,\dots,d_k)$. Assume that $F$ is smooth at
all points of $C$,
 then there exist isolated rational curves
 with this multidegree  on a general CI of type
$\lbrack n\vert\vert g_{ij}\rbrack$.
\endproclaim
\flushpar\demo{Proof}
 Let $\Bbb M$ parametrise the space of all smooth rational curves
 of multidegree
  $(d_1, \dots ,d_k)$, and and let $G\subseteq\text{Hilb}(\Bbb P^{n_1}\times \dots \Bbb P^{n_k})$ be
the parameter space of CI threefolds of type $\lbrack n\vert\vert
g_{ij}\rbrack$.
Let 
$$\Bbb J_0 =\lbrace (C,F) \in \Bbb M \times G \vert 
\text{ $C\subseteq \ F$, and $F$ is nonsingular along $C$ } \rbrace $$
The incidence variety $\Bbb J_0$ is equipped with projection maps
$$\pi_1\co \Bbb J_0 \longrightarrow M\quad
\text{and }\pi_2\co \Bbb J_0 \longrightarrow G.$$
A standard dimension count gives
$$\text{dim } \Bbb J_0 \geq \text{dim } G$$

Let $(C,F)\in \Bbb J_0$ be as in the statement of the theorem, and
study the map
 $d\pi_2\co  T_{\Bbb J,(C,F)}\longrightarrow T_{G,F}$.
Then $\text{ker }d\pi_2=H^0(N_{C/F})=0$. This implies
that $d\pi_2$ is injective. It is also surjective, since
$\text{dim } \Bbb J_0\geq \text{dim } G$ implies that 
$\text{dim }T_{\Bbb J_0,(C,F)}\geq \text{dim } \Bbb J_0\geq \text{dim }  G
= \text{dim } T_{G,F}$ ($G$ is smooth at $F$).
The surjectivity of $d\pi_2$ implies that $\pi_2$ is dominant since
$G$ is irreducible. (The injectivity implies that the component of 
$\Bbb J_0$ has dimension equal to the dimension of $G$).
\enddemo

\subheading{4. Existence of isolated rational curves on $\Bbb P^1$-surjective CICYs}
\vskip.3cm
A key tool in this paragraph will be

\proclaim{Proposition 4.1}
Let  $d$ be a positive integer.  Then:
\vskip.2cm\flushpar
i.) (\cite\Mo) 
There exists a nonsingular quartic $K3$-surface $ S_1$ in $\Bbb P^3$,
containing an isolated rational curve $C_1$ of degree $d$, for every positive integer $d$.
 \vskip.2cm\flushpar
ii.) (\cite\O) 
There exists a nonsingular $K3$-surface $S_2$ of type (2,3) in $\Bbb P^4$, 
containing an isolated nonsingular rational curve $C_2$ of degree $d$,
for every positive integer $d$.
  \vskip.2cm\flushpar
iii.) (\cite\Kn) 
There exists a nonsingular $K3$-surface $S_3$ of type $(2,2,2)$ in $\Bbb P^5$,
containing an  isolated nonsingular rational curve $C_3$ of degree $d$,
for every positive integer $d$.
\endproclaim

\demo{Proof}
See \cite\Mo{} , \cite\O , and Theorem 8.1 of \cite\Kn   .
\enddemo

\remark{Remark 4.2}
One can prove that $S_2$ can be taken not only to be smooth, but also to be the intersection of a smooth
hyperquadric and a smooth hypercubic. It is easy to see that the hypercubic can be
taken to be smooth, given that $S_2$ is. From the way Oguiso constructed the surface
$S_2$, it is clear that the rank of the Picard group of this hypersurface is 2,
with the rational curve, and a hyperplane section as generators. If the hyperquadric
were singular, necessarily with a single singular point, then $S_2$ would
contain a plane cubic curve, with arithmetic genus 1.
A study of such a curve contradicts the fact that this curve must be a linear integer
combination of the generators of $\text{Pic}(S_2)$. Hence $S_2$ cannot be singular.
We thank Kristian Ranestad for pointing this out to us. It is also obvious that
$S_3$ can be taken to be the intersection of three smooth hyperquadrics. We
will assume this.
\endremark

We will use Proposition 4.1 and some additional argument to prove the main
result of this section, which is:

\proclaim {Theorem 4.3}
Let $d$ be a positive integer. For each of the $16$ $\Bbb P^1$-surjective CICY-types
 listed in Lemma 3.2, there exists a smooth $3$-fold $F$ of 
that type containing an isolated rational curve $C$ of bidegree $(d,0)$.
 \endproclaim

\demo{Proof}
Let a positive integer $d$, and one of the $16$ CICY-types be given.
We will produce a  smooth $K3$-fibration  
$F$ over $\Bbb P^1$ of this CICY-type such that the fibre 
over the point $(0,1)$ is one of the $K3$-surfaces $S_i$ of 
Proposition 4.1, containing a  curve $C_i$ of degree  $d$ for one $i\in \{ 1,2,3\}$.
In addition the fibration will be chosen in such a way that $C_i$ does not
deform  "sideways" in $F$, not even infinitesimally.
In the language of paragraph 2: the induced infinitesimal deformation
of $S_i$ corresponds  to an element outside $\text{ker } \eta_n$ for
$n=3,4,$ or $5$. Assume that such an $F$ has been found, and 
consider (following \cite\O) the exact sequence:

$$0\longrightarrow N_{C_i/S_i}=\Cal O(-2)\overset{h}\to\longrightarrow N_{C_i/F}
   \overset{g}\to\longrightarrow N_{S_i/F}\vert_{C_i}=\Cal O\longrightarrow 0$$

Since $h$ is injective, the bundle $N_{C_i/F}$ is $\Cal O_{\Bbb
P^1}(-1)\oplus \Cal O_{\Bbb P^1}(-1)$ or $\Cal O_{\Bbb P^1}\oplus
\Cal O_{\Bbb P^1}(-2)$. Assume
$N_{C_i/F}=\Cal O_{\Bbb P^1}\oplus \Cal O_{\Bbb P^1}(-2)$.  Then $
g$ gives an isomorphism 
$ g^\prime\co  H^0(N_{C_i/F})\longrightarrow H^0(N_{S_i/F}\vert_{C_i}) = \Bbb C$.
On the other hand, since $ N_{S_i/F}= \Cal O_{S_i}$, we have that the
restriction map 
$\text{Res}\co  H^0(N_{S_i/F})\longrightarrow H^0(N_{S_i/F}\vert_{C_i})= \Bbb C$
is an isomorphism. Looking at the combined isomorphism
$(g^\prime)^{-1}\circ \text{Res}\co   H^0(N_{S_i/F})\longrightarrow H^0(N_{C_i/F}),$
 we see that any non-zero  element of $ H^0(N_{S_i/F})$ is mapped
to a non-zero
infinitesimal deformation of $C_i$ along the
deformation  of $S_i$ of $F$.
This means that $H^0(N_{S_i/F})$ is contained in $\text{ker} (f\circ dP_n)$,
in the language of paragraph 2, or that $H^0(N_{S_i/F})$ is in the image
 $\phi_n(\text{ker } \eta_n)$ (Here a natural identification is made
 between $H^0(N_{S_i/\Bbb P^n})$ and $T_{R_n,\lbrack S_i \rbrack}$, and
 $ H^0(N_{S_i/F})$ is in a natural way a subset of $H^0(N_{S_i/\Bbb P^n})$).
 This is a contradiction, and hence $N_{C_i/F}$ must be $\Cal O_{\Bbb P^1}(-1)\oplus \Cal O_{\Bbb P^1}(-1)$.
In other words $C_i$ is isolated in $F$.
Hence, in order to prove the theorem it is enough to construct
$K3$-fibrations in each of the $16$ cases satisfying the following
 properties:
 \vskip.2cm\flushpar
a. The fibre over  the point $(0:1)\in \Bbb P^1$ is one of the $S_i$
   from Proposition 4.1.
 \vskip.2cm\flushpar
b. The induced element of 
   $ H^0(N_{S_i/F})\subseteq H^0(N_{S_i/\Bbb P^n})\subseteq
T_{R_n,\lbrack S_i \rbrack}$
   is not contained in $\text{ker}(f\circ dP_n)$.
 \vskip.2cm\flushpar
c. The constructed $F$ is smooth.

  \vskip.2cm              
In each of the  16 cases we start with a rational curve $C_i$ inside one of the
three $K3$-surfaces  $S_i$ ($i\in\{1,2,3\}$) constructed by Oguiso and others,  
depending on whether we are in
Group 1, 2 or 3, using the terminology of
paragraph 3. The curve $C_i$ is now isolated in $S_i$.  Let the 
 $K3$-surfaces $S_i$ be given by: 

$$
\align
   S_1=Z(q_0(x_0,\dots,x_3))&\subseteq \Bbb P^3\\
   S_2=Z(k_0(x_0,\dots,x_4),c_0(x_0,\dots,x_4))&\subseteq \Bbb P^4      \\
   S_3=Z(k_{01}(x_0,\dots,x_5),k_{02}(x_0,\dots,x_5),k_{03}(x_0,\dots,x_5))&\subseteq \Bbb P^5
\endalign 
$$ 
where $q_0$ is a polynomial of degree $4$, $k_0$ and $c_0$ are polynomials of degree
$2$ and $3$ respectively, and $k_{01},k_{02},k_{03}$ are polynomials of degree $2$.
In the following we will write $\bar x_i$ as a short form for $x_0,\dots,x_i$
in order to make the exposition more transparent.
The surface $ S_1$ will after this convention be the zero set of  $(q_0(\bar x_3))$. 
In general we will use $l,k,c,q$ as names for polynomials of degree $1,2,3,4$ respectively,
in order to make the degrees manifest in the construction. As seen above we have reserved
the subscript zero, for the polynomials defining a special $K3$ surface. A polynomial
with a subscript without a zero will unless otherwise stated be a general polynomial
of the assigned degree.
We treat the three different groups appearing in  Lemma 3.2 separately.

In 12 of the 16 cases we will see directly that the linear span of 
a relevant subset of $W_n$ is all of $W_n$. See the remark after
Proposition 2.4. These cases include the 9 cases that will be mentioned in 
Remark 5.4. Five of these cases will be used to prove Corollary 1.4 in chapter 5.
The remaining cases are irrelevant to the proof of Corollary 1.4, but in order to prove
Theorem 1.3 we will show that in the 4 cases not among the 12 mentioned, the  linear spans of the corresponding relevant 
 subsets are not contained
in $\text{ker } \eta_n$  even if they are not equal to all
of $W_n$. A useful observation  will then be: Different sets of
equations sometimes give rise to the same $K3$-surface. The sets
$f=g=0$ and $f=g+kf=0$, for homogeneous polynomials $f,g,k$ of
degrees $d,d',d'-d$, with $d \leq d'$ give the same $K3$-surface.
The first order deformation $f+tm=g+kf=0$ of the latter set
corresponds to the first order deformation $f+tm=g-tkm=0$ of the
first set. Therefore, if $(m,0) \in \text{ker } \eta_n$, then 
$(m,-km) \in \text{ker } \eta_n$ and $(0,km) \in \text{ker } \eta_n$ also,
for any $k$. 
We treat the three different groups appearing in  Lemma 3.2 separately.

Consider Group 1.
  As a typical case, we study:
 $$
 \align
 X_{4\vert 41,1\vert 02}:\hskip5cm&\\
 q_0(\bar x_3)  + x_4c_1(\bar x_4) &= 0,\\
x_4 z^2 + l_1(\bar x_4)yz + l_2(\bar x_4) y^2  &= 0,
\endalign
$$

The $K3$-surface $S_1$ now sits as the fibre of $F$ over $P = (0:1)$ in 
$\Bbb P^1 = \text{Proj }\Bbb C[y,z]$. 
A local parameter of $\Bbb P^1$ at $P$ is $t=y/z$, and letting the $K3$-surface
vary in the family, we obtain as a first order deformation (modulo $t^2$):

$$
\align
q_0(\bar x_3) - tl_1(\bar x_3,0)c_1(\bar x_3,0) &= 0.
\endalign
$$

Recall the map $\eta_n=f\circ dP_n\circ \phi_n\co W_n\rightarrow \Bbb C$ introduced
at the end of paragraph 2. We study the case $n=3$.

The deformations above correspond to applying the map $\phi_3$
(and $\eta_3$) to a subset of $W_3=\Bbb A^{35}$ of the following
form:
$$
\align
&\{l_1(\bar x_3,0)c_1(\bar x_3,0)\}
\endalign
$$
   We see that
the linear span of the subset inside $\Bbb A^{35}$  is all of this space.
Recall: Any homogeneous polynomial is in the span of all monomials of the
degree in question.
We also see that the condition that $\eta_3$, restricted to some subdomain of
$\Bbb A^{35}$ is zero, is equivalent to the condition that $\eta_3$, restricted (only) to
all of the linear span of this subdomain of $\Bbb A^{35}$ is zero.  This is true,
simply since $\text{ker}(\eta_3)$ is a linear subspace. But since the linear
span of the set in question is all of
$\Bbb A^{35}$,  Proposition 2.4 gives
that the set is not contained in the kernel of $\eta_3$. Hence,  if
$l_1, c_1$ are
chosen generically enough $C_1$ will not survive in the first order
deformations described. The three other cases in Group 1 are treated
in a practically identical way. In all cases the linear span of the
relevant subdomain is $\Bbb A^{35}$.

In Group 2, recall that the equation of $S_2$ be  $k_0(\bar x_4) = c_0(\bar x_4) = 0$ 
for smooth hypersurfaces $k_0(\bar x_4)=0$ and  $c_0(\bar x_4) = 0$.
 We only treat 3 of the 7 types in Group 2 in detail: 
 Let the
equations of
$F$ be:
$$
\align
 X_{4\vert 32,1\vert 20}:\hskip5cm&\\
k_0(\bar x_4) &=0 \\
 c_0(\bar x_4) z^2 + c_1(\bar x_4) yz + c_2(\bar x_4) y^2 &=0 \\
 X_{4\vert 32,1\vert 02}:\hskip5cm &\\
c_0(\bar x_4) &= 0\\
k_0(\bar x_4) z^2 + k_1(\bar x_4) yz + k_2(\bar x_4) y^2 &=0\\
  X_{5\vert 321,1\vert 002}:\hskip5cm&\\
k_0(\bar x_4)  + x_5l_1(\bar x_5)  &= 0\\
c_0(\bar x_4) + x_5k_1(\bar x_5) &=0\\
x_5 z^2+ l_2(\bar x_5) yz + l_3(\bar x_5) y^2 &= 0.
 \endalign
$$
respectively. The remaining 3 types essentially behave as the third
among these 3 listed types. Working locally, as in Group 1, we
obtain:
$$
\align
c_0(\bar x_4)  + tc_1(\bar x_4)  = k_0(\bar x_4) &= 0, \\
&\\
k_0(\bar x_4)  + tk_1(\bar x_4)  = c_0(\bar x_4) &= 0,   
\endalign
$$
and
$$
\align
k_0(\bar x_4)  -tl_2(\bar x_4,0) l_1(\bar x_4,0)  &=0\\ 
c_0(\bar x_4)
-tl_2(\bar x_4,0) k_1(\bar x_4,0) &=  0,
\endalign
$$

respectively.

 We recall that the equation of $S_2$ is  $k_0(\bar x_4) = c_0(\bar x_4) = 0$.
In analogy with case Group 1 we apply the linear map  $\phi_4$  from the parameter
space $\Bbb A^{15} \times \Bbb A^{35}$  of pairs of homogeneous quadric polynomials and cubic
polynomials in $5$ variables $\bar x_4$ to the tangent space $T_{R_4,\lbrack S_2\rbrack}$. The deformations above correspond to applying $\phi_4$
on subsets of the forms:

$$
\align
&\{(0, c_1(\bar x_4))\},  \\
&\{(k_1(\bar x_4),0)\},         \\
&\{\{(l_2(\bar x_4,0) l_1(\bar x_4,0), l_2(\bar x_4,0) k_1(\bar
x_4,0) \},
\endalign
$$
respectively.
 In the third case (and the 4 cases corresponding to the types not
listed), one easily sees that the linear span of
the subsets is all of $\Bbb A^{15}\times \Bbb A^{35}$. Put
$k_1=0$, and thereafter; $l_1=0$. This gives
$(l_2l_1, 0)$ and $(0, l_2k_2)$. The rest is as in Group 1.
Now consider the second of the first two cases where the linear
span is not all  $\Bbb A^{15}\times\Bbb A^{35}$.
In fact the desired deformations only
correspond to $\Bbb A^{15}\times{0}$. On the other hand, if
all of $\Bbb A^{15}\times{0}$ were contained in $\text{ker}(\eta_4)$, 
then for all  $k_1(\bar x_4)$, the

$$
k_0(\bar x_4)  + tk_1(\bar x_4)  = c_0(\bar x_4) = 0
$$
would give deformations such that $C_2$ deforms in the family. But then

$$
k_0(\bar x_4)    = c_0(\bar x_4) + tk_1(\bar x_4)l_1(\bar x_4) = 0
$$
would give the same thing, for all (the same) $k_1(\bar x_4)$  and all
linear forms $l_1(\bar x_4)$, and ${0}\times\Bbb A^{35}$ 
would be contained in $\text{ker}(\eta_4)$.
Hence the linear span, $\Bbb A^{15}\times\Bbb A^{35}$, 
of the union of $\Bbb A^{15}\times{0}$ and ${0}\times\Bbb A^{35}$,
would be contained in $\text{ker} (\eta_4)$, again a contradiction.

The last case  can  be handled by a modified version of the strategy working for
the other cases considered above.
Look at the maps
$$f\circ dP_4 \circ \phi_4\co \{0\}\times\Bbb A^{35}\longrightarrow
   T_{R_4,\lbrack S\rbrack}\longrightarrow T_{D_4,\lbrack S\rbrack}\longrightarrow \Bbb C$$
Here the first map is certainly not surjective, but the composite map
is. This is true since all smooth quadrics in $\Bbb P^4$ are projectively
equivalent.

 That is: Infinitesimal
deformations of $K3$-surfaces from $S$ of type:

$$       
           c_0(\bar x_4)  + tc_1(\bar x_4)  = k_0(\bar x_4) = 0,
$$
are as general as desired, when intrinsic properties, like containing a
smooth rational curve, are concerned. See also \cite\O{} for another formulation
of the argument.

 Finally we treat Group 3.
  We only treat 3 of the 5 types in Group 3 in detail:
Let the equations of $S_3$ be  $k_{01}(\bar x_5) = k_{02}(\bar x_5) = k_{03}(\bar x_5) = 0$.
Let the equations of $F$ be:
$$
\align
 X_{5\vert 222,1\vert 200}:\hskip5cm&\\
k_{01}(\bar x_5)z^2 + k_1(\bar x_5) yz + k_2(\bar x_5) y^2 &=0\\
 k_{02}(\bar x_5) =
k_{03}(\bar x_5) &= 0\\
 X_{5\vert 222,1\vert 110}:\hskip5cm&\\
k_{01}(\bar x_5)z + k_1(\bar x_5) y &=0\\
 k_{02}(\bar x_5)z + k_2(\bar x_5) y&=0\\
k_{03}(\bar x_5) &= 0\\
  X_{7\vert 22211,1\vert 00011}:\hskip5cm&\\
k_{01}(\bar x_5) + x_6l_1(\bar x_6)+ x_7l_2(\bar x_7) &=0\\
k_{02}(\bar x_5) + x_6l_3(\bar x_6)+ x_7l_4(\bar x_7) &=0\\
k_{03}(\bar x_5) + x_6l_5(\bar x_6) + x_7l_6(\bar x_7) &=0 \\
x_6z + l_7(\bar x_7)y &=0\\
x_7z + l_8(\bar x_7)y &= 0.
\endalign
$$
\flushpar
 The two remaining types are essentially treated as the last of the
 3 types listed. From a local viewpoint, these cases turn into:
$$
k_{01}(\bar x_5) + tk_1(\bar x_5) = k_{02}(\bar x_5) = k_{03}(\bar x_5) = 0,
$$
and
$$
k_{01}(\bar x_5) + tk_1(\bar x_5) = k_{02}(\bar x_5) + tk_2(\bar x_5)=
k_{03}(\bar x_5) = 0,
$$
and

$$
\align
k_{01}(\bar x_5) -t l_7(\bar x_5,0,0)l_1(\bar x_5,0) -
t l_8(\bar x_5,0,0)l_2(\bar x_5,0,0) &=0 \\
k_{02}(\bar x_5) -t l_7(\bar x_5,0,0)l_3(\bar x_5,0)-
t l_8(\bar x_5,0,0)l_4(\bar x_5,0,0) &=0\\
k_{03}(\bar x_5) -t l_7(\bar x_5,0,0)l_5(\bar x_5,0) -
t l_8(\bar x_5,0,0)l_6(\bar x_5,0,0) &= 0.
\endalign
$$

We recall that the equation of $S_3$ is
$$k_{01}(\bar x_5) = k_{02}(\bar x_5) = k_{03}(\bar x_5) = 0.$$
In analogy with Group 1 we apply the linear map  $\phi_5$  from the
parameter space $(\Bbb A^{21})^3$  of triples of homogeneous quadric polynomials in $6$ variables $\bar x_5$ to the tangent space $T_{R_5\lbrack S_3\rbrack}$.
The deformations above correspond to applying $\phi_5$ on the
subsets: 
 
$$
\align
&\{(k_1(\bar x_5) ,0,0\}, \{k_1(\bar x_5) , k_2(\bar x_5),0\}\\
&\{(-l_7(\bar x_5,0,0)l_1(\bar x_5,0) - \\ &l_8(\bar
x_5,0,0)l_2(\bar x_5,0,0),-l_7(\bar x_5,0,0)l_3(\bar x_5,0)-\\
&l_8(\bar x_5,0,0)l_4(\bar x_5,0,0),  -l_7(\bar x_5,0,0)l_5(\bar
x_5,0) -\\ &l_8(\bar x_5,0,0)l_6(\bar x_5,0,0))\}.
\endalign
$$

The linear span of the third subset is $(\Bbb A^{21})^3$.
 In the first case the linear
span is $\Bbb A^{21}\times{0}\times{0}$.  
On the other hand, if both $\Bbb A^{21}\times{0}\times{0}$, and $ {0} \times
\Bbb A^{21}\times{0}$, and ${0}\times{0}\times\Bbb A^{21}$, 
all were contained in $\text{ker}(\eta_5)$, then the
linear span of the union of these $3$ sets, i.e.
$\Bbb A^{21}\times\Bbb A^{21} \times\Bbb A^{21}$, would be,
again a contradiction. Hence one of these $3$ sets is not contained in the
union, and by interchanging the roles of $k_{01}(\bar x_5)$, 
$k_{02}(\bar x_5)$, $k_{03}(\bar x_5)$ if necessary we obtain an
element outside $\text{ker}(\eta_5)$. 
This argument certainly applies to give a proof for the second case 
also, since a desired deformations in the first case also is included in 
the second case.

 To complete the proof of Theorem 4.3 it is enough to prove that in each of
the $16$ cases listed above, a general F as constructed, is smooth. To do so
we will use the following version (extension) of Bertini's theorem, as
given in [10].

\proclaim{Lemma 4.4 ([10])}
On an arbitrary ambient variety V, if a linear system has no fixed
components, then the general member has no singular points outside of the
union of the base locus of the system and the singular locus of the
ambient variety.
\endproclaim

Repeated usage of this result  in
each of the $16$ cases will solve the problem. We will be able to put
ourselves in a position where the ambient varieties always are smooth.

As an example take the case $X_{4|32, 1|20}$:
  $$c_0(\bar x_4) z^2+ c_1(\bar x_4)yz + c_2(\bar x_4)y^2 =  k_0(\bar x_4)  = 0.$$
Here the fixed smooth hypersurfaces $Z(k_0(\bar x_4) )$ and $Z(c_0(\bar x_4) )$ cut out
the smooth surface $S_2$ in  $\Bbb P^4$.  On the smooth ambient variety
$V_1=Z(k_0(\bar x_4) )$ in $\Bbb P^4 \times \Bbb P^1$ we study the linear system
$$\frak L_1 = \{ c_0(\bar x_4)z^2+ c_1(\bar x_4)yz + c_2(\bar x_4)y^2 | \ c_1(\bar x_4) 
\text{and $c_2(\bar x_4)$ arbitrary cubics} \}$$
This linear system has no fixed components, and the base locus of 
$\frak L_1$ is
\flushpar
$Z(c_0(\bar x_4), k_0(\bar x_4), y)$ in $\Bbb P^4 \times \Bbb P^1$, that is $S_2$. But since $S_2$ is the 
scheme-theoretical complete intersection in $V_1$ of any element $l$  of $\frak L_1$ with the
hypersurfaces $c_0(\bar x_4)$ and  $y$, this element $l$ cannot be singular at a point in
$S_2$ without $S_2$ being so. Hence all elements of $\frak L_1$ are smooth at all points
of the base locus of $\frak L_1$. Since a general element is smooth outside the
base locus, by Lemma 4.4, a general element of  $ \frak L_1$ is smooth. Hence a
general $F$ of type $X_{4|32, 1|20}$ is smooth.
This was a relatively easy type to handle. Other easy types, which can be
handled in a similar way, are:
$ X_{3|4, 1|1},  X_{4|32, 1|02} , X_{5|222, 1|200}$ . In these cases it suffices to apply
Lemma 4.4. once.
As another (and slightly more involved) example take the case
$X_{5|411, 1|011}$:
$$ q_0(\bar x_3) + x_4c_1(\bar x_4) + x_5c_2(\bar x_5) = x_4z +l_1(\bar x_5)y = x_5z + l_2(\bar x_5)y =
 0.$$
First, let the ambient variety $V$ be all of $\Bbb P^5 \times \Bbb P^1$, which is smooth. The
linear systems
$$
 \frak L_1 =\{ q_0(\bar x_3) + x_4c_1(\bar x_4) + x_5c_2(\bar x_5) | \ c_1(\bar x_4) 
\text{and $c_2(\bar x_5)$ arbitrary cubics}\}
$$
and    
$$
\frak L_0 =\{ q_0(\bar x_3) + x_4c_1(\bar x_4) | \  c_1(\bar x_4)  \text{ arbitrary cubic}\}
$$
have no fixed components, and the base loci are $q_0(\bar x_3) = x_4 = x_5 = 0$ and
$q_0(\bar x_3) = x_4 = 0$, respectively. These loci are isomorphic to  
$S_1 \times\Bbb P^1$, and
a $x_5$-cone over $S_1 \times  \Bbb P^1$, respectively, for the smooth quartic 
$S_1$ in $\Bbb P^3$ from Proposition 4.1.  
Hence the first base locus, say  $B_1$, is smooth, and 
the second, say $B_0$, is singular only along $(0,..,0,1) \times  \Bbb P^1$. 
But $B_1$ is also the scheme-theoretical intersection of any member of 
$ \frak L_1$ and $Z( x_4, x_5)$. 
Since $B_1$ is smooth of dimension $3$, no member of $\frak L_1$ can then be 
singular at any point of $B_1$, and the general members of $ \frak L_1$ 
is then smooth everywhere, in virtue of Lemma 4.4. By an analogous argument
a general member of $ \frak L_0$ can only be singular along 
$(0,..,0,1) \times  \Bbb P^1$. Hence there exists an element 
$q_0(\bar x_3) + x_4c_1(\bar x_4) + x_5c_2(\bar
x_5)$ of 
$ \frak L_1$ , such
that  $Z(q_0(\bar x_3) + x_4c_1(\bar x_4) + x_5c_2(\bar x_5))$ is smooth,
and $Z(q_0(\bar x_3) + x_4c_1(\bar x_4))$ is singular only along
$(0,..,0,1) \times  \Bbb P^1$. Fix $Z(q_0(\bar x_3) + x_4c_1(\bar x_4) + x_5c_2(\bar x_5))$  as the new ambient
variety $V_2$.  We denote  $Z(q_0(\bar x_3) + x_4c_1(\bar x_4), x_5, y )$ by $G$. 
This is smooth, since $Z(q_0(\bar x_3) + x_4c_1(\bar x_4) )$  is smooth for all
points with $x_5 = 0$. From now on $c_1(\bar x_4)$ and $c_2(\bar x_5)$ are fixed.
Next, consider the linear system on  $V_2$:
 $$
\frak L_2 = \{ x_5z + l_2(\bar x_5)y  | \ l_2(\bar x_5) \text{ arbitrary} \}.$$
This system has no fixed component, and the base locus is $V_2\cap Z( x_5,y)$.
This is the smooth quartic hypersurface $G$  in  $\Bbb P^4$.  It is  clear that  all
members of $ \frak L_2$ are smooth on $G$ (On each point of $G$, $y=0$, 
and hence $z$ can be
taken to be one, and hence  we may express properly dehomogenized
versions of $q_0(\bar x_3) + x_4c_1(\bar x_4) + x_5c_2(\bar x_5)$ and $x_5z + l_2(\bar x_5)y$ in terms of
local parameters around any point $p$ of $G$, such that these expressions are
independent modulo the square of the maximal ideal of $p$).
Hence a general element of $ \frak L_2$ will give a smooth submanifold of $V_2$.
Fix one such element, that is, fix $l_2(\bar x_5)$. This gives a smooth submanifold
$V_3$ given as 
$$q_0(\bar x_3) + x_4c_1(\bar x_4) + x_5c_2(\bar x_5) =  x_5z + l_2(\bar x_5)y = 0.$$
Next, consider the linear system on  $V_3$:
$$
\frak L_3 = \{ x_4z + l_1(\bar x_5)y  | \  l_1(\bar x_5) \text{arbitrary} \}.
$$
This linear system has no fixed component, and the basis locus is given as

\flushpar
$Z(q_0(\bar x_3), x_4 , x_5 , y)$. This is $S_1$, a smooth surface. But 
$S_1 = V(y| V_3, l| V_3)$,
scheme-theoretically, for each $l$ in $\frak L_3 $, so no such $l$ in 
$\frak L_3$ can be singular
on $S_1$. Hence a general member of $\frak L_3$ is smooth.
This gives that a general $F$ as constructed, of the type $X_{5|411, 1|011}$ is a
smooth 3-fold. Other cases that can be treated in a very similar way are:
$X_{5|321, 1|002} , X_{5|321, 1|011}, X_{5|321, 1|101}$.

Intermediate cases (where we need to use Lemma 4.4 twice) are:
$$X_{4|41, 1|02},  X_{4|41, 1|11},  X_{4|32, 1|11},  X_{5|222, 1|200}.$$

The cases $X_{6|2221, 1|0002}, X_{6|2221, 1|0011}, X_{7|22211, 1|00011}$ are
similar to each other, but no essentially new technique is required.
As an example we treat the case $X_{7|22211, 1|00011}$ :
 $$
\align
k_{01}(\bar x_5) + x_6l_1(\bar x_6) + x_7l_2(\bar x_7) =&0\\
  k_{02}(\bar x_5) + x_6l_3(\bar x_6) + x_7l_4(\bar x_7)=&0\\
 k_{03}(\bar x_5) + x_6l_5(\bar x_6) + x_7l_6(\bar x_7) =&0\\
 x_6z + l_7(\bar x_7)y =&0\\
 x_7z + l_8(\bar x_7)y=&0.
\endalign
$$
Using Lemma 4.4 several times as in the first part of the proof in the case
of type $X_{5|411, 1|011}$ , we obtain that for a generic fixed choice of
$l_1(\bar x_6)$, $l_2(\bar x_7)$, $l_3(\bar x_6)$ , $l_4(\bar x_7)$ ,$l_5(\bar x_6)$ , $l_6(\bar x_7)$ both the $(2,2,2)$ four-fold
$$
\align
k_{01}(\bar x_5) + x_6l_1(\bar x_6) + x_7l_2(\bar x_7) &=0\\
k_{02}(\bar x_5) + x_6l_3(\bar x_6) + x_7l_4(\bar x_7)&=0\\
k_{03}(\bar x_5) + x_6l_5(\bar x_6) + x_7l_6(\bar x_7) &= 0
\endalign
$$
 in $\Bbb P^7$ and the $(2,2,2)$ threefold
$$
\align
k_{01}(\bar x_5) + x_6l_1(\bar x_6) =&0\\
  k_{02}(\bar x_5) + x_6l_3(\bar x_6)  =&0\\
   k_{03}(\bar x_5) + x_6l_5(\bar x_6)  = &0
\endalign
$$

in $\Bbb P^6$
are smooth, and so are their respective counterparts in 
$\Bbb P^7\times \Bbb P^1$, given by
the same equations. We stick to such a fixed generic choice of the $6$ linear
functions. Let the new ambient space $V_1$ be defined by:
$$
\align
k_{01}(\bar x_5) + x_6l_1(\bar x_6) + x_7l_2(\bar x_7) =&0\\
  k_{02}(\bar x_5) + x_6l_3(\bar x_6) + x_7l_4(\bar x_7)=&0\\
k_{03}(\bar x_5) + x_6l_5(\bar x_6) + x_7l_6(\bar x_7)=&0
\endalign
$$
 in $\Bbb P^7 \times P^1$.
Let $ \frak L_1$ be the linear system 
$$\{x_7z + l_8(\bar x_7)y | \ l_8(\bar x_7) \text{ is arbitrary}\}$$ 
on $V_1$.
The base locus of $\frak L_1$ is $V_1$ is $Z(\bar x_7, y)$, which is the smooth 3-fold 
in $\Bbb P^6$
referred to. This is the scheme-theoretical intersection of $Z(\bar x_7, y)$ on $ V_1$ and
every element of $ \frak L_1$; hence every element of $ \frak L_1$ is smooth along this base
locus, and concequently a generic element is smooth everywhere. Choose
one such fixed element, that is fix  $l_8(\bar x_7)$, and let the new ambient variety
$V_2$, be cut out by this element on $V_1$.
At last study the linear system 
$$ \frak L_2 = \{x_6z + l_6(\bar x_7)y | \  l_6(\bar x_7) \text{ is arbitrary}\}$$
on $V_2$. The base locus is $$Z(y, k_{01}(\bar x_5), k_{02}(\bar x_5), k_{03}(\bar x_5)) = S_3,$$ which is
also the scheme-theoretical  complete intersection of $y$ and any element
of on $V_2$. Hence any such section is smooth on $S_3$, and a generic section is
smooth everywhere. Hence a generic 3-fold $F$ constructed of type
$X_{7|22211, 1|00011}$ is smooth.
   \enddemo

 \subheading{5.  Determinantal contractions}
\vskip.3cm
In this section we will prove Corollary 1.4. First we will
 give a way to relate a certain 
CICY in one multiprojective space to a CICY in another
multiprojective space. The method of determinantal contractions
is introduced in \cite{\CandDaLuSc}. The construction
is used in \cite{\CandGrHu} to show that the moduli space of
CICY varieties is connected.  
We use the same notation as introduced in the paragraph 1.

Let $X$ be a CICY of type:

$$
\Draw
\NewNode(\RectDiNoden,\MoveToRect){   \Units(1pt,1pt)
   \GetNodeSize   \SetMinNodeSize
   \Move(-\Val\Va,-\Val\Vb)\Move(-1.0,-1.0)   \Va*2;   \Vb*2;
   \Line(0,\Val\Vb) \Line(0,2.0) \Line(\Val\Va,0)\Line(2.0,0)
   \Line(0,-\Val\Vb)\Line(0,-2.0) \Line(-\Val\Va,0)\Line(-2.0,0)
   \Move(18.0,0) \Line(0,\Val\Vb) \Line(0,2.0)
   \Move(2.0,0) \Line(0,-\Val\Vb)\Line(0,-2.0)}
\RectDiNoden(a)(--$\matrix\format\c&\quad\c&\quad\c&\quad\c&\quad\r\\ N&M&a_0&\cdots&a_n\\ n&0&1&\cdots&1\endmatrix$--) 
\EndDraw
\tag{5.1}$$
 \       

\vskip1cm
 Here we have introduced a block notation. The $N$ is a column consisting
 of the dimensions of the first $k-1$ projective spaces. Likewise $M$ is
 a block matrix where the rows are the multidegrees of the first $m-n-1$ polynomials
 in the $k-1$ first projective spaces. The $a_i$'s are $(k-1)$-dimensional vectors.
 The zero in the last row is a short
 way of writing $m-n-1$ zeros.
In the following we denote the last $n+1$ polynomials by $q_0,\cdots, q_n$
We can write the $q_i$'s as
$$
\matrix
   q_0 = q_{00}x_0+\cdots + q_{0n}x_n\\
  \vdots \\
 q_n = q_{n0}x_0+\cdots + q_{nn}x_n
\endmatrix
$$
 where $x_0,\dots,x_n$ are the projective coordinates in the projective
 space singled out in the last row of (5.1).      We use  $P$ for the product of the remaining
 $k-1$ projective spaces in order to keep the notation simple.
We can write this on matrix form
$$
\pmatrix
  q_0\\
  \cdot\\
 \cdot\\
 \cdot\\
 q_n
\endpmatrix = 
A\pmatrix
  x_0\\
 \vdots\\
 x_n
\endpmatrix              
$$
where
$$
A=
\pmatrix
  q_{00}&\hdots&q_{0n}\\
  \vdots&      & \vdots \\
  q_{n0}&\hdots&q_{nn}\\
\endpmatrix      \tag{5.2}
$$
That some point $(a,b) \in P\times \Bbb P^n$ is in the common zero set of the $q_i$'s,
is equivalent to $a$ being in the zero set of $\text{det}(A)$, since not all the
$x_i$'s can vanish. In other words, the projection of $F$ on the factor
$P$ gives us a new variety $F_s$ of type:

$$
\Draw
\NewNode(\RectDiNodem,\MoveToRect){   \Units(1pt,1pt)
   \GetNodeSize   \SetMinNodeSize
   \Move(-\Val\Va,-\Val\Vb)\Move(-1.0,-1.0)   \Va*2;   \Vb*2;
   \Line(0,\Val\Vb) \Line(0,2.0) \Line(\Val\Va,0)\Line(2.0,0)
   \Line(0,-\Val\Vb)\Line(0,-2.0) \Line(-\Val\Va,0)\Line(-2.0,0)
   \Move(16.0,0) \Line(0,\Val\Vb) \Line(0,2.0)
   \Move(2.0,0) \Line(0,-\Val\Vb)\Line(0,-2.0)}
\RectDiNodem(a)(--$\matrix\format\c&\quad\c&\quad\r\\ N & M     &\text{mult}\vert A\vert\endmatrix$--) 
\EndDraw
$$
where  $\text{mult}\vert A\vert$ denotes the 
multidegree of the determinant of $A$.
This is a singular variety in general. 

\flushpar\proclaim{Proposition 5.1}
Let $F$ be a CICY of dimension $m$ 
in a multiprojective space. Suppose that there exists
a determinantal contraction
$f\co F\longrightarrow F_s$.
Let $$F^r_s=\{ p \vert\quad \text{\ dim}f^{-1}(p)\geq r-1\}.$$ Then
$F^r_s$ is the zero set of the $(n-r+2)$- minors of $A$, and
  the fiber over each point $q\in F^r_s-F^{r+1}_s$
is isomorphic to $\Bbb P^{r-1}$.
\endproclaim

\flushpar\demo{Proof}
 We     use the notation introduced above, i.e., 
 $F$ is of type

 $$
\Draw

\NewNode(\RectDiNodek,\MoveToRect){   \Units(1pt,1pt)
   \GetNodeSize   \SetMinNodeSize
   \Move(-\Val\Va,-\Val\Vb)\Move(-1.0,-1.0)   \Va*2;   \Vb*2;
   \Line(0,\Val\Vb) \Line(0,2.0) \Line(\Val\Va,0)\Line(2.0,0)
   \Line(0,-\Val\Vb)\Line(0,-2.0) \Line(-\Val\Va,0)\Line(-2.0,0)
   \Move(18.0,0) \Line(0,\Val\Vb) \Line(0,2.0)
   \Move(2.0,0) \Line(0,-\Val\Vb)\Line(0,-2.0)}
\RectDiNodek(a)(--$\matrix\format\c&\quad\c&\quad\c&\quad\c&\quad\r\\ N&M&a_0&\cdots&a_n\\ n&0&1&\cdots&1\endmatrix$--) 
\EndDraw
$$

\flushpar 
and $F_s$
is
of type
 $$
\Draw
\NewNode(\RectDiNodej,\MoveToRect){   \Units(1pt,1pt)
   \GetNodeSize   \SetMinNodeSize
   \Move(-\Val\Va,-\Val\Vb)\Move(-1.0,-1.0)   \Va*2;   \Vb*2;
   \Line(0,\Val\Vb) \Line(0,2.0) \Line(\Val\Va,0)\Line(2.0,0)
   \Line(0,-\Val\Vb)\Line(0,-2.0) \Line(-\Val\Va,0)\Line(-2.0,0)
   \Move(16.0,0) \Line(0,\Val\Vb) \Line(0,2.0)
   \Move(2.0,0) \Line(0,-\Val\Vb)\Line(0,-2.0)}
\RectDiNodej(a)(--$\matrix\format\c&\quad\c&\quad\r\\ N & M     &\text{mult}\vert A\vert\endmatrix$--) 
\EndDraw
$$
and is  defined by the zero set of the first $m-n-1$ polynomials defining
$F$ and the determinant of the following matrix:

 $$
A=
\pmatrix
  q_{00}&\hdots&q_{0n}\\
  \vdots&      & \vdots \\
  q_{n0}&\hdots&q_{nn}\\
\endpmatrix      \tag{5.3}
$$

 We
have the following commutative diagram:

$$
\CD
X             @>>>         P\times \Bbb P^n\\
@V\bar{pr}_1VV              @Vpr_1VV\\
X_s           @>>>         P
\endCD
$$
where $pr_1$ is the projection on the first factor and $\bar{pr}_1$
its restriction.

This means that two points $p=q\times s$ and
$p^\prime = q^\prime\times s^\prime$ are mapped
to the same point in $F_s$ if and only if
$q=q^\prime$.
Let $f\in \Cal P_1\times\dots\times\Cal P_k$
 be polynomials defining
$F$, and let $q$ be a point in $F_s$.
Then  $q\times v$ is in the fiber over $q$
if and only if $Av=0$ where $A$ is the matrix
associated to $f$ as above. In other words the
 points in
the fiber are naturally identified with the
kernel of the map $A\co \Bbb A^{n+1} \rightarrow \Bbb A^{n+1}$, affinely.
Hence, if $q\in F^r_s-F^{r+1}_s$ then the fiber has affine
dimension $r$. This gives $pr_1^{-1}(q)\cong \Bbb P^{r-1}$.\qed

 \enddemo

\remark{Remark 5.2}
Let $Y$ be a threefold with an isolated
singularity $p$. Assume that there exists 
a nonsingular variety $Y^\prime$ and a morphism
$$Y^\prime\overset{\pi}\to\rightarrow Y,$$
\flushpar
such that $\pi^{-1}(p)\cong \pen$, and $Y^\prime-\pi^{-1}(p)\cong Y-\{p\}$.
Then we say that $Y^\prime\overset{\pi}\to\rightarrow Y$
 is a small resolution of the singularity $p$ (\cite{\Atiyah}).
 Furthermore, if $Y$ contains several isolated singularities,
 and we can find a $Y^\prime$ such that $Y^\prime$ is, locally
 around each singularity on $Y$, a small resolution, then
 we say that $Y^\prime\overset{\pi}\to\rightarrow Y$ is
 a small resolution of $Y$.
\endremark
\vskip.2cm
\flushpar\proclaim{Lemma 5.3}
Let $V\subseteq \Bbb P^n$ be a smooth subvariety, and let
$S=Z(Q_0,L_0)\cap V$ be a smooth codimension $2$ subvariety
of $V$, where $\text{deg}(Q_0)=s$, and $\text{deg}(L_0)=t$.
Let $\frak L$ be the linear system $\{Q_0L_1-Q_1L_0\}$ on $V$,
where $L_1$ varies through all $t$-ics and $Q_1$ varies
through all $s$-ics. Then  for a general element $l$ of $\frak L$,
$\text{Sing }(l)=Z(Q_0,Q_1,L_0,L_1)$. In general $Z(Q_0, Q_1,L_0,L_1)\cap V $
will be a smooth manifold of degree $s^2t^2\text{deg } V$.
\endproclaim
\flushpar\demo{Proof}
By
 Lemma 4.4 a general $V$ has no singular points outside the base locus of
$\frak L$. This base locus is $S$. For each point $p$ in $V$, let $W$
be the subspace of $\Bbb A^{n+1}$ spanned by the gradients of a set
of generators of $I(V)$. For $p$ in $V$ and $l$ in $\frak L$ we have
that $p$ belongs to $\text{Sing }(l)$ iff

$$\nabla l(p)=(\frac{\partial l}{\partial x_0}(p),
\dots \frac{\partial l}{\partial x_n}(p))\in W.\tag{5.4}$$
Set $l=Q_0L_1-Q_1L_0$. Then
$$
\frac{\partial l}{\partial x_i}(p)=
Q_0(p)\frac{\partial L_1}{\partial x_i}(p)
+L_1(p)\frac{\partial Q_0}{\partial x_i}(p)
-Q_1(p)\frac{\partial L_0}{\partial x_i}(p)
-L_0(p)\frac{\partial Q_1}{\partial x_i}(p),$$ 
for $\ i= 0,\dots, n$

For points in the base locus $Z(q_0,L_0)\subseteq V$ this
reduces to:
 $$
 \frac{\partial l}{\partial x_i}(p)=
L_1(p)\frac{\partial Q_0}{\partial x_i}(p)
-Q_1(p)\frac{\partial L_0}{\partial x_i}(p), \ i= 0,\dots, n
$$
Hence (5.4) reduces to
$\lbrack L_1(p)\ Q_1(p) \rbrack M \in W$, where the
rows  of $M$ is a matrix with rows  that are the gradients of
$Q_0$ and $-L_0$ at $p$, respectively. The subvariety $S$ is smooth of 
codimension $2$ at $p$, so the rows are independent. Hence,  for a general element $l$ of $\frak L$,
$\text{Sing }(l)=Z(Q_0,Q_1,L_0,L_1)$ 
The last assertion of the lemma is immediate.
\enddemo

\subheading{ Application to complete intersection Calabi--Yau manifolds in projective
space}  

In the previous section we constructed isolated rational curves of
multidegree $(d,0)$ for every integer $d$, and for every $\Bbb P^1$-surjective 
CICY-type in biprojective space. Some of these
CICY-types allow a determinantal contraction to  mildly singular
Calabi-Yau varieties in a projective space of dimension $4,5,6$ or $7$.
We will use this to show existence of isolated rational curves of every
degree on generic CICYs in these projective spaces. There are
$5$ types $X_{4\vert5}$, $X_{5\vert42}$, $X_{5\vert33}$, $X_{6\vert332}$,$X_{7\vert2222}$.
Consider first  a CICY $F$ of  the type $X_{4\vert41,1\vert 11}$ constructed
in the
previous section.

$$
{\pmatrix\format\r&\quad\r\\
q_0(\bar x_3)+x_4c_1(\bar x_4) & q_1(\bar x_4)\\
x_4 & l_1(\bar x_4)
\endpmatrix}
{\pmatrix
z\\
y
\endpmatrix}
=0
$$
 In the following denote the $2\times 2$ matrix by $A$. 
The threefold $X_{4\vert41,1\vert 11}$ has a determinantal contraction
to the quintic $ F_s=Z(\text{det}(A))\subseteq \Bbb P^4$. Explicitly
$F_s$  is the zero-set of the polynomial:
$(q(\bar x_3)+x_4c_1(\bar x_4))l_1(\bar x_4)- x_4 q_1(\bar x_4)$.
By letting $Q_0=q(\bar x_3)+x_4c_1(\bar x_4)$ and $L_0=x_4$ and 
$V=\Bbb P^3\times \Bbb P^1$ and using Lemma 5.3 we get that $F_s$ has exactly the $16$
 singular points given by the vanishing of the four polynomials
$q_0(\bar x_3), x_4, q_1(\bar x_4), l_1(\bar x_4)$.

 Let $ \pi\co F \longrightarrow F_s$ be the contraction map. As
noted before the fiber over each singular point is a $\Bbb P^1$.
By construction $F$ contains a rational curve $C_1$  of bidegree $(d,0)$.
Furthermore, if $C_1$ does not intersect any of the $\Bbb P^1$ contracted
by $\pi$, then  the image of $C_1$ by $\pi$ is an isolated curve on
$F_s$. 
This is possible to achieve. First observe that the intersection
between $S_1$ in $F$ and a $\Bbb P^1$ that is contracted is transversal,
i.e. gives a point in $S$. Furthermore, varying  $q_1(\bar x_4)$
and $l_1(\bar x_4)$ moves the $16$ intersection points freely in $S_1$.
But   $q_1(\bar x_4)$
and $l_1(\bar x_4)$ were chosen to be generic (by construction of $F$
in the previous paragraph). Hence, $C_1$ does not intersect any of
the $\Bbb P^1$'s contracted by $\pi$. Denote the image of $C_1$ by
$D$. By construction 
$N_{D/F_s}\cong N_{C_1/F}\cong \oto{-1}{-1}$. Then Corollary
1.4 follows from Theorem 3.4 in case of quintic threefolds in $\Bbb P^4$.

The other four CICY-types:

$X_{5\vert 33}:$ 

Here we consider the following $F$ of type $X_{5\vert321,1\vert 011}$
with matrix $A$ given by:

$$
A=
\pmatrix
k_0(\bar x_4)+ x_5 l_1(\bar x_5) & k_1(\bar x_5)\\
x_5     &  l_2(\bar x_5)
\endpmatrix
$$

$X_{5\vert 42}:$ 

Here we consider the following $F$ of type $X_{5\vert321,1\vert 101}$
with matrix $A$ given by:

$$
A=
\pmatrix
c_0(\bar x_4)+ x_5 k_1(\bar x_5) & c_1(\bar x_5)\\
x_5     &  l_2(\bar x_5)
\endpmatrix
$$

$X_{6\vert 322}:$ 

Here we consider the following $F$ of type $X_{6\vert2221,1\vert 1001}$
with matrix $A$ given by:

$$
A=
\pmatrix
k_{01}(\bar x_5)+ x_6 l_1(\bar x_6) & k_1(\bar x_6)\\
x_6     &  l_4(\bar x_5)
\endpmatrix
$$

$X_{7\vert 2222}:$ 

Here we consider the following $F$ of type $X_{7\vert22211,1\vert 00011}$
with matrix $A$ given by:

$$
A=
\pmatrix
 x_6 & l_7(\bar x_7)\\
x_7    &  l_8(\bar x_7)
\endpmatrix
$$

As another example, let us study the last case in detail. 
We will study the contraction  $\pi; F\longrightarrow F_s$.
Let $V$ be defined by:

$$
\align
X_{7\vert 22211,1\vert 00011}:\hskip5cm&\\
k_{01}(\bar x_5) + x_6l_1(\bar x_6)+ x_7l_2(\bar x_7) &=0\\
k_{02}(\bar x_5) + x_6l_3(\bar x_6)+ x_7l_4(\bar x_7) &=0\\
k_{03}(\bar x_5) + x_6l_5(\bar x_6) + x_7l_6(\bar x_7) &=0 
\endalign
$$
Take $Q_0=x_6$ and $L_0=x_7$. We see that $Q_0$ and $L_0$
cut out 
$$S_3=Z(k_{01}(\bar x_5), k_{02}(\bar x_5), k_{03}(\bar x_5),x_6,x_7)$$
inside $V$, and that $S_3$ is smooth of codimension 2 in $V$.
Hence Lemma 5.3 can be used. 
The threefold $F_s= V\cap Z(x_6l_8(\bar x_7)- x_7l_7(\bar x_7)$ will only be
singular at  the $8$ points
 given by
the vanshing of the $7$ polynomials
$$k_{01}(\bar x_5),k_{02}(\bar x_5), k_{03}(\bar x_5),l_7(\bar x_7),
l_8(\bar x_7), x_6, x_7.$$
 for general $l_7(\bar x_5)$ and
 $l_8(\bar x_5)$.

Furthermore, for general $l_7(\bar x_7)$ and $l_8(\bar x_7)$ the isolated
 bidegree $(d,0)$ curve $C_3$ in $S_3$ does not intersect any of
 the $\Bbb P^1$ contracted by $\pi$. Hence, the image of $C_3$ by
$\pi$ is a curve $D$ of degree $d$.
By construction $N_{D/F_s}\cong N_{C_1/F}\cong \oto{-1}{-1}$. Then Corollary
1.4 follows from Theorem 3.4 in this case.

 \remark{Remark 5.4}
In order to prove Corollary 1.4 we used five determinantal contractions.
There are $9$ of the $16$  CICY-types of Lemma 3.2 that have determinantal contractions.
The four ``unused" CICY-types with such contractions are $X_{5|411,1|011}$, $X_{5|222,1|110}$,
$X_{4|32,1|11}$, $X_{6|3211,1|0011}$. The first two can be used to give additional proofs
of Corollary 1.4. in case of $X_{5|42}$. These two additional proofs
use the rigid curves on $S_2$ and $S_3$ respectively. Hence 
the $X_{5|42}$ case can be proved using any of the three surfaces
$S_1,S_2,S_3$. The contraction $X_{4|32,1|11}$ can be used to give an
additional proof of Corollary 1.4 in the $X_{4|5}$ case (using $S_2$),
while the contraction $X_{6|3211,1|0011}$ can be used in the $X_{6|322}$-case.
For the remaining CICY-types we have no choice of determinantal contraction.
\endremark

\centerline{References}

\references 

\ref{Atiyah} M. F. Atiyah, {\it  On analytic surfaces with double points,}
                     Proc. Roy. Soc., {\bf A247}(1958), 237--244. 

\ref{CandDaLuSc}  P. Candelas, A. M. Dale, C. A. L\"utken, R. Schimmrigk,
 {\it   Complete intersection Calabi--Yau manifolds,}                
            Nucl. Phys.  {\bf B298} (1988),  493--525.

\ref{CandGrHu} P. Candelas, P.S. Green, T. H\"ubsch,
              {\it  Rolling among Calabi--Yau vacua, }
              Nucl. Phys. {\bf B330}  (1990), 49--102.

\ref{CandLuSc}    P. Candelas, C. A. L\"utken, R. Schimmrigk,     
             {\it  Complete intersection Calabi--Yau manifolds II
                    Three generation manifolds,}

\ref{Cl}     H. Clemens,
             {\it  Homological equivalence, modulo algebraic
              equivalence, is not finitely generated,}
              Publ. Math. I.H.E.S. {\bf 58}   (1983),19--38.

\ref{GrH}    P. S. Green, T. H\"ubsch,
             {\it  Possible phase transitions among Calabi--Yau compactifications, }
              Phys. Rev. Lett. {\bf 61} (1988),  1163--1166.

\ref{GrHu} P. S. Green, T. H\"ubsch,
             {\it  Calabi--Yau Manifolds as Complete Intersections
                    in Products of Complex Projective Spaces,}
             Comm. Math. Phys. {\bf 109}  (1987),  99--108.

\ref{GrHa}   P. Griffiths, J. Harris, {\it Principles of Algebraic
             Geometry,} Wiley, 1978.

\ref{JohKlei}   T. Johnsen, S. L. Kleiman,
                {\it Rational curves of degree at most 9 on a general quintic threefold,        }
               Comm. Alg.  {\bf 24(8)}  (1996), 2721--2753.

\ref{Kleim}   S. L. Kleiman,
                {\it Bertini and his two fundamental theorems,}
                 Rendic. del circ. Matem. di Palermo, serie II,
                Supplemento (1987).

\ref{Katz}     S. Katz,
               {\it  On finiteness of rational curves on quintic
               threefolds,}
              Compositio Math. {\bf 60}   (1986), 151--162.

\ref{Kley}   H. P. Kley,
             {\it Rigid curves in quintic threefolds,}
                         Preprint (1996)

\ref{Kn}     A. L. Knutsen,
             {\it On Degrees and Genera of Smooth Curves on Projective K3
             surfaces,}
             Math.AG/9805140, (1998)                       

\ref{Ko}     K. Kodaira,
             {\it  On the structure of compact complex analytic surfaces,
              I,} American J. Math. {\bf 86}   (1964), 751--98.

\ref{Mo}     S. Mori,
               {\it  On degrees and genera of curves on smooth
               quartic surfaces in $\Bbb P^3$,}
              Nagoya Math. J. {\bf 96}   (1984), 127--32.

\ref{Nijsse} P. G. J. Nijsse, 
              {\it Clemens' conjecture for octic and nonic curves,}
               Indag. Math. (N.S) 6, (1995) no. 2, 213-221.

\ref{O}     K. Oguiso,
               {\it  Two remarks on Calabi-Yau Moishezon
               threefolds,} J. fur die reine und angew. Math. {\bf 452}
            (1994),153--62.

\ref{So}    D. E. Sommervoll, {\it Rational curves on Complete intersection Calabi--Yau threefolds,} Dr.sci.-thesis, University of
               Oslo, 1997.

\endreferences

\end{document}